  \documentclass[11pt]{amsart}
   \usepackage[breaklinks,colorlinks=true,allcolors=orange]{hyperref}

  \usepackage[margin=2.8cm]{geometry}
  \usepackage{latexsym}
  \usepackage{comment}
  \usepackage{amsmath,amsthm,amsfonts,amssymb,graphicx,epsfig,latexsym,float}
  \usepackage{tikz-cd} 
  \definecolor{MyGray}{RGB}{128, 128, 128}

  \usepackage{epsf,enumerate}
  \usepackage{overpic}  
  \usepackage[normalem]{ulem}

  \makeatletter
  \def\paragraph{\@startsection{paragraph}{4}%
    \z@\z@{-\fontdimen2\font}%
    {\normalfont\bfseries}}
  \makeatother

  \def\int{\operatorname{int}}
  \newcommand{\diagclosure}[1]{\langle\!\langle #1\rangle\!\rangle}

  \newtheorem{thm}{Theorem}[section]
  \newtheorem{lemma}[thm]{Lemma}
  
  \newtheorem{cor}[thm]{Corollary}
  \newtheorem{prop}[thm]{Proposition}
  {}
  {}

  \newtheorem{thmintro}{Theorem}

  \newtheorem{propintro}[thmintro]{Proposition}
  
  \theoremstyle{remark}
  \newtheorem{example}[thm]{Example}

  \newtheorem{remark}[thm]{Remark}
  
  \newtheorem*{definition*}{Definition}
  \newtheorem*{remark*}{Remark}
  \newtheorem{question}[thm]{Question}

  \def\R{{\mathbb R}}
  
  \def\Z{{\mathbb Z}}
  \def\N{{\mathbb N}}

  \def\e{{\varepsilon}}

  \def\C{{\mathcal C}}
  \def\Cper{{\mathcal C}_{\text{per}}}
  \def\Cinf{{\mathcal C}_\infty}

  \def\Sper{{S}_{\text{per}}}
  \def\PS0{\tilde{S}_0}
  \def\Sinf{{S}_\infty}
  \DeclareMathOperator{\Ends}{Ends}
  \DeclareMathOperator{\interior}{int}
  \DeclareMathOperator{\mcg}{Map}
  \newcommand{\Cis}{\C_{\text{p}}}
  \newcommand{\Cint}{\C_{\text{int}}}
  \newcommand{\Ca}{\C_{\text{a}}}

  \renewcommand{\L}{\mathcal{L}}
  \newcommand{\st}{\;|\;}
  \newcommand{\ssm}{\smallsetminus}
  \newcommand{\Teich}{Teichm\"uller }
  
\begin{document}


  \title[Towards Nielsen-Thurston classification]{Towards
  Nielsen-Thurston classification\\ for surfaces of infinite type:\\
  well-tempered homeomorphisms}
  
  \author{Mladen Bestvina}
  \address{Dept of Math, Univ of Utah, Salt Lake City, UT 84112, USA}
  \email{bestvina@math.utah.edu}
  \author{Federica Fanoni}
  \address{CNRS, Univ Paris Est Creteil, Univ Gustave Eiffel, LAMA UMR8050,
  F-94010 Creteil, France}
  \email{federica.fanoni@u-pec.fr} 
  \author{Jing Tao}
  \address{Dept of Math, Univ of Oklahoma, Norman, OK 73019, USA}
  \email{jing@ou.edu}
  
  \date{\today}

  \begin{abstract}
    We introduce and study tempered mapping classes of surfaces of infinite
    type. These are maps for which curves under iteration do not accumulate
    onto geodesic laminations with non-proper leaves, but only on unions of
    possibly intersecting curves or proper lines. Assuming an additional
    finiteness condition on the accumulation set, we prove a
    Nielsen-Thurston-type classification theorem. We prove that for such
    maps there is a canonical decomposition of the surface into invariant
    subsurfaces on which the first return is either periodic or a
    translation.
  \end{abstract}
   
  \thispagestyle{empty}

  \maketitle

\section{Introduction}

  In this paper, we study mapping classes -- that is, isotopy classes of
  self-homeomorphisms -- of a connected, orientable surface $S$ of
  \emph{infinite} type. For finite-type surfaces, the Nielsen–Thurston
  classification gives a complete picture. Every mapping class is periodic,
  reducible, or pseudo-Anosov. Periodic classes are represented by
  isometries of suitable hyperbolic metrics via Nielsen realization;
  pseudo-Anosov classes preserve a pair of transverse measured laminations,
  expanding one and contracting the other; reducible classes decompose the
  surface along a multicurve, reducing the dynamics to subsurfaces on which
  the first-return maps are again periodic or pseudo-Anosov. 
  
  The long-term goal of this project is to extend this understanding to
  surfaces of infinite type. Here, we take a step in that direction: we
  introduce the notion of \emph{tempered} mapping classes and prove a
  structure theorem for the subclass of \emph{well-tempered} mapping
  classes.
	
  To motivate our definitions, we briefly recall several approaches to the
  Nielsen-Thurston theorem and assess their prospects in the infinite
  setting. Thurston's original proof \cite{thurston_geometry,flp_travaux}
  is to find a fixed point in the compactified \Teich space. Bers' proof
  \cite{bers_extremal} also uses \Teich space, but from the point of view
  of extremal quasiconformal maps. Casson's proof \cite{cb_automorphisms}
  produces an invariant geodesic lamination by iterating a curve, while
  Nielsen's approach \cite{Nielsen_UntersuchungenI,
  Nielsen_UntersuchungenII, Nielsen_UntersuchungenIII}, completed by Miller
  \cite{miller} and Handel-Thurston \cite{ht_new}, analyzes the dynamics on
  the circle at infinity of lifts of the map to the universal cover. The
  Bestvina-Handel proof \cite{BH} produces invariant train tracks from
  efficient spines of the surfaces and relies on Perron-Frobenius theory. 

  In the infinite-type setting, however, \Teich space is much more
  complicated: it is infinite-dimensional, has uncountably many connected
  components, and some mapping classes do not preserve any component. So it
  is hard to imagine how Thurston's or Bers' approach could be adapted.
  Likewise, train-track methods would require infinite-sized matrices
  without a suitable analog of Perron-Frobenius theory. By contrast, the
  viewpoints of Casson and Nielsen-Handel-Thurston seem more amenable. In
  this paper, we take Casson's approach. 

  Let us say more carefully what Casson's method requires, because the
  difficulties in infinite type are already visible here. Given $f: S \to
  S$ and a curve $\alpha$, one takes a subsequence $f^{n_i}(\alpha)$, pulls
  each iterate tight in some hyperbolic metric, and attempts to extract a
  geodesic lamination $\lambda$ as a limit. For $\lambda$ to be
  useful--specifically for $\overline{\bigcup_n f^n(\lambda)}$ to be an
  $f$-invariant lamination without transverse intersections, the
  convergence must be ``robust'', meaning that more and more strands of
  $f^{n_i}(\alpha)$ should converge to the leaves of $\lambda$. On
  finite-type surfaces, this robustness is essentially automatic; it fails
  only when $\alpha$ intersects a subsurface on which $f$ is periodic, up
  to isotopy. On an infinite-type surface, it can fail in ways that are
  genuinely new and geometrically interesting.

  Three examples are worth keeping in mind. The first example comes from
  infinite-order isometries, such as the translation $(x,y) \mapsto
  (x+1,y)$ on $\mathbb{R}^2 \smallsetminus \Z^2$, which sends every curve
  to infinity. This is an instance of a general \emph{translation}--a
  homeomorphism that generates an infinite cyclic group acting properly on
  the surface. For a more complicated example, consider the same surface
  $\mathbb{R}^2 \smallsetminus \Z^2$ but consider a map that agrees with
  the translation on $y \le 1/3$, the inverse of the translation on $y\geq
  2/3$, and a continuous interpolation of the two by horizontal
  translations in the strip $1/3< y< 2/3$. In this case, a curve that
  intersects the strip essentially will converge to a line homotopic into
  the strip, but the convergence is not robust. A different phenomenon is
  an \emph{irrational rotation}. Start with a homeomorphism of the equator
  in $S^2$ whose minimal invariant set is a Cantor set $C$, and extend it
  to a homeomorphism $f$ of $S^2 \smallsetminus C$. Here, a convergent
  subsequence $f^{n_i}(\alpha)$ limits onto a line in a non-robust way, and
  the union of all such limits is a non-locally finite collection of lines,
  some of which intersect; in particular, it is not a geodesic lamination.
  These are not exotic counterexamples; they represent genuinely different
  dynamical behaviors, and any classification theorem must account for all
  of them.

  With these examples in mind, we arrive at the central definition of
  \emph{tempered} maps, which are precisely those mapping classes where
  robust convergence never occurs for any curve; in particular, a tempered
  mapping class does not exhibit pseudo-Anosov behavior anywhere on the
  surface. Concretely, $f: S \to S$ is tempered if for any two curves
  $\alpha,\beta$ on $S$, there is a uniform bound on the intersection
  numbers $i(f^n(\alpha),\beta)$ for all $n\in\Z$ (see Section
  \ref{sec:tempered}). As shown in Section \ref{sec:tempered}, if $f$ is
  tempered, then a curve $\alpha$ either leaves every compact set, or has
  limit set $\L(\alpha)$ consisting of curves or properly embedded lines
  (rather than more complicated laminations). A tempered map is \emph{well
  tempered} if, in addition, the limit set of every curve is a {\it finite}
  collection of lines and curves. These notions are well-defined on the
  mapping class of $f$. 

  On a finite-type surface, tempered is the same as periodic. On an
  infinite-type surface, the class is considerably richer. The basic
  examples are periodic maps, translations, and irrational rotations, with
  the first two being well tempered, and irrational rotations tempered but
  not well-tempered. More complicated examples can be built by gluing such
  maps together in a controlled way--for instance, by gluing two
  translations along a fixed line, as in the second example above--and this
  construction can be extended to infinitely many maps.
  
  Our main theorem shows that well-tempered mapping classes are exactly
  those arising from this kind of construction. In particular, they admit a
  satisfying Nielsen–Thurston-type theory.

  \begin{thmintro} \label{thm:mainintro}
    Let $\phi$ be a well-tempered mapping class of a surface $S$ of
    infinite type. Then there exists a representative $f \in \phi$ and an
    \(f\)-invariant 1-manifold \(\lambda \subset S\) without compact
    components such that, for the closure \(X\) of each component of
    \(S\ssm \lambda\), $f$ returns to $X$, and its first return is isotopic
    either to a periodic map or a translation.
  \end{thmintro}
  
  It is natural to ask whether Theorem \ref{thm:mainintro} admits an
  analogue for arbitrary tempered mapping classes. The full strength of the
  theorem can fail in several ways. First, one must allow for
  irrational-rotation behavior. Second, the fact that well-tempered maps
  are reduced along a 1-manifold is special to that setting; it ensures
  that the natural invariant pieces are always subsurfaces with boundary.
  For general tempered maps, by contrast, these pieces may have boundaries
  that accumulate in a complicated way, and hence need not be subsurfaces
  (see the \emph{Strategy of proof} section for more details). Even
  assuming local finiteness of limit sets of curves is not enough to
  circumvent this issue (see Section \ref{sec:examples}). Third, the
  existence of first-return maps is itself a nontrivial feature of the
  well-tempered setting and need not exist in general.

  Despite these complications, one might still hope that every tempered
  mapping class can be decomposed into pieces modeled on periodic maps,
  translations, and irrational rotation. Indeed, in the first version of
  this article, we conjectured such a picture, but it turns out to be
  incomplete. In addition to the three types above, there is a fourth type
  of dynamics, which we call an \emph{almost translation}. This is a map
  $f$ such that, for every two curves \(\alpha\) and \(\beta\), the set
  \(\{n\st i(f^n(\alpha),\beta)\neq 0\}\) has zero density in \(\Z\). 

  Taking these complications into account, we show in forthcoming work that
  for every tempered mapping class $\phi$, there exists a representative $f
  \in \phi$ and an $f$-invariant lamination $\lambda$ such that, for every
  complementary component of $S \setminus \lambda$, whenever a first return
  map exists, it is isotopic either to a periodic map, an almost
  translation, or an irrational rotation. We also show there that almost
  translations do occur. Whether the leaves of $\lambda$ can always be
  taken to be proper remains open.

  Beyond the tempered case, if $\phi$ is not tempered, one could still try
  to construct an invariant geodesic lamination using the Casson technique.
  This approach was carried out successfully for the class of
  \emph{irreducible end-periodic} maps by Handel and Miller
  \cite{ccf_endperiodic}. In our vision, a comprehensive structure theorem
  for an arbitrary mapping class $\phi$ should begin by decomposing $S$
  into invariant pieces on which the first return is either tempered or
  preserves two transverse geodesic laminations. Therefore, a structure
  theorem for tempered maps is an essential part of the broader program,
  and Theorem \ref{thm:mainintro} should be regarded as a first step.

  \subsection*{Strategy of proof}

  Given a map $f: S \to S$, our strategy for decomposing $S$ into
  $f$-invariant pieces is to group curves according to their dynamical
  behavior under $f$, and then analyze the subsurface \emph{spanned} by
  curves of a given type, when that subsurface exists. For instance, let
  $\Cper$ be the collection of $f$-periodic curves. If the subsurface
  \(\Sper\) spanned by $\Cper$ exists, then it is $f$-invariant (up to
  isotopy) and the action of $f$ on each component of \(\Sper\) is isotopic
  to a periodic map (see Section \ref{sec:structure}). 

  We apply the same idea to the collection $\C_\infty$ of \emph{wandering}
  curves, namely curves whose iterates leave every compact set of $S$. With
  some work, one can show that if the subsurface $S_\infty$ spanned by
  $\C_\infty$ exists, then the action of $f$ on each component of
  $S_\infty$ is isotopic to a translation. This may be viewed as an
  analogue of Brouwer’s plane translation theorem (\cite{brouwer_beweis})
  and it is proved in Section \ref{sec:translation}.

  \begin{thmintro}\label{thm:translationintro}

    Suppose $\phi$ is a mapping class of a surface $S$ such that every
    curve of $S$ is wandering. Then there is a representative $f \in \phi$
    and a hyperbolic metric on $S$ such that $f$ is an isometric
    translation.

  \end{thmintro}
  
  Thus, once the existence of \(\Sinf\) has been established, then Theorem
  \ref{thm:translationintro} reduces the proof of the main theorem to
  showing that the first return map on each component of \(S\ssm \Sinf\)
  exists and is periodic.

  This strategy makes the existence of spanning subsurfaces a central
  technical issue. For a finite collection of curves, the span is always
  well defined: put the curves in general position, take the union of their
  regular neighborhoods, and cap off any components homeomorphic to a disk
  with at most one puncture. For an infinite collection, one may take the
  span $S_n$ of the first $n$ curves and study the limiting object as $n
  \to \infty$. On a finite-type surface, Euler characteristic
  considerations force the topology of $S_n$ to stabilize, and this
  procedure essentially recovers the Nielsen-Thurston decomposition of
  $f$. On an infinite-type surface, however, the limiting object need not
  be a subsurface of $S$, but may instead be a complementary component of a
  geodesic lamination. Indeed, there are examples of tempered maps for
  which the collection of periodic or wandering curves fail to span a
  subsurface (see Section \ref{sec:examples}). Under the additional
  assumption that $f$ is well tempered, however, this bad behavior
  disappears, allowing the argument to go through. 

  \begin{propintro} \label{introprop:Subsurface}
    If $f$ is well tempered, then the collection of $f$-periodic curves
    spans a subsurface \(\Sper\), and the collection of $f$-wandering
    curves spans a subsurface \(\Sinf\). 
  \end{propintro}
  
  Proposition \ref{introprop:Subsurface} is proved in Section
  \ref{sec:decomposition} and follows from Lemmas \ref{lem:SinftySubs} and
  \ref{lem:SperSubs}.

  Note that periodic curves and wandering curves together do not, in
  general, span the whole surface -- said differently, \(\Sinf \cup
  \Sper\) is not always \(S\). We therefore need to study the complement of
  this union, denoted \(S_0\), and show that it contains no essential
  nonperipheral curves, and that the first-return map on each of its
  components exists and has finite order. In fact, Section
  \ref{sec:structure} shows a more technical statement than Theorem
  \ref{thm:mainintro} -- namely, Theorem \ref{thm:decomposition} --
  from which our main theorem is deduced.

  \subsection*{Connection to actions on graphs}
	  
  A related, and perhaps more tractable, variant of the Nielsen-Thurston
  classification theorem is to study the action of mapping classes on
  Gromov-hyperbolic graphs naturally associated to a surface, and to
  classify mapping classes by their dynamic types as graph isometries. For
  example, in the finite-type setting, pseudo-Anosov mapping classes act
  loxodromically on the curve graph of $S$, while all other mapping classes
  act elliptically.  
		
  In the infinite-type setting, the curve graph always has bounded
  diameter, reducing the classification of isometries to a trivial problem.
  However, many infinite-type surfaces admit other interesting
  Gromov-hyperbolic graphs of infinite diameter. One example is the ray
  graph, introduced by Bavard \cite{bavard_hyperbolicite} for the plane
  minus a Cantor set $\mathbb{R}^2 \smallsetminus C$, generalized to
  surfaces with one isolated end by Aramayona-Fossas-Parlier
  \cite{afp_arc}, and extended further by Bar-Natan-Verberne
  \cite{bv_grand}, who define the \emph{grand-arc graph}. In
  \cite{bavard_hyperbolicite}, Bavard also constructed a map of
  $\mathbb{R}^2 \smallsetminus C$ that acts loxodromically on the ray
  graph, and Abbott-Miller-Patel \cite{amp_infinite} generalized this
  construction to a large class of surfaces. However, a general
  classification of isometries remains open, and even the following
  question is unknown.
		
  \begin{question}
          
    Are there mapping classes of $\mathbb{R}^2 \smallsetminus C$ that act
    as parabolic isometries of the ray graph? More generally, are there
    mapping classes that act as parabolic isometries of the grand arc
    graphs?
          
  \end{question}
  
  By our work above, well-tempered mapping classes are elliptic isometries
  of the ray graph, and our future work on tempered mapping classes
  suggests that they are also elliptic isometries.

  \subsection*{Acknowledgments}
  
  We thank the American Institute of Mathematics for the hospitality during
  the workshop ``Surfaces of infinite type'', where this project began. We
  are also grateful to Peter Feller, Diana Hubbard, and Hannah Turner for
  sharing their preprint \cite{fht_dehn}, which was useful for refining our
  work on tempered maps. We thank Sanghoon Kwak for many useful comments on
  the first version of this paper, and the referee for their careful
  reading and comments. M.B.\ gratefully acknowledges support from the
  National Science Foundation under grant DMS-1905720.  F.F.\ gratefully acknowledges support from the ANR grants MAGIC
(ANR-23-TERC-0007) and GALS (ANR-23-CE40-0001). J.T.\ gratefully
  acknowledges support from the National Science Foundation under grants
  DMS-1651963 and DMS-2304920. F.F.\ also thanks Peter Feller for useful
  conversations.
  
  We acknowledge the use of ChatGPT to search for typos and minor mistakes.

\section{Background}

  \subsection*{Standing assumptions} 
  
  By a surface $S$ we will mean a connected orientable 2-manifold, without
  boundary unless otherwise stated.  When we need to consider surfaces with
  boundary, we will usually say \emph{bordered} surface for emphasis.  By a
  \emph{hyperbolic metric} on $S$ we always mean a complete hyperbolic
  metric without funnels or half-planes; equivalently, the hyperbolic
  surface coincides with its convex core, which is also equivalent to the
  metric being of the \emph{first kind}; see \cite{ar_structure}. For
  bordered surfaces, hyperbolic metrics are assumed to contain no funnels
  or half-planes and to have totally geodesic boundary.

  A surface $S$ has \emph{finite type} if $\pi_1(S)$ is finitely generated
  and \emph{infinite type} otherwise. By the classification of surfaces,
  $S$ is of finite type if and only if $S$ is homeomorphic to a closed
  genus $g$ surface with finitely many points removed.   A \emph{pair of
  pants} is a (possibly bordered) surface homeomorphic to either a
  three-holed sphere, or a two-holed once-punctured sphere, or a
  twice-punctured disk.
  
  Associated to $S$ are its \emph{ends}, which can be \emph{planar} or
  \emph{nonplanar}. A \emph{puncture} is an isolated planar end. We refer
  to \cite{av_big} for definitions and properties of ends of a surface.
    
  The \emph{mapping class group} of $S$ is the group $\mcg(S)$ of isotopy
  classes of orientation-preserving homeomorphisms of $S$.
  
  \emph{Curves} on $S$ are assumed to be simple, closed, and not homotopic
  to a point or a puncture. We will sometimes call such a curve
  \emph{essential} for emphasis. If $S$ has boundary, then a curve is
  called \emph{peripheral} if it is homotopic to a boundary component. When
  $S$ is endowed with a hyperbolic metric with geodesic boundary, every
  essential curve is realized by a unique geodesic $\alpha^*$ representing
  its homotopy class.
    
  A \emph{line} on $S$ is the image of a proper embedding of $\R \to S$. We
  will usually assume a line is \emph{essential}, meaning it does not bound
  a (topological) half-plane. We will usually consider lines up to proper
  isotopy relative to their ends. Similarly to curves, essential lines
  admit geodesic representatives:
  
  \begin{lemma}\label{lem:essential-to-geodesic}
    Let \(\ell\) be an essential line on a surface \(S\) endowed with a
    hyperbolic metric (of the first kind). Then \(\ell\) is isotopic to a
    unique geodesic \(\ell^*\).
  \end{lemma}

  \begin{proof}[Sketch of proof]
    Fix a pants decomposition of $S$, which yields an exhaustion of $S$ by
    finite-type subsurfaces with compact totally geodesic boundary (see
    \cite{ar_structure}). Pick a lift \(\tilde{\ell}\) of \(\ell\) in the
    universal cover. Since $\ell$ is proper and the metric is of the first
    kind, by looking at the lifts of the pants curves we can conclude that
    \(\tilde{\ell}\) has two distinct limit points on the circle at
    infinity. Let \(\tilde{\ell^*}\) be the unique geodesic joining the two
    limit points and let $\ell^*$ be the projection of $\tilde{\ell^*}$ to
    $S$. Arguing as in the proof of Proposition
    \ref{prop:straightening-laminations}, one shows that \(\ell\) is
    isotopic to \(\ell^*\). Uniqueness follows because the limit points of
    $\tilde \ell$ determine $\tilde \ell^*$ and hence $\ell^*$. \qedhere
  \end{proof}
  
  An \emph{arc} on $S$ is the image of an embedding $[0,1] \to S$; in
  particular in our definition arcs are compact. If $S$ has boundary, an
  arc with endpoints on $\partial S$ is called \emph{essential} if it
  doubles to an essential curve in the double of $S$.
  
  A \emph{ray} is the image of $[0,\infty)$ under a proper embedding. 
  
  A \emph{geodesic lamination} is a closed subset of $S$ that is a union of
  complete, simple, and pairwise disjoint geodesics (with respect to some
  hyperbolic metric on $S$), called the \emph{leaves} of the lamination.
    
  The \emph{geometric intersection number} $i(\alpha,\beta)$ is the minimum
  number of intersections between representatives of $\alpha$ and $\beta$
  in their homotopy classes, where each of $\alpha$ and $\beta$ could be
  either a line, curve, or a leaf of a lamination. Note that
  $i(\alpha,\beta)$ could be infinite, but we have the following
  characterization of lines.

  \begin{lemma} \label{lem:proper}

    Let $S$ be a surface equipped with a hyperbolic metric. Then a complete
    simple geodesic $\ell \subset S$ is a line or a curve if and only if
    $i(\ell,\alpha) < \infty$ for all curves $\alpha$.

  \end{lemma}

  \begin{proof}

    Clearly, if $i(\ell,\alpha) = \infty$ for some curve $\alpha$, then
    $\ell$ cannot be proper. Conversely, suppose $i(\ell,\alpha) < \infty$
    for all curves $\alpha$. Let $K \subset S$ be a compact subsurface.
    Since $\sum_{\alpha \subset \partial K} i(\ell, \alpha) < \infty$,
    $\ell \cap K$ has a finite number of components. If a component $\tau
    \subset \ell \cap K$ had infinite length, then $\tau$ would accumulate
    onto some geodesic lamination $\lambda \subset K$, or a boundary
    component $\gamma$ of $K$. In either case, any curve $\alpha$
    intersecting $\lambda$ (or $\gamma$) would satisfy $i(\ell,\alpha) \ge
    i(\tau,\alpha) = \infty$. Thus $\ell \cap K$ is always a finite union
    of arcs, and hence $\ell$ is proper. \qedhere 

  \end{proof}
  
  \subsection{Subsurfaces}

  A \emph{subsurface} of $S$ is a closed subset $X \subset S$ that is a
  two-dimensional (possibly disconnected) submanifold such that each
  boundary component of $X$ is either an essential curve or an essential
  line of $S$. We also require that no connected component of $X$ can be
  homotoped into another. In particular, no component of a subsurface can
  be a closed disk with at most one puncture. Moreover, $X$ is a subsurface
  if and only if $Y = \overline{S \setminus X}$ is also a subsurface. We
  will usually consider a subsurface up to isotopy.
  
  A connected subsurface $X$ is \emph{essential} if its double either has
  negative Euler characteristic or is of infinite type. Topologically, this
  rules out annuli and \emph{strips} (i.e.\ closed disks with two points
  removed from the boundary). 
  
  Two subsurfaces in $S$ are \emph{disjoint} if they have disjoint
  representatives. Similarly, a curve or line or geodesic lamination is
  disjoint from a subsurface if they have disjoint representatives. Note
  that in this convention, the boundary of a subsurface is considered
  disjoint from the subsurface. 

  Fix a hyperbolic metric on $S$. Because we will usually consider a
  subsurface up to isotopy, it is convenient to pick out a canonical
  representative within each isotopy class. To this end, we will introduce
  the notion of an \emph{(almost) geodesic representative} of a subsurface
  $X \subset S$, subject to the hyperbolic metric on $S$, which serves as
  an \emph{almost} canonical representative of the isotopy class of $X$.

  We begin with representatives of the interiors of subsurfaces. First
  suppose that $X$ is connected. A \emph{geodesic representative of the
  interior of $X$} is an open set $X^\circ$ defined as follows:

  \begin{itemize}
    \item If $X$ is essential, then $X^\circ$ is the open set properly
      homotopic to $\int(X)$ whose metric completion is a hyperbolic
      surface with totally geodesic boundary, homeomorphic to $X$.
    \item If $X$ is an annulus or a strip with core curve/line $\alpha$,
      then $X^\circ$ is the interior of any closed regular neighborhood of
      $\alpha^*$, whose boundary components are proper curves or lines.
  \end{itemize}
  
  For a (not necessarily connected) subsurface $X$, a \emph{geodesic
  representative of its interior} (also denoted $X^\circ$) is the union,
  over the connected components of $X$, of geodesic representatives of
  their interiors. By construction, $X^\circ$ is an open subset properly
  homotopic to $\int(X)$.

  To get a representative for $X$, we proceed as follows. For each
  component $Y$ of $X$, let $\overline{Y^\circ}$ be the closure of
  $Y^\circ$, and set $\partial Y^\circ := \overline{Y^\circ} \smallsetminus
  Y^\circ$. When $Y$ is essential, then $\partial Y^\circ$ is a disjoint
  union of simple closed geodesics and geodesic lines; otherwise, $\partial
  Y^\circ$ consists of a pair of almost geodesic curves or lines, each
  homotopic to the core curve/line of $Y$. Set $\partial X^\circ :=
  \bigcup_{Y \subset X} \partial Y^\circ$, and let $\partial_{sa} X^\circ$
  consist of those components $\alpha$ of $\partial X^\circ$ such that both
  boundary components of a regular neighborhood of $\alpha$ are contained
  in $X^\circ$ (informally, boundary components which will be
  ``self-adjacent''). We define an \emph{almost geodesic representative} of
  $X$ by $$X^*=\overline{X^\circ}\ssm \bigcup_{\alpha\in \partial_{sa}
  X^\circ} N(\alpha)$$ where $N(\alpha)$ is an open regular neighborhood
  of $\alpha$.

  An almost geodesic representative is well-defined up to our choices of
  regular neighborhoods of the components in $\partial_{sa}(X^\circ)$ and
  the representatives of the inessential components of $X$ (and the
  hyperbolic metric on $S$). In particular, if $\partial_{sa}(X^\circ) =
  \emptyset$ and $X$ has no inessential components, then $X^*$ is a totally
  geodesic subsurface representing $X$. Moreover, $\overline{X^\circ}$ is
  homeomorphic to $X^*$ if and only if $\partial_{sa}X^\circ = \emptyset$.

  \begin{figure}[htp!]
    \begin{center}
      \begin{overpic}{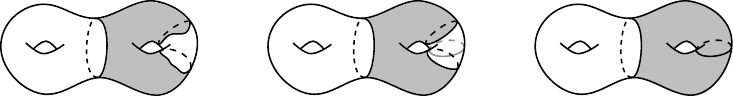}
      \put(16,3){\(X\)}
      \put(53,3){\(X^*\)}
      \put(88,3){\(\overline{X^\circ}\)}
      \put(49,12){\(\alpha\)}
      \put(86,12){\(\alpha\)}
      \put(64,5){\(\beta\)}
      \put(101,5){\(\beta\)}
      \end{overpic}
      \caption{An example of a subsurface $X$ so that $\overline{X^\circ}$
      is not homeomorphic to $X^*$. The curves \(\alpha\) and \(\beta\) are
      the geodesic representatives of the boundary components of \(X\).}
    \end{center}
  \end{figure}

  \begin{lemma}\label{lem:almost-geod-rep}
    For any hyperbolic metric on \(S\), any subsurface is properly
    homotopic to its almost geodesic representative.
  \end{lemma}
  
  \begin{proof}[Sketch of proof]
    If $\partial_{sa}X^\circ=\emptyset$ and $X$ has no inessential
    components, then the almost geodesic representative agrees with the
    geodesic representative of $X$, and the claim follows by properly
    homotoping all boundary components of \(X\) at once to their geodesic
    representatives (see Lemma \ref{lem:essential-to-geodesic} and
    \cite[Lemma 2.6]{cb_automorphisms}). If $\partial_{sa}X^\circ\neq
    \emptyset$ or $X$ has annular or strip components, then the same
    argument applies with a minor modification to handle pairs of parallel
    boundary components. \qedhere
  \end{proof}
  
  A consequence is that ``inclusion up to proper homotopy'' defines a
  partial order on the set of equivalence classes of subsurfaces defined up
  to proper homotopy:
  
  \begin{cor}
    Let \(S\) be a surface equipped with a hyperbolic metric. If \(X_1\)
    and \(X_2\) are subsurfaces such that each is properly homotopic
    into the other, then \(X_1\) and \(X_2\) are properly homotopic.
  \end{cor}
  
  \begin{proof}
  
    We first prove that inessential components of $X_1$ and $X_2$ are the
    same, up to proper homotopy. Indeed, suppose $Y$ is an inessential
    component of $X_1$. Since $X_1$ is properly homotopic into $X_2$, $Y$
    properly homotopes into some component $Z$ of $X_2$. Since $X_2$
    properly homotopes into $X_1$, $Z$ is properly homotopic into some
    component $Y'$ of $X_1$. Thus $Y$ properly homotopes into $Y'$, and by
    the definition of subsurface (no component homotopes into another), we
    must have $Y=Y'$. In particular, $Z$ properly homotopes into a strip or
    an annulus and is hence a strip or an annulus itself. By repeating the
    same argument for the inessential components of $X_2$ we get the claim.

    We may therefore assume that $X_1$ and $X_2$ have no inessential
    components. With respect to the hyperbolic metric on $S$, $X_1^\circ$
    and $X_2^\circ$ are uniquely defined. We claim they are the same.
    Otherwise, after swapping indices if necessary, there exists a
    component $Y$ of $X_1^\circ\ssm X_2^\circ$. Since $X_1$  properly
    homotopes into $X_2$, such $Y$ can only be a curve or a line, so
    $Y\subset \partial_{sa}X_1^\circ$. Choose an essential curve or line
    $\alpha\subset X_1^\circ$ that intersects $Y$ essentially. Again,
    because $X_1$ properly homotopes into $X_2$, $\alpha$ can be homotoped
    into $X_2$ and hence must be disjoint from $Y$, a contradiction. This
    shows $X_1^\circ = X_2^\circ$, so $X_1$ and $X_2$ have the same almost
    geodesic representative. Lemma \ref{lem:almost-geod-rep} then implies
    that $X_1$ and $X_2$ are properly homotopic.  \qedhere
  \end{proof}

  \subsection{Some results on isometries of surfaces}

  We say that a collection $\C$ of curves \emph{fills} a (sub)surface $X$
  if every nonperipheral curve in $X$ has positive intersection number
  with some curve in $\C$. The following well-known statement follows from
  the proof of \cite[Theorem 2.7]{cb_automorphisms}.

  \begin{lemma}\label{lem:periodic-finite-type}
    Let $F$ be a possibly bordered finite-type surface of negative Euler
    characteristic, whose boundary is either empty or compact. Let $f$ be a
    mapping class of $F$. If $F$ is filled by a finite collection of
    $f$-periodic curves, then $f$ has finite order.
  \end{lemma}

  \begin{remark}
    The order of $f$ in Lemma \ref{lem:periodic-finite-type} can be
    arbitrarily larger than the periods of the filling curves. For
    instance, consider the torus $\R^2 /\Z^2$ and the curves $\alpha$ and
    $\beta$ given by the projection of the lines $y=nx$ and $y=-nx$, for
    some positive integer $n$. Puncture the torus at the points that are
    the projections of $\left(\frac{1+2k}{2n},0\right)$ and
    $\left(\frac{k}{n},\frac{1}{2}\right)$, for $k=0,\dots, n-1$, which do
    not meet $\alpha$ and $\beta$. Then the map induced by $(x,y)\mapsto
    \left(x+\frac{1}{2n},y+\frac{1}{2}\right)$ has order $2n$, but it fixes
    both $\alpha$ and $\beta$.
  \end{remark}

  We will also need a realization result for periodic mapping classes of
  infinite-type surfaces; the proof is essentially borrowed from
  Afton-Calegari-Chen-Lyman (\cite{accl_Nielsen}).

  \begin{prop}\label{prop:Nielsenrealization}
    Let $S$ be an infinite-type surface with boundary, and $f$ a
    homeomorphism of $S$ such that $f^n$ is homotopic to the identity for
    some $n > 0$. For each compact boundary component $\alpha$ of $S$,
    choose a positive number $\ell(\alpha)$ that is $f$-invariant, i.e.\
    $\ell(f(\alpha))=\ell(\alpha)$. Then there exists a hyperbolic metric
    on $S$ such that each compact boundary component $\alpha$ has length
    $\ell(\alpha)$ and $f$ is isotopic to a periodic isometry of $S$.
  \end{prop}

  \begin{proof}
    Suppose first that $S$ has only compact boundary components. Using
    Nielsen's realization theorem for finite-type surfaces  with compact
    boundary \cite{nielsen_abbildungsklassen} and the argument of
    Afton-Calegari-Chen-Lyman (see the proof of \cite[Theorem
    2]{accl_Nielsen}), one obtains the desired metric and periodic
    representative in this case.

    Now suppose $S$ has noncompact boundary components. Let $DS$ be the
    double of $S$ along its noncompact boundary and extend $f$ to a
    homeomorphism $\hat{f}$ of $DS$. Let $i$ be the involution of $DS$
    with quotient $S$, which commutes with $\hat{f}$ by construction. By
    \cite{accl_Nielsen} and the discussion above, there is an
    $i$-invariant metric on $DS$ realizing the prescribed lengths for
    the compact components, and a periodic isometry $g$ isotopic to
    $\hat{f}$ that commutes with $i$. Since $i$ is an isometry, the
    noncompact boundaries of $S$ are realized as geodesics in $DS$, so
    that the metric descends to a metric on $S$ with totally geodesic
    boundary. Since $g$ commutes with $i$, it descends to a periodic
    isometry of $S$ isotopic to $f$. \qedhere
  \end{proof}

\section{Limit sets}

  Throughout this section, $S$ is a surface without boundary equipped with
  a hyperbolic metric, though everything we say will be independent of the
  particular choice of metric. Indeed, the following holds (see
  \cite[Theorem 3.6]{saric_train} for the infinite-type case). Recall
  that all hyperbolic metrics are of the first kind.
  
  \begin{thm}\label{thm:geods_correspondence} 
    Let $m$ and $m'$ be two hyperbolic metrics on $S$. The identity map on
    $S$ induces a natural identification between the boundaries at infinity
    of the universal covers of $(S,m)$ and $(S,m')$. This identification
    yields a homeomorphism between the spaces of complete geodesics in the
    universal covers of $(S,m)$ and $(S,m')$, which descends to a
    homeomorphism between complete $m$-geodesics and complete
    $m'$-geodesics on $S$.
  \end{thm}

  Given a complete $m$-geodesic $\ell$, we call the corresponding
  $m'$-geodesic the \emph{$m'$-straightening} of $\ell$. For a collection
  of $m$-geodesics, its \emph{$m'$-straightening} is the union of the
  $m'$-straightenings of its geodesics. A consequence of Theorem
  \ref{thm:geods_correspondence} is the following: 
    
  \begin{lemma}
    Let $m$ and $m'$ be two hyperbolic metrics on $S$. If two complete
    $m$-geodesics $\ell_1$ and $\ell_2$ are simple and disjoint, then the
    $m'$-straightenings $\ell_1'$ and $\ell_2'$ are also simple and
    disjoint. If $\lambda$ is an $m$-geodesic lamination, its
    $m'$-straightening is an $m'$-geodesic lamination.
  \end{lemma}
  
  The reason why the lemma holds is that (self-)intersections correspond to
  linked pairs of endpoints at infinity of lifts and convergence of
  geodesics corresponds to convergence of pairs of endpoints at infinity of
  lifts. Moreover, we have:
  
  \begin{prop}\label{prop:straightening-laminations}
    Let $m$ and $m'$ be hyperbolic metrics on $S$ and $\lambda$ an
    $m$-geodesic lamination. Then there is a homeomorphism of $S$, isotopic
    to the identity, that maps $\lambda$ leaf-by-leaf to its
    $m'$-straightening $\lambda'$.
  \end{prop}

  \begin{proof}
    Pick an $m$-geodesic pants decomposition $\mathcal{P}$ that contains
    no leaf of $\lambda$ (this can be done by choosing a pants
    decomposition and performing elementary moves until none of its curves
    lie in $\lambda$). In particular, $\mathcal{P}$ is in minimal position
    with respect to $\lambda$. Let $\mathcal{P}'$ be the
    $m'$-straightening of $\mathcal{P}$. There is an ambient isotopy $H_1$
    of $S$ sending $\mathcal{P}$ to $\mathcal{P}'$. Let $\lambda_1 :=
    H_1(\lambda)$. Then $\mathcal{P}'$ is in minimal position with respect
    to both $\lambda_1$ and $\lambda'$, and contains no leaf of either. 

    Now pick a curve $\alpha$ in $\mathcal{P}'$, and lift $\alpha$,
    $\lambda_1$, and $\lambda'$ to the universal cover of $(S,m')$. Choose
    and orient a lift $\beta$ of $\alpha$. Given a lift $\ell_1$ of a leaf
    in $\lambda_1$, its $m'$-straightening is the unique lift $\ell_1'$ of
    the $m'$-straightened leaf in $\lambda'$ with the same endpoints; in
    particular, $\ell_1$ intersects $\beta$ if and only if $\ell_1'$
    intersects $\beta$. Moreover, if $\ell_1$ and $\ell_2$ are two such
    lifts intersecting $\beta$, and $\ell_1'$ and $\ell_2'$ are their
    $m'$-straightenings, then the intersection points $\ell_1 \cap \beta$
    and $\ell_2 \cap \beta$ occur in the same order along $\beta$ as
    $\ell_1' \cap \beta$ and $\ell_2' \cap \beta$. Hence there is an
    order-preserving homeomorphism \[ \tilde{\lambda_1}\cap\beta \to
    \tilde{\lambda'}\cap\beta \] sending $\ell\cap\beta$ to $\ell'\cap
    \beta$, for every leaf $\ell\subset \tilde{\lambda_1}$, where $\ell'$
    is its $m'$-straightening. Extending this to an equivariant isotopy on
    a small neighborhood of $\beta$, and doing so equivariantly for every
    lift of every curve $\alpha$ in $\mathcal{P}'$ on pairwise disjoint
    neighborhoods, produces an equivariant isotopy of the universal cover.
    This descends to an ambient isotopy $H_2$ of $S$, sending $\lambda_1$
    to a lamination $\lambda_2$, such that
    $\lambda_2\cap\alpha=\lambda'\cap \alpha$ for all $\alpha$ in
    $\mathcal{P}'$.

    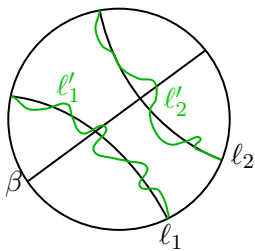
\begin{figure}[htp!]
    \begin{center}
    \begin{overpic}{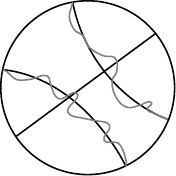}
    \put(-1,18){$\beta$}
    \put(68,-6){$\ell_1$}
    \put(100,30){$\ell_2$}
    \put(22,61){\color{MyGray}$\ell_1'$}
    \put(70,55){\color{MyGray}$\ell_2'$}
    \end{overpic}
    \caption{Lifting a curve in $\mathcal{P}'$}
    \end{center}
    \end{figure}

    Now let $P$ be a pair of pants of $\mathcal{P}'$. Since $\lambda_2$ and
    $\lambda'$ are transverse to $\partial P$, $\lambda_2\cap P$ and
    $\lambda'\cap P$ are unions of arcs with endpoints on $\partial P$. We
    claim that there is an isotopy of $P$, fixing $\partial P$ pointwise,
    sending $\lambda_2\cap P$ to $\lambda'\cap P$. Partition  the arcs of
    $\lambda_2 \cap P$ into (at most three) homotopy classes relative to
    the boundary. Do the same for the arcs of $\lambda' \cap P$; note that
    the same classes occur for both. For a class of $\lambda_2 \cap P$,
    choose a rectangle $R \subset P$, whose two opposite sides lie in
    $\partial P$ and whose two other sides are distinct arcs of
    $\lambda_2\cap P$, such that $R$ contains all arcs of $\lambda_2 \cap
    P$ in that class. Construct a corresponding rectangle $R' \subset P$
    for the same class of $\lambda' \cap P$, arranged so that $R$ and $R'$
    have the same four corners on $\partial P$. Then an isotopy of $P$
    fixing $\partial P$ pointwise can be chosen to send $R$ to $R'$ and the
    arcs of $\lambda_2\cap R$ to the corresponding arcs of $\lambda'\cap
    R'$. Doing this for each homotopy class of arcs yields the desired
    isotopy of $P$. 

    Repeating the procedure in each pair of pants produces an ambient
    isotopy $H_3$ of $S$, fixing each curve of $\mathcal{P}'$ pointwise,
    and sending $\lambda_2$ to $\lambda'$. The composition $H_3 \circ H_2
    \circ H_1$ is the required isotopy.\qedhere

    \begin{figure}[htp!]
    \begin{center}
    \begin{overpic}[scale=1.2]{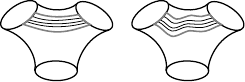}
    \put(20,35){$R$}
    \put(72,35){$R'$}
    \end{overpic}
    \caption{Isotopying rectangles}
    \end{center}
    \end{figure}

  \end{proof}

  Let $\lambda \subset S$ be a closed subset. We say that a sequence of
  curves $\{\alpha_n\st n\in \N\}$ \emph{converges} to $\lambda$, and that
  $\lambda$ is the \emph{limit} of $\{\alpha_n\st n\in \N\}$, if for every
  finite-type subsurface $K \subset S$ with compact boundary, $\alpha_n^*
  \cap K$ converges to $\lambda \cap K$ with respect to the Hausdorff
  distance. We will also denote this by $\alpha_n \to \lambda$. The same
  definition applies to sequences of geodesic arcs. Since convergence is
  tested on finite-type subsurfaces, the standard argument (see e.g.\
  \cite{cb_automorphisms}) implies: 

  \begin{lemma}\label{lem:limit}
    Let $\lambda \subset S$ be a closed subset. If $\lambda$ is the limit
    of a sequence of curves $\{\alpha_n\st n\in \N\}$, then $\lambda$ is a
    geodesic lamination. If $\lambda$ is the limit of a sequence of
    geodesic arcs, it is a union of pairwise disjoint geodesic arcs,
    geodesic rays and/or complete simple geodesics.
  \end{lemma}

  The \emph{limit set} of a sequence of curves $\{\alpha_n\st n\in \N\}$ is
  the set $\L(\{\alpha_n\st n\in \N\})$ consisting of all complete
  geodesics that occur as leaves of the limit of some subsequence of
  $\{\alpha_n\st n\in \N\}$. Similarly we can define the limit set of a
  sequence of geodesic arcs. Changing the metric on $S$ does not change the
  limit set, in the following sense. 

  \begin{lemma}
    Let $m$ and $m'$ be two hyperbolic metrics on $S$. For a sequence
    $\{\alpha_n\}$ of curves, let $\L(\{\alpha_n\st n\in \N\},m)$ and
    $\L(\{\alpha_n\st n\in \N\},m')$ be their respective limit sets in the
    two metrics. Then a subsequence $\{\alpha_{n_j}\}_{j\in\N}$ converges
    to $\lambda \subset \L(\{\alpha_n\st n\in\N\},m)$ if and only if
    $\{\alpha_{n_j}\}_{j\in\N}$ converges to $\lambda' \subset
    \L(\{\alpha_n\st n\in \N\},m')$, where $\lambda'$ is the
    $m'$-straightening of $\lambda$.
  \end{lemma}

  \begin{proof}
    As convergence of geodesics corresponds to convergence of pairs of
    endpoints at the boundary at infinity, this follows from Theorem
    \ref{thm:geods_correspondence}.\qedhere
  \end{proof}

  This justifies our notation $\L(\{\alpha_n\st n\in\N\})$, which makes no
  reference to the metric on $S$. 
    
  \subsection{Limit sets under a homeomorphism}

  For a homeomorphism $f$ of $S$ and a curve $\alpha$, we say $\alpha$ is
  $f$-\emph{periodic} if there exists $n \ge 1$ such that $f^n(\alpha)$ is
  isotopic to $\alpha$. The smallest such $n$ is called the
  $f$-\emph{period} of $\alpha$. We say $\alpha$ is $f$-\emph{forward
  wandering} if $\{f^n(\alpha)^*\}_{n \ge 0}$ leaves every compact set of
  $S$, and $f$-\emph{backward wandering} if it is $f^{-1}$-forward
  wandering. The curve $\alpha$ is $f$-\emph{wandering} if it is both
  $f$-forward and $f$-backward wandering. These properties are independent
  of the hyperbolic metric on $S$ as well as perturbing $f$ by an isotopy.
  Thus, we will also adopt the same definition for a mapping class of $S$.
  We will usually suppress the reference to $f$ when the context is clear.
  
  Note that a curve is $f$-wandering if and only if for every compact set
  $C \subset S$, $f^n(\alpha)$ can be homotoped to be disjoint from $C$ for
  all sufficiently large $|n|$. With this in mind, we can extend the
  definition of $f$-wandering to lines or subsurfaces of $S$, by which we
  mean that for every compact set $C \subset S$, the images under $f^n$ can
  be homotoped away from $C$ for any $n$ with sufficiently large $|n|$. 
  
  For a curve $\alpha$, let $\L^+(\alpha) := \L(\{f^n(\alpha)\st n \ge
  0\})$ be the \emph{forward limit} of $\alpha$ under $f$, and
  $\L^-(\alpha) := \L(\{f^n(\alpha)\st n \le 0\})$ the \emph{backward
  limit}. Set $\L(\alpha):=\L^+(\alpha)\cup\L^-(\alpha)$. When $\alpha$ is
  $f$-periodic, then \[ \L^{\pm}(\alpha) =  \{\alpha^*, f(\alpha)^*, \dots
  , f^{n-1}(\alpha)^*\}, \] where $n \ge 1$ is the $f$-period of $\alpha$.
  Moreover, $\alpha$ is forward wandering if and only if $\L^+(\alpha) =
  \emptyset$, and it is backward wandering if and only if $\L^-(\alpha) =
  \emptyset$. See Corollary \ref{cor:tempered-limits}.  
  
  We establish some basic facts.

  \begin{lemma}\label{lem:tempered-intersection}
    
    Let $f$ be a homeomorphism and $\alpha$ a curve. Suppose $\lambda$ is
    the limit of a sequence of iterates $\{f^{n_j}(\alpha)\st j\in\N\}$.
    Then for any curve $\beta$, $i(\beta,\lambda) \ne 0$ if and only if
    $i(\beta,f^{n_j}(\alpha)) \ne 0$ for all sufficiently large $j$.
    Moreover, if $i(\beta,\lambda) \ge k$, then $i(\beta,f^{n_j}(\alpha))
    \ge k$ for all sufficiently large $j$. 

  \end{lemma}

  \begin{proof}

    If we look at an annular neighborhood $K$ of $\beta$,
    $f^{n_j}(\alpha)^*\cap K$ converges to $\lambda \cap K$ in the
    Hausdorff metric, so the first statement follows. If $\beta$ intersects
    $\lambda$ at least $k$ times, then we can find $k$ arcs in $K \cap
    \lambda$. Choose sufficiently small neighborhood about each arc so that
    they are pairwise disjoint. Then for all sufficiently large $j$,
    $f^{n_j}(\alpha) \cap K$ will pass through each of these neighborhoods,
    so $i(\beta,f^{n_j}(\alpha)) \ge k$. \qedhere 

  \end{proof}

  \begin{lemma}\label{lem:duality}
    Let \(f\) be a homeomorphism. Then for any two (not necessarily
    distinct) curves \(\alpha\) and \(\beta\) we have:
    \[i(\alpha,\L^+(\beta))\neq 0 \; \Longleftrightarrow \;
    i(\L^-(\alpha),\beta)\neq 0.\]
  \end{lemma}

  \begin{proof}
    Suppose \(i \left( \alpha,\L^+(\beta) \right)\neq 0\). Let
    \(n_j\to\infty\) be a sequence so that $f^{n_j}(\beta) \to \lambda
    \subset \L^+(\beta)$ with \(i(\alpha,\lambda)\neq 0\). By Lemma
    \ref{lem:tempered-intersection}, for every \(j\) large enough,
    \[i(\alpha,f^{n_j}(\beta))\neq 0\]
    and thus
    \[i(f^{-n_j}(\alpha),\beta)\neq 0.\] In particular, $\alpha$ is not
    backward wandering. By taking a further subsequence so that
    $f^{-n_{j_k}}(\alpha) \to \lambda' \subset \L^-(\alpha)$, we have,
    again by Lemma \ref{lem:tempered-intersection}, \(i(\lambda',\beta)\neq
    0\), i.e.\ \(i(\L^-(\alpha),\beta)\neq 0\). By replacing $f$ by
    $f^{-1}$ we get the other implication. \qedhere 

  \end{proof}

  \begin{lemma}\label{lem:not-wandering} 

    Let $f$ be a homeomorphism. Then a curve $\alpha$ is forward wandering
    if and only if \(i(\alpha,\L^-(\beta)) = 0\) for all curves $\beta$.
    Similarly, $\alpha$ is backward wandering if and only if 
    \(i(\alpha,\L^+(\beta)) = 0\) for all $\beta$.

  \end{lemma}

  \begin{proof}
    
    If $\alpha$ is not forward wandering, then we can choose a curve
    $\beta$ with $i(\L^+(\alpha),\beta) \ne 0$. By Lemma
    \ref{lem:duality}, this implies $i(\alpha,\L^-(\beta)) \ne 0$. Conversely, if
    \(i(\alpha,\L^-(\beta))\neq 0\), then again by Lemma \ref{lem:duality},
    \(i(\L^+(\alpha),\beta)\neq 0\). In particular
    \(\L^+(\alpha)\neq\emptyset\), so \(\alpha\) is not forward wandering.
    \qedhere

  \end{proof}
  
  \begin{cor}\label{cor:not-wandering}
    Let $f$ be a homeomorphism. If $\alpha$ is a periodic curve and $\beta$
    is a wandering curve, $i(\alpha,\beta)=0$.
  \end{cor}
  \begin{proof}
    If $i(\alpha,\beta)\neq 0$, $i(\L^+(\alpha),\beta)\neq 0$ (as
    $\alpha\in \L^+(\alpha)$), so $\beta$ is not wandering by Lemma
    \ref{lem:not-wandering}.\qedhere
  \end{proof}

\section{Tempered maps}

  \label{sec:tempered}

  As usual, assume $S$ is a surface without boundary endowed with a
  hyperbolic metric. We say a homeomorphism $f$ of $S$ is \emph{tempered} if
  for every finite-type subsurface $K \subset S$ with compact boundary and
  every curve $\alpha$ there are $N=N(K,\alpha)\in \N$ and finitely many
  isotopy (relative to $\partial K$) classes of arcs from $\partial K$ to
  $\partial K$ and curves in $K$ such that, for every $n\in\N$, the
  intersection $f^n(\alpha)^*\cap K$ has at most $N$ components and each is
  in one of the given isotopy classes. Being tempered can be characterized by
  looking at intersections of curves, as the following lemma shows. 

  \begin{lemma} \label{lem:tempered-properties}
   
    Let $f$ be a homeomorphism of $S$. The following statements are
    equivalent. 
    
    \begin{enumerate}
      \item $f$ is tempered.
      \item $f^{-1}$ is tempered.
      \item For curves $\alpha$ and $\beta$, $i(f^n(\alpha),\beta)$ is
        uniformly bounded for all $n \ge 0$.
      \item For curves $\alpha$ and $\beta$, $i(f^n(\alpha),\beta)$ is
        uniformly bounded for all $n \le 0$.
      \item For curves $\alpha$ and $\beta$, $i(f^n(\alpha),\beta)$ is
        uniformly bounded for all $n$.
    \end{enumerate}
    In particular, being tempered is independent of the hyperbolic
    metric on $S$.
  \end{lemma}

  \begin{proof}

    The equivalence of (3), (4), and (5) are immediate as the action of $f$
    preserves geometric intersection number. That (1) implies (3) follows
    from taking the compact surface to be an annular neighborhood of the
    curve $\beta$. On the other hand, if $f$ is not tempered, then there
    exists a finite-type subsurface $K$ with compact boundary, a curve
    $\alpha$, and a subsequence $\{n_j\} \subset \N$ such that
    $f^{n_j}(\alpha)^* \cap K$ has $N_j$ components with $N_j \to \infty$
    and $j \to \infty$, or includes components with increasing
    multiplicity. Then there exists a curve $\beta\subset K$ such that
    $i(f^{n_j}(\alpha),\beta) \to \infty$. This shows the equivalence of
    (1) and (3), as well as (2) and (4) by symmetry, whence the equivalence
    of all statements. Property (5) is independent of the hyperbolic metric
    on $S$, so the same is true of being tempered. \qedhere 

  \end{proof}

  We say a mapping class $\phi \in \mcg(S)$ is \emph{tempered} if it has a
  tempered representative. From the characterization of Lemma
  \ref{lem:tempered-properties}, if $\phi$ has a tempered representative,
  then all of its representatives are tempered. We call a tempered
  homeomorphism $f$ \emph{well tempered} if the limit set \(\L(\alpha)\) is
  finite for every curve \(\alpha\). This property is also inherited by the
  mapping class of $f$. 
  
  As a point of reference, on a finite-type surface the following are
  equivalent: being tempered, well tempered, periodic, and isotopic to an
  isometry for some hyperbolic structure. In the infinite-type case, the
  situation is more complicated. Being periodic implies being isotopic to
  an isometry (with respect to some hyperbolic structure), which implies
  being well tempered, and well-tempered maps are tempered by definition.
  None of these implications is an equivalence. 

  For instance, consider a translation of a surface $S$, which -- as
  mentioned in the introduction -- is a map $f$ that generates an infinite
  cyclic group acting properly on $S$. If $S$ has infinite type, then the
  quotient surface $X=S/\langle f \rangle$ has negative Euler
  characteristic or has infinite type, so we can lift a hyperbolic metric
  from $X$ to $S$ so that $f$ acts by an (infinite-order) isometry. In this
  case, every curve on $S$ is wandering, so $f$ is (well) tempered. 
  
  Another example is what we call an \emph{irrational rotation}. Start with
  a homeomorphism of the circle with minimal invariant subset a Cantor set
  $C$ (see the Denjoy construction \cite[Section 12.2]{kh_introduction}).
  View the circle as the equator of the two-sphere $S^2$, and extend the
  map to a homeomorphism of $S=S^2 \smallsetminus C$. This map is tempered but
  \emph{not} well tempered: every curve has an infinite limit set, consisting
  of lines which may intersect. More complicated examples of (well) tempered
  maps can be found in Sections \ref{sec:examples} and \ref{sec:structure}.
  
  On the other hand, a Dehn twist is not tempered: if $\tau$ is the twist about
  a curve $\alpha$, let $\beta$ be a curve intersecting $\alpha$
  essentially. If $K$ is a finite-type subsurface with compact boundary
  containing $\alpha$ and $\beta$, the curves $\tau^n(\beta)\cap
  K=\tau^n(\beta)$ are all distinct. Another example of a map which is not
  tempered is a homeomorphism which restricts to a pseudo-Anosov on some
  finite-type subsurface (with compact boundary).

  The goal of the remainder of this section is to describe properties of
  limit sets of curves and lines under tempered mapping classes. Our first
  result is that if a curve or line is not periodic, then its limit set can
  only contain lines.

  \begin{prop} \label{prop:tempered-limits}
    
    Let $f$ be a tempered homeomorphism and $\alpha$ a nonperiodic curve. If
    $\alpha$ is not forward wandering, then $\L^+(\alpha)$ is a non-empty
    collection of lines. Similarly, if $\alpha$ is not backward wandering,
    then $\L^-(\alpha)$ is a non-empty collection of lines.

  \end{prop}

  \begin{proof}
    
    Assume $\alpha$ is not forward wandering, and for convenience, write
    $\alpha_n = f^n(\alpha)^*$. If $\L^+(\alpha) = \emptyset$, then for
    every compact subsurface $K \subset S$, the set $\bigcup_{n \ge 0}
    \alpha_n \cap K$ has no accumulation points. Equivalently, $\bigcup_{n
    \ge 0} \alpha_n \cap K$ consists of finitely many arcs or curves. But
    two geodesics that coincide along a non-trivial arc must be equal, and
    $\alpha$ is not $f$-periodic, we must have $\alpha_n \cap K =
    \emptyset$ for all sufficiently large $n$. But this contradicts the
    assumption that $\alpha$ is not forward wandering. By Lemma
    \ref{lem:limit}, $\L^+(\alpha)$ is a non-empty union of complete simple
    geodesics. If $\L^+(\alpha)$ contains a non-proper leaf $\ell$, then by
    Lemma \ref{lem:proper} there exists a curve $\beta$ such that
    $i(\beta,\ell) = \infty$.  Let $\lambda \subset \L^+(\alpha)$ be the
    limit of some subsequence $\{\alpha_{n_j}\}$ containing $\ell$. By
    Lemma \ref{lem:tempered-intersection}
    $i(f^{n_j}(\alpha),\beta)=i(\alpha_{n_j},\beta) \to \infty$,
    contradicting the characterization of being tempered given by Lemma
    \ref{lem:tempered-properties}. If instead $\L^+(\alpha)$ contains a
    compact leaf $\beta$, then since $\L^+(\alpha)$ has no non-proper
    leaves, $\beta$ arises as the limit of some subsequence
    $\{\alpha_{n_i}\}$. By the tempered condition, such a sequence must
    eventually be constant, but this contradicts the nonperiodicity of
    $\alpha$. Therefore $\L^+(\alpha)$ has no compact or non-proper leaves,
    so it is a non-empty union of lines. By replacing $f$ by $f^{-1}$ we
    get the second statement. \qedhere 

  \end{proof}
  
  As $\L^{\pm}(\alpha)$ is well defined for mapping classes, we obtain the
  following immediate corollary.
  
  \begin{cor} \label{cor:tempered-limits}
    For a tempered mapping class $f$ and a curve or line $\alpha$, the
    following statements hold. 
    \begin{enumerate}
      \item $\L^+(\alpha)=\emptyset$ if and only if $\alpha$ is
        forward wandering. 
      \item $\L^-(\alpha)=\emptyset$ if and only if $\alpha$ is
        backward wandering. 
      \item $\L(\alpha)=\emptyset$ if and only if $\alpha$
        is wandering. 
      \item $\L(\alpha)$ contains a curve if and only if
        $\alpha$ is a periodic curve, in which case
        $\L(\alpha)=\L^+(\alpha)=\L^-(\alpha)$ is the orbit of $\alpha$
        under $f$.
    \end{enumerate}
  \end{cor}

  Another consequence of being tempered is that limits of iterates of curves can
  intersect only finitely many times.

  \begin{lemma}
  
    Let $f$ be a tempered homeomorphism. For any two (not necessarily distinct)
    curves $\alpha$ and $\beta$, if $f^{n_j}(\alpha) \to \lambda \subset
    \L(\alpha)$ and $f^{m_k}(\beta) \to \lambda' \subset \L(\beta)$, for
    some sequences $\{n_j\}$ and $\{m_k\}$, then $i(\lambda,\lambda') <
    \infty$.

  \end{lemma}

  \begin{proof}

    If $i(\lambda,\lambda') = \infty$, then for every $N$ there are $j,k$
    such that $i(f^{n_j}(\alpha),f^{m_k}(\beta))\geq N$, by Lemma
    \ref{lem:tempered-intersection}. Thus, $i(\alpha, f^{m_k-n_j}(\beta))\geq
    N$, so $f$ is not tempered by Lemma \ref{lem:tempered-properties}.\qedhere
    
  \end{proof}
  
\section{Subsurfaces and spanning sets of curves}

  \label{sec:subsurfaces}
  Throughout this section, we fix a hyperbolic structure on $S$.  Given a
  collection of curves $\C$, we say a curve $\alpha$ is a
  \emph{concatenation} of curves in $\C$ if it is homotopic to a finite
  concatenation of geodesic arcs coming from the geodesic
  representatives of curves in $\C$. The collection $\C$ is \emph{closed
  under concatenation} if every curve obtained as a concatenation of curves
  in $\C$ is again an element of $\C$. 
  
  We say that a subsurface $\Sigma$ is \emph{spanned} by $\C$ if:
  \begin{enumerate}
    \item every curve in $\C$ is contained in $\Sigma$, up to isotopy; and
    \item every curve $\gamma \subset \Sigma$ is a concatenation
      of the curves in $\C$
  \end{enumerate}
  and $\Sigma$ is minimal (with respect to inclusion, up to proper
  homotopy) among all subsurfaces satisfying (1) and (2).

\begin{figure}[htp!]
  \begin{center}
  \begin{overpic}[scale=.9]{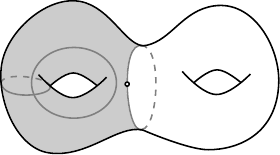}
    \put(-10,22){$\Sigma$}
  \end{overpic}
  \caption{A subsurface $\Sigma$ not spanned by its
    curves.}\label{fig:crowned-torus}
  \end{center}
  \end{figure}

  Note that if $\C$ spans a subsurface $\Sigma$, then necessarily $\C$
  fills $\Sigma$.  On the other hand, a (sub)surface $\Sigma$ need not
  coincide with the subsurface spanned by its curves. For instance, on a
  once-punctured surface of genus two we can consider the subsurface
  $\Sigma$, homeomorphic to a one-holed torus with a point removed from the
  boundary, as in Figure \ref{fig:crowned-torus}. Here, the subsurface
  spanned by the curves of \(\Sigma\) is a one-holed torus with compact
  boundary. 

  We say that $X$ is \emph{spanned by its curves} if it coincides with the
  subsurface spanned by its curves.

  When $\C$ is a finite collection, there is always a subsurface $\Sigma$
  (of finite type) spanned by $\C$. Namely, put the curves in $\C$ in
  general position, take a regular neighborhood of their union, and fill in
  any complementary disks with at most one puncture.
  
  In general, one can show that if a subsurface spanned by a collection
  $\C$ of curves exists, then it is unique up to proper homotopy. To prove
  this fact, we introduce three subcollections of curves, which we will
  also use later:
  \begin{align*}
    &\Cis :=\{\gamma\in\C\st \forall\gamma'\in\C,\; i(\gamma,\gamma')=0\};\\
    &\Cint :=\C\ssm\Cis;\\
    &\Ca :=\{\gamma\in\Cis\st \gamma \text{ is not a concatenation of
    curves in }\Cint\};
  \end{align*}
  If a subsurface spanned by $\C$ exists, the curves in \(\Cis\) are
  peripheral. Among them, the curves in $\Ca$ are precisely those that
  contribute annular components to the spanned subsurface, while the curves
  in $\Cis\ssm\Ca$ occur as boundary components of non-annular components.
  The curves in \(\Cint\) instead are in the interior of the subsurface.

  \begin{lemma}\label{lem:spanned-unique}
    Let $\C$ be a collection of curves in $S$. If $\Sigma_1$ and $\Sigma_2$
    are subsurfaces spanned by $\C$, then they are properly homotopic.
  \end{lemma}
  
  \begin{proof}
    Fix a hyperbolic structure on $S$, and assume that all curves in $\C$
    are realized as geodesics. Let $\Sigma_i^\circ$ be the geodesic
    representative of the interiors of $\Sigma_i$, chosen so that
    $\partial\Sigma_1^\circ$ and $\partial\Sigma_2^\circ$ are in minimal
    position. Let $\Sigma$ be an almost geodesic representative of
    $\Sigma_1^\circ\cap\Sigma_2^\circ$. Then $\Sigma$ is properly homotopic
    into both $\Sigma_1$ and $\Sigma_2$. We claim that $\Sigma$ contains
    all curves in $\C$ and that all curves in $\Sigma$ are concatenations
    of curves in $\C$. Then, by minimality of $\Sigma_1$ and $\Sigma_2$,
    they both are properly homotopic to $\Sigma$ and hence to each other.
  
    It remains to prove the claim. 
  \begin{enumerate}
    \item Let $\gamma\in\C$. 
      \begin{itemize}
        \item If $\gamma\notin \Cis$, then $\gamma$ is a
          nonperipheral curve of a subsurface contained both in
          $\Sigma_i$. Hence it is contained in both $\Sigma_i^\circ$, and
          therefore $\gamma\subset \Sigma$, up to homotopy. 
        \item If $\gamma\in \Ca$, then $\gamma$ corresponds to an annular
          component contained in both $\Sigma_i$. Hence it is contained in
          both $\Sigma_i^\circ$, and therefore $\gamma \subset \Sigma$, up to
          homotopy. 
        \item If $\gamma\in\Cis\ssm \Ca$, then there are finitely many
          curves in $\C\ssm\Cis$, say $\beta_1,\dots,\beta_k$, such that
          $\gamma$ is a concatenation of arcs of the $\beta_i$. In
          particular, $\gamma$ is contained in the subsurface spanned by
          $\{\beta_1, \ldots, \beta_k\}$, which is contained, up to
          homotopy, in both $\Sigma_i$ and hence in $\Sigma_i^\circ$.
          Therefore, $\gamma\subset \Sigma$, up to homotopy.
      \end{itemize}
    \item Let $\gamma\subset \Sigma$. Then $\gamma$ is up to
      homotopy contained in both $\Sigma_i$. Since each $\Sigma_i$ is
      spanned by $\C$, $\gamma$ is  a concatenation of curves in $\C$.
      \qedhere 
  \end{enumerate}
  \end{proof}

\begin{figure}[htp!]
  \begin{center}
  \begin{overpic}[scale=.9]{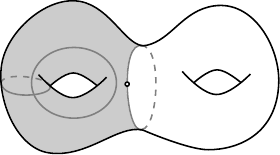}
  \put(5,15){$\alpha$}
  \put(20,6){$\beta$}
    \put(-10,22){$\Sigma$}
  \end{overpic}
  \caption{A (non-minimal) subsurface $\Sigma$ whose curves are
    concatenations of $\alpha$ and $\beta$.}
    \label{fig:crowned-torus-non-minimal}
  \end{center}
  \end{figure}

  Note that uniqueness requires the minimality assumption, even for finite
  collections of curves on finite-type surfaces. For example, in Figure
  \ref{fig:crowned-torus-non-minimal}, let $\C$ be a pair of curves
  $\alpha$ and $\beta$ that intersect once. Then a one-holed torus with a
  point removed from the boundary satisfies conditions (1) and (2), but it
  is not minimal. The subsurface spanned by $\C$ is a one-holed torus (with
  compact boundary).
  
\begin{figure}[htp!]
  \begin{center}
  \includegraphics[scale=.9]{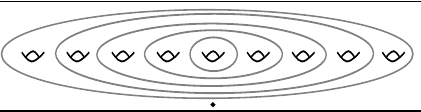}
  \caption{The beginning of a collection of curves not spanning a
    subsurface.}\label{fig:notspanning}
  \end{center}
  \end{figure}

  If $\C$ is an infinite collection of curves, then a subsurface spanned by
  $\C$ does not always exist -- see for instance Figure
  \ref{fig:notspanning}. A simple condition that ensures the existence of a
  spanned subsurface is that $\C$ be \emph{locally finite}, meaning that
  every compact set of $S$ intersects only finitely many curves in $\C$.
  The main goal of this section is to give a condition on $\C$ that
  guarantees it spans a subsurface. In fact, under the additional
  assumption that $\mathcal{C}$ is closed under concatenation, we give a
  necessary and sufficient condition for a collection $\C$ to span a
  subsurface (Proposition \ref{prop:span}). This assumption holds for some
  natural collections of curves, such as curves which are periodic or
  wandering with respect to a mapping class:
  
  \begin{lemma}\label{lem:periodic&wandering}

    Let $f$ be a mapping class on $S$. Every curve obtained as a
    concatenation of $f$-periodic curves is periodic, and every curve
    obtained as a concatenation of wandering curves is
    wandering. 

  \end{lemma}

  \begin{proof}

    We first consider periodic curves. A curve is a concatenation of curves
    $\gamma_1,\dots, \gamma_n$ if and only if it is contained in the
    subsurface $F$ spanned by $\{\gamma_1,\ldots,\gamma_n\}$. Note that $F$
    is connected and has negative Euler characteristic. Since each
    $\gamma_i$ is $f$-periodic, we can take a suitable power so that $f^k$
    fixes $\gamma_i$ for all $i$. Since $F$ is spanned by $\{\gamma_i\}$,
    $f^k(F)$ is spanned by $\{f^k(\gamma_i)\}=\{\gamma_i\}$, so $f^k(F)$ is
    isotopic to $F$. Since $F$ has negative Euler characteristic, $f^k$
    restricted to $F$ is isotopic to a periodic map of $F$ (see Lemma
    \ref{lem:periodic-finite-type}), so every curve in $F$ is
    $f$-periodic. 
    
    The proof for wandering curves is similar. Let $F$ be the subsurface
    spanned by a finite collection $\{\gamma_1,\ldots,\gamma_n\}$ of
    wandering curves. It suffices to show that for all compact subsurface
    $K \subset S$, $f^{\pm k}(F)$ can be homotoped to be disjoint from $K$
    for all sufficiently large $k$. Since each $\gamma_i$ is wandering, we
    can find $k_0$ such that $f^{\pm k}(\gamma_i)^* \cap K = \emptyset$ for
    all $k \ge k_0$ and all $i$. Then $f^{\pm k}(F)$ is spanned by
    $\{f^{\pm k}(\gamma_i)\}$, so $f^{\pm k}(F)$ can be homotoped away from
    $K$. \qedhere

  \end{proof}
  
  As mentioned before, if $\C$ is a finite collection of curves then there
  is an explicit way to construct the subsurface spanned by $\C$. Given any
  collection of curves, our goal is to describe a procedure that yields the
  subsurface spanned by them, if it exists, and to understand conditions
  under which it does exist.

  Recall the sets $\Cis, \Cint$ and $\Ca$ defined before:
  \begin{align*}
    &\Cis :=\{\gamma\in\C\st \forall\gamma'\in\C,\; i(\gamma,\gamma')=0\}\\
    &\Cint :=\C\ssm\Cis\\
    &\Ca :=\{\gamma\in\Cis\st \gamma \text{ is not a concatenation of
    curves in }\Cint\}
  \end{align*}
  Enumerate the curves in \(\Cint=\{\gamma_0,\gamma_1,\dots\}\).
  Up to reordering, we can assume that $i(\gamma_0,\gamma_1) \ne 0$. Then
  for every \(i\geq 1\), let \[\Cint^{i}:=\{\gamma_j\st j\leq i
  \;\mbox{and}\; \exists\;j'\leq i \text{ such that }
  i(\gamma_j,\gamma_{j'})\neq 0\}.\] Note that $\Cint^{1} =
  \{\gamma_0,\gamma_1\}$, $\Cint^{i} \subset \Cint^{i+1}$, and
  \(\Cint=\bigcup_{i\geq 1}\Cint^{i}\). Moreover, if
  $\Cint^i\neq\Cint^{i+1}$, to obtain $\Cint^{i+1}$ we add curves which
  either intersect curves in $\Cint^i$ or intersect each other.  In
  particular, each component of the subsurface spanned by $\Cint^{i}$ has
  negative Euler characteristic. 
    
  We now let \(F^{i}=F(\Cint^i)\) be the geodesic representative of the
  \emph{interior} of the subsurface spanned by \(\Cint^{i}\) and define
    \[F=F(\Cint):=\bigcup_{i\geq 1}F^{i}.\]
    By construction every connected component of \(F^{i}\) is an open subset
    (homotopic to a subsurface), and thus $F$ is an open subset of $S$. We
    will call $F$ the \emph{set spanned} by $\Cint$.
    
    We also define \(F(\C)\) to be the union
     \[F(\C):=F(\Cint)\cup\bigcup_{\alpha\in\Ca}A(\alpha),\]
     where, for every \(\alpha\), \(A(\alpha)\) is a regular neighborhood
     of the geodesic representative of \(\alpha\), chosen so that they are
     pairwise disjoint and disjoint from \(F(\Cint)\). 
    
    We will call a collection $\C$ of curves \emph{good} if:
    \begin{enumerate}
      \item[(a)] the curves in \(\Ca\) when realized as geodesics, form a
        locally finite collection, i.e.\ they do not accumulate anywhere;
        and
      \item[(b)] for every \(p\in\partial F\) there is a closed disk \(B\)
        centered at \(p\) such that \(B\cap F\) is either \(B\) minus a
        diameter, or a connected component of \(B\) minus a diameter. 
    \end{enumerate}

    \begin{figure}[htp!]
    \begin{center}
    \begin{overpic}{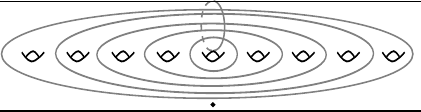}
    \put(49,-4){$p$}
    \end{overpic}
    \vspace{.2cm}
    \caption{A collection of curves which is not good.}\label{fig:notspanning2}
    \end{center}
    \end{figure}
    
    For instance, the collection of curves depicted in Figure
    \ref{fig:notspanning} is not good, because the isolated curves
    accumulate somewhere. Similarly, the collection in Figure
    \ref{fig:notspanning2} is not good either, because condition (b) is not
    satisfied at the point \(p\) in the drawing (onto which the curves
    accumulate). In both examples, the collection of curves does not span a
    subsurface. Indeed, we will show that the condition of being good is
    equivalent to $\C$ admitting a spanning subsurface (provided $\C$ is
    closed under concatenation).
    
    Let $K \subset S$ be a finite-type subsurface with compact totally
    geodesic boundary. Each $F^{i}$ intersects $K$ in a disjoint union of
    open subsurfaces and disks with at most one puncture. By a
    \emph{rectangle} we mean a connected component $R$ of $K \cap F^{i}$
    that is homeomorphic to a closed disk and whose boundary is split into
    \emph{horizontal sides} lying in $\partial F^{i}$ and \emph{vertical
    sides} lying in $\partial K$. All components of $K \cap F^{i}$ that are
    not rectangles are called \emph{essential}. Let $K^{i}_e$ be the union
    of all the essential components of $K \cap F^{i}$.  
    
    \begin{figure}[htp!]
      \begin{center}
        \includegraphics[scale=.8]{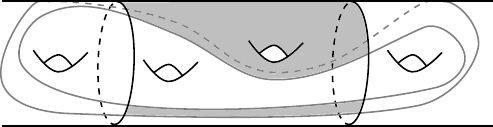}
        \caption{A compact subsurface $K$, whose boundary is in black, and
        the boundary of a subsurface $F^i$ in gray. The shaded part
        is the intersection, formed by an essential component and a
        rectangle.}
      \end{center}
    \end{figure}

    \begin{lemma}\label{lem:stable-essential}
      For all $i \ge 2$, $K^{i}_e \subset K^{i+1}_e$, and for large $i$ it
      stabilizes. More precisely, there exists an open subsurface
      $K_e\subset K$ such that for all sufficiently large $i$, $K^{i}_e$ is
      isotopic to $K_e$.
    \end{lemma}

    \begin{proof}
      This follows from Euler characteristic considerations. \qedhere
    \end{proof}

    So for any finite-type subsurface $K$ with compact totally geodesic
    boundary, there exists $I$ such that for all $i \ge I$, the topology of
    the essential components $K^{i}_e$ stabilizes, and further increasing
    $i$ only adds new rectangles or adjoins previous ones.  We define:
      \[ \mathcal{R}(K) =  \{ R \st R \subset K \cap F^{i} \text{
      is a rectangle for some $i \ge I$ and is not contained in any
      $K^i_e$} \}.\] Two parallel rectangles in $\mathcal R(K)$ are
    declared \emph{equivalent} if they are contained in the same connected
    component of $F^i\cap K$ for some $i$. The following lemma is not hard
    to prove, and we leave it as an exercise.

    \begin{lemma}
      The relation defined above is an equivalence relation on
      $\mathcal{R}(K)$. Moreover, suppose $R_1$, $R_2$, $R_3$ are three
      parallel rectangles with $R_2$ contained in the rectangle between
      $R_1$ and $R_3$. If $R_1$ and $R_3$ are equivalent, then all three
      are equivalent.
    \end{lemma}

    This allows us to naturally organize parallel rectangles into
    equivalence classes. As $i$ increases, the number of equivalence
    classes may increase. We have the following proposition.
    
    \begin{prop}\label{prop:span}
      Suppose \(\C\) is a collection of curves closed under concatenation.
      The following statements are equivalent:
      \begin{enumerate}
      \item There is a subsurface spanned by \(\C\). 
      \item The curves in \(\Ca\) do not accumulate anywhere and for every
        compact subsurface \(K\) with totally geodesic boundary, the number
        of equivalence classes in \(\mathcal{R}(K)\) is finite.
      \item \(\C\) is good.
      \end{enumerate}
    \end{prop}

    \begin{proof}

      [(1) \(\Rightarrow\) (2)] Let $\Sigma$ be a surface spanned by $\C$.
      By contradiction, suppose the curves in $\Ca$ accumulate somewhere.
      Then we can find a sequence of curves $\gamma_j\in\Ca$ and a curve
      $\alpha$ such that $i(\gamma_j,\alpha)\neq 0$ for all $j$. Since
      $\Sigma$ is a subsurface, $\Sigma\cap \alpha$ consists of finitely
      many arcs, with at least one, denoted $a$, intersecting at least
      three curves, say $\gamma_{j_1},\gamma_{j_2}$ and $\gamma_{j_3}$. Up
      to replacing \(a\) with a subarc and changing the roles of the
      curves, we can assume that $\gamma_{j_1}\cap a$ and $\gamma_{j_3}\cap
      a$ are reduced to a point and there is at least one intersection
      point of \(\gamma_{j_2}\) with \(a\) between $\gamma_{j_1}\cap a$ and
      $\gamma_{j_3}\cap a$. 
      
      Then the pair of pants containing $\gamma_{j_1}\cup a\cup
      \gamma_{j_3}$  is in $\Sigma$ and contains a curve $\beta$ which
      intersects $\gamma_{j_2}$ essentially. As $\beta \subset \Sigma$, it
      is a concatenation of curves in $\C$, and since $\C$ is closed under
      concatenation we have $\beta\in\C$, contradicting the definition of
      $\Ca$. This shows that the curves in $\Ca$ cannot accumulate
      anywhere.

      Next, suppose there is a finite-type subsurface $K$ with compact
      totally geodesic boundary such that $\mathcal{R}(K)$ has infinitely
      many equivalence classes. Then we can find infinitely many parallel
      rectangles $R_i$ that are not equivalent. Let $\gamma_i$ be a curve
      in $F$ passing through $R_i$. By construction, each $\gamma_i$ is a
      concatenation of subarcs of curves in $\C$, hence $\gamma_i \subset
      \Sigma$. Since $\Sigma$ is a subsurface, $\Sigma\cap K$ is a finite
      union of connected components. In particular, there is some component
      $X$ of $\Sigma\cap K$ containing $\gamma_i\cap K$ for infinitely many
      $i$. This implies that we can find distinct curves $\gamma_{i_1}$ and
      $\gamma_{i_2}$ such that the rectangle $R$ between $\gamma_{i_1}\cap
      K$ and $\gamma_{i_2}\cap K$ is included in $\Sigma$. This implies that $R$ is
      contained in $F^l\cap K$ for some $l$, so $R_{i_1}$ and $R_{i_2}$ are
      equivalent, a contradiction.

      [(2) \(\Rightarrow\) (3)] We just need to verify property (b) in the
      definition of ``good''. By construction, $\partial F$ is a union of
      limits of simple closed geodesics. We claim that for every
      finite-type subsurface $K$ with compact totally geodesic boundary,
      $\partial F\cap K$ has finitely many connected components. Otherwise,
      as $K$ is of finite type, there are infinitely many components $l_j$
      of $\partial F \cap K$ parallel to a given arc from $\partial K$ to
      $\partial K$. Then there are infinitely many rectangles between the
      $l_j$ which are not in $F$ and therefore infinitely many rectangles
      that are not equivalent, a contradiction.

      Therefore, for every $p\in\ell\subset \partial F$ we can find a
      sufficiently small ball $B$ such that $\ell\cap B$ coincides with
      $\partial F\cap B$ and it is a single arc from $\partial B$ to
      $\partial B$. In particular, $B\ssm \partial F=B\ssm \ell$ has two
      components, each of which is either entirely contained in $F$ or in
      $S\ssm F$. Moreover, since $\ell\subset \partial F$, at least one of
      the components intersects $F$, and hence is contained in $F$. So $B$
      is the required disk.
      
      [(3) \(\Rightarrow\) (1)] The condition on boundary points implies
      that $\partial F$ contains no non-proper component and that its
      components do not accumulate anywhere. If $\gamma\subset \partial F$
      is a compact boundary component, we can find some finite-type
      subsurface $K$ with compact boundary components such that
      $\gamma\subset K$. As $K^i_e$ stabilizes, $\gamma$ is eventually a
      boundary component of some $F^i$, and therefore is an essential
      curve. If $\gamma\subset \partial F$ is a non-compact boundary
      component, since it is the limit of simple closed geodesics, it is
      geodesic, and therefore it is a line. Moreover, since each component
      of $\partial F$ is a simple closed geodesic or the limit of simple
      closed geodesics, for any $\ell\subset \partial F$, either
    \begin{itemize}
      \item for every $p\in\ell$ there is a closed disk $B$ centered at $p$
        such that $B\cap F$ is a connected component of a disk minus a
        diameter, or
      \item for every $p\in\ell$ there is a closed disk $B$ centered at $p$
        such that $B\cap F$ is a disk minus a diameter.
    \end{itemize}
    Let $T_1$ be the collection of lines and curves in $\partial F$ for
    which the first condition holds, and $T_2$ be the collection of lines
    and curves for which the second holds.

    For every $\ell\in T_2 \cup \Ca$ we can find a regular neighborhood
    $N(\ell)$ such that:
    \begin{itemize}
      \item $N(\ell)$ is open if $\ell\in T_2$ and closed if $\ell\in \Ca$;
      \item $N(\ell)$ is disjoint from $\partial F\cup \Ca\ssm\{\ell\}$; 
      \item if $\ell_1\neq \ell_2$, then $N(\ell_1)\cap
        N(\ell_2)=\emptyset$, and
      \item the neighborhoods do not accumulate anywhere. 
    \end{itemize}
    We can find such neighborhoods by the properness of components of
    $\partial F$, together with the non-accumulation assumption on $\Ca$
    (we can fix a compact exhaustion of the surface and construct the
    neighborhoods piece by piece).

    Define $$\Sigma:=\left(F\cup\bigcup_{\ell\in T_1}\ell \cup
    \bigcup_{\gamma\in\Ca}N(\gamma) \right)\ssm \bigcup_{\ell\in
    T_2}N(\ell).$$ Then $\Sigma$ is a two-dimensional submanifold. Its
    boundary components are either homotopic to boundary components of $F$
    (hence are essential curves or lines), or homotopic to curves in $\Ca$
    (hence are essential curves). Therefore, $\Sigma$ is a subsurface.

    It is clear that $\Sigma$ contains all curves in $\C$. Moreover, if
    $\gamma$ is a curve in $\Sigma$, either it is contained in an annular
    component, and therefore is in $\C$, or it is homotopic into some $F^i$
    (by the same argument for compact components of $\partial F$), and
    therefore $\gamma$ is a concatenation of curves in $\C$. Finally,
    $\Sigma$ is properly homotopic to $F\cup\bigcup_{\gamma\in\Ca}
    N(\gamma)$. Since any surface containing all curves in $\C$ must
    contain $F\cup\bigcup_{\gamma\in\Ca}\gamma$ up to proper isotopy,
    $\Sigma$ satisfies the minimality hypothesis, and hence is a
    subsurface spanned by $\C$. \qedhere

  \end{proof}
  
  \begin{cor}\label{cor:spanned-disjoint}
    If $X$ and $Y$ are subsurfaces spanned by two collections of curves
    $\C_X$ and $\C_Y$ and $i(\alpha,\beta)=0$ for every $\alpha\in
    \C_X,\beta\in\C_Y$, then $X$ and $Y$ are disjoint.
  \end{cor}
  
  \begin{proof}
    This follows from the explicit construction of a subsurface spanned by
    a set of curves given in the proof of Proposition \ref{prop:span},
    together with the uniqueness of spanned subsurfaces (Lemma
    \ref{lem:spanned-unique}). \qedhere
  \end{proof}
 
  We record a consequence of having a non-good collection of curves. This
  will be our main tool for showing certain collections of curves are good.

  \begin{lemma}\label{lem:not-good}
    Let $\C$ be a collection of curves closed under finite concatenations.
    Suppose $\C$ is not good. Then there is a simple closed geodesic
    $\alpha$ and two infinite sequences of simple closed geodesics
    $\beta_j,\gamma_j$, $j \in \N$, with the following properties:
    \begin{enumerate}
      \item $\alpha, \gamma_j \notin \C$ but $\beta_j\in \C$, for all
        $j$.
      \item every $\beta_j$ intersects $\alpha$ and the
        $\beta_j$'s are
        pairwise disjoint. 
      \item For every $j$, there is a subarc $\tau_j \subset \alpha$ joining
        $\beta_{j}$ to $\beta_{j+1}$ so that the regular neighborhood
        of the union
        $\tau_j\cup\beta_j\cup\beta_{j+1}$ 
        is a pair of pants whose third boundary component is $\gamma_j$.
        Furthermore, $\tau_j$'s are pairwise disjoint and $\tau_j$ is disjoint
        from $\beta_{k}$ for $k<j$.
    \end{enumerate}
  \end{lemma}
 
  \begin{proof}
    By Proposition \ref{prop:span}, either there is a compact subsurface
    $K$ with totally geodesic boundary such that $\mathcal{R}(K)$ contains
    infinitely many equivalence classes, or the curves in \(\Ca\)
    accumulate somewhere. 
    
    In the first case, we can find a sequence $R_i\in \mathcal{R}(K)$ which
    are all parallel and pairwise not equivalent. Let $\alpha$ be a
    boundary curve of $K$ containing a vertical side of each $R_i$. Up to
    passing to a subsequence, we can assume that there is an oriented arc
    $[p,q] \subset \alpha$ and a sequence of points $v_i\in R_i\cap [p,q]$
    converging monotonically to $q$. For every $i$ we can moreover find a
    curve $\delta_i\in \C$ which is a boundary component of some surface
    $F^{l_i}$ and intersecting $K$ in a collection of arcs that includes a
    horizontal arc of $R_i$. This boundary curve is in $\mathcal C$ since
    the collection is closed under concatenation. In this way all
    $\delta_i$'s are pairwise disjoint since the $F^{l_i}$'s are nested.

    After a subsequence, we can assume that the intersections
    $R_i\cap\alpha$ converge to a point $q\in\alpha$ and that the
    $\delta_i$'s do not contain $q$. After a further subsequence there will
    be an oriented arc $[p,q]\subset\alpha$ that contains a vertical side
    of each $R_i$ and the convergence to $q$ is monotonic along $[p,q]$.
    Further, any arc of any $\delta_i\cap K$ that intersects $[p,q]$ is
    contained in some rectangle. Inductively, we can construct a
    subsequence $\delta_{i_j}$ such that
    \begin{enumerate}[(i)]
     \item there is a strand of $\delta_{i_{j+1}}$ below (i.e.\ closer
       to $q$) all strands of $\delta_{i_j}$, and
     \item between the lowest strand of $\delta_{i_j}$ and the first
       strand of $\delta_{i_{j+1}}$ below $\delta_{i_j}$
        there are at least two other rectangles.
    \end{enumerate}
    Let $\tau_j$ be the subarc of $[p,q]$ connecting the two strands of
    $\delta_{i_j}$ and $\delta_{i_{j+1}}$ in (ii). Let $\gamma_{j}$ be the
    third boundary component of a small regular neighborhood of
    $\delta_{i_j}\cup\tau_{j}\cup\delta_{i_{j+1}}$. Since $\tau_{j}$
    intersects another $\delta_{k}$, $\gamma_{j}$ is essential. Moreover,
    if $\gamma_{j}\in\mathcal C$, then the rectangles associated with the
    two strands of $\delta_k$'s would be equivalent, which is a
    contradiction. Now rename $\beta_j:=\delta_{i_j}$. 
    
    If the curves in \(\Ca\) accumulate somewhere, we can find a curve
    $\alpha$ intersecting infinitely many $\delta_i\subset \Ca$. As before,
    up to passing to a subsequence, we can find an arc $[p,q] \subset
    \alpha$, with some $v_i \in \delta_i\cap [p,q]$ converging to $q$
    monotonically. We then define $\delta_{i_j}$, $\gamma_j$ and $\tau_j$
    as above. Since each $\gamma_j$ intersects a curve in $\Ca$,
    $\gamma_j\notin\C$.

    Note that in both cases there may be strands of $\beta_k$ that
    intersect $\tau_j$ when $k>j+1$. \qedhere

  \end{proof}
  
  \begin{figure}[htp!]
  \begin{center}
  \begin{overpic}{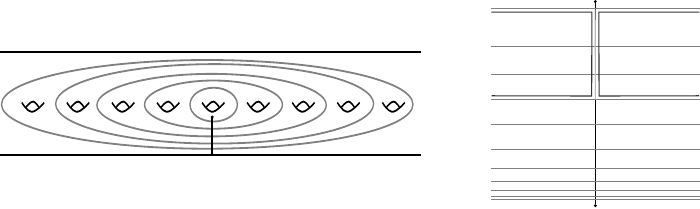}
  \put(29.5,5){$a$}
  \put(84,31){$p$}
  \put(84,-2){$q$}
  \put(101,28){$\beta_{1}$}
  \put(101,14){$\beta_{2}$}
  \put(95,26){$\gamma_1$}
  \put(86,24){$\tau_1$}
  \end{overpic}
  \vspace{.2cm}
    \caption{On the left-hand side, the arc $a=[p,q]$ for the non-good
    collection of curves of Figure \ref{fig:notspanning}; on the right-hand
    side, the general situation in the proof of Lemma \ref{lem:not-good}.}
  \end{center}
  \end{figure}

  \begin{remark}\label{rmk:tau^n}
    In the situation of the previous lemma, and given a homeomorphism $f$,
    fix an index \(j\). Look at \(f^n(\alpha), f^n(\beta_{j})\) and
    \(f^n(\beta_{j+1})\) and consider an ambient isotopy which sends them
    to their geodesic representatives. The image through this isotopy of
    \(f^n(\tau_j)\) is a geodesic arc of $f^n(\alpha)^*$ that we denote by
    $\tau_j^n$. In particular, $f^n(\gamma_j)^*$ is homotopic to
    $f^n(\beta_{j})^* \cup \tau_j^n \cup f^n(\beta_{j+1})^*$.
  \end{remark}
  
  \section{Diagonal closure of curves and lines}
  
  In this section we discuss two subsurfaces associated to a finite
  collection of curves and lines: the subsurface spanned by the collection
  and the diagonal closure of the collection.
  
  Given a finite collection of lines and curves $L$, the subsurface
  \emph{spanned} by $L$, denoted  $\langle L \rangle$, is the subsurface
  obtained by putting the curves and lines in minimal position, taking a
  regular neighborhood of their union, and then filling in all
  complementary disks with at most one puncture.
  
  Given two rays $r_1$ and $r_2$, we say that they \emph{cobound a
  half-strip} if there is a proper embedding $\varphi:\R_{\geq 0}\times
  [0,1]\to S$ such that $\varphi(\R_{\geq 0}\times \{0\})=r_1$ and
  $\varphi(\R_{\geq 0}\times \{1\})=r_2$. We say that two lines $\ell_1$
  and $\ell_2$ are \emph{asymptotic} in some direction if they contain rays
  $r_1\subset \ell_1$ and $r_2\subset \ell_2$ which cobound a half-strip.
  
  \begin{figure}[htp!]
  \begin{center}
  \begin{overpic}{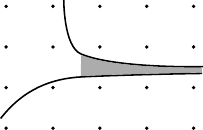}
  \put(35,50){$\ell_1$}
  \put(15,10){$\ell_2$}
  \end{overpic}
  \caption{Two asymptotic lines with the half-strip cobounded by two
    rays.}\label{fig:asymptotic}
  \end{center}
  \end{figure}
  
  Let $L$ be a finite collection of curves or proper lines such that any
  two elements of \(L\) intersect only finitely many times, and no line in
  $L$ ends in an isolated planar end. We denote by $\diagclosure{L}$ the
  subsurface of $S$ obtained as follows. First put the curves and lines in
  $L$ in minimal position. Next form the subsurface $\langle L \rangle$
  spanned by $L$. Then, for every two (not necessarily distinct) lines and
  every direction in which they are asymptotic, adjoin to $\langle L
  \rangle$ the corresponding half-strip cobounded by the asymptotic rays.
  Finally fill in any additional complementary disks with at most one
  puncture. We call $\diagclosure{L}$ the \emph{diagonal closure} of $L$.
  The following lemma gives a characterization of $\diagclosure{L}$. 

  \begin{lemma}\label{lem:lines-outside-surfaces}
    Let $L=\{\ell_1,\dots,\ell_N\}$ be a finite collection of curves and
    lines such that no element ends in an isolated planar end of $S$, and
    any two elements of \(L\) intersect only finitely many times. Then for
    any line $\ell$ that does not end in an isolated planar end, the
    following statements are equivalent.
   \begin{enumerate}
     \item $\ell$ can be homotoped into $\diagclosure{L}$.
     \item For any curve $\alpha$, if $i(\alpha,\ell)\neq 0$, then
       $i(\alpha,\ell_j)\neq 0$ for some $j$.
     \item For every compact subsurface $K$, $\ell$ can be homotoped to
       $\ell'$ such that $\ell'\cap K\subset \langle L\rangle$.
   \end{enumerate}
   Moreover there is a finite union of curves $\mu$ such that if $\ell$ is
    a line in $\diagclosure{L}$, then $i(\mu,\ell)\neq 0$.
  \end{lemma}
  
\begin{figure}[htp!]

  \begin{center}
  \begin{overpic}{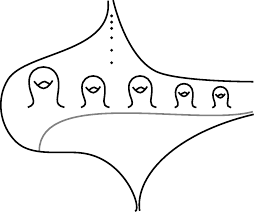}
  \put(70,56){$F$}
  \put(35,31){$a$}
  \put(50,22){$C$}
  \end{overpic}
  \caption{The situation in the proof of Lemma
    \ref{lem:lines-outside-surfaces}.}
    \label{fig:contractible}
  \end{center}
  \end{figure}

  \begin{proof} We first prove the equivalence of the three conditions.
  
    [(1) $\Rightarrow$ (3)] We may assume $\ell\subset \langle
    L\rangle$. Since $\diagclosure{L}$ is formed from $\langle L
    \rangle$ by adjoining finitely many half-strips and disks, we can first
    homotope $\ell$ away from all the adjoined disks. For any compact
    subsurface $K$, there is a homotopy that pushes all the half-strips off
    of $K$. Let $\ell'$ be the image of $\ell$ under this homotopy. Then
    $\ell' \cap K \subset \langle L \rangle$.

    [(3) $\Rightarrow$(2)] Let $\alpha$ be a curve with $i(\alpha,\ell)\neq
    0$. Let $K$ be an annular neighborhood about $\alpha$ and $\ell'$ a
    line homotopic to $\ell$ so that $\ell'\cap K\subset \langle L\rangle$.
    Since $\ell'$ intersects $\alpha$ essentially, $\alpha$ crosses
    $\langle L\rangle$ essentially. Therefore, there is some $j$ such that
    $i(\alpha,\ell_j)\neq 0$.

    [(2) $\Rightarrow$ (1)] We prove the contrapositive: if $\ell$ cannot
    be homotoped into $\diagclosure{L}$, then there is some curve $\alpha$
    that intersects $\ell$ but disjoint from all $\ell_j$'s. 

    Assume all lines are geodesic and $\ell$ is in minimal position with
    respect to $\diagclosure{L}$. Note that $\diagclosure{L}$ is the union
    of a compact surface with finitely many half-strips. In particular the
    boundary of $\diagclosure{L}$ is a finite union of circles and lines.
    The same therefore holds for $\Sigma:=\overline{S\ssm
    \diagclosure{L}}$. Moreover, $\Sigma$ contains no pairs of boundary
    rays cobounding a half-strip in $S$, since any such half-strip would be
    included in the construction of $\diagclosure{L}$. 
    
    Let $a$ be a component of $\ell\cap\Sigma$ and $F$ the component of
    $\Sigma$ containing $a$. If $a$ is nonseparating in $F$, we can find a
    curve contained in $F$ intersecting $a$ once and we are done. If $a$ is
    separating and both components of $F\ssm a$ are not contractible,
    choose two non-contractible curves $\alpha_1$ and $\alpha_2$ (which may
    be homotopic to a puncture), one in each component of $F\ssm a$, and a
    simple arc $b$ connecting them. The regular neighborhood of
    $\alpha_1\cup b\cup\alpha_2$ is a pair of pants whose third boundary
    component gives the required curve $\alpha$. Thus we may assume that
    $a$ is separating and at least one component $C$ of $F\ssm a$ is
    contractible. By the classification of contractible surfaces with
    boundary, $C\cup a$ is a closed disk with points removed from its
    boundary. Since $\partial \Sigma$ has only finitely many components,
    $C\cup a$ is in fact a closed disk with finitely many points removed
    from its boundary. As $\ell$ cannot be homotoped into
    $\diagclosure{L}$, one such point is not an end of $a$. It follows that
    $\Sigma$ contains a half-strip cobounded by two boundary rays, a
    contradiction (see Figure \ref{fig:contractible}). This finishes the
    proof of this case.

    Finally, we prove the last statement. Write 
    \[\diagclosure{L} = K\sqcup S_1\sqcup\dots\sqcup S_k,\]
    where \(K\) is a compact surface and the $S_i$ are pairwise disjoint
    half-strips. For every $j=1,\dots,k$, choose a ray $r_j$ going through
    $S_j$. As no line ends in a puncture, there is a sequence of essential
    separating curves, pairwise not homotopic, converging to the end of
    $r_j$ contained in $S_j$. Choose one such curve $\mu_j$ sufficiently
    far out so that $\mu_j$ intersects the ray $r_j$ and is disjoint from
    $K$. Let $$\mu=\bigcup_{j=1}^k\mu_j.$$ If $\ell$ is a line in
    $\diagclosure{L}$, then by properness $\ell$ has to exit $K$ along some
    ray $r_j$, hence it intersects the corresponding $\mu_j$. This shows
    $i(\mu,\ell) \ne 0$ as desired.\qedhere

  \end{proof}
 
\begin{figure}[htp!]
  \begin{center}
  \begin{overpic}[width=.5\textwidth]{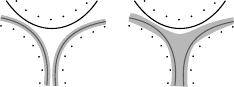}
  \put(20,30){$\ell$}
  \put(75,30){$\ell$}
  \put(18,-5){$\langle L\rangle$}
  \put(72,-5){$\diagclosure{L}$}
  \end{overpic}
  \vspace{.5cm}
  \caption{$L$ is the union of the gray-green lines and $\ell$ is the black
    one. Then $\ell$ is not properly homotopic into $\langle L\rangle$, but
    it is in $\diagclosure{L}$.}\label{fig:diagclosure}
  \end{center}
  \end{figure}

  \begin{remark}
    At first one might think that condition (2) in Lemma
    \ref{lem:lines-outside-surfaces} is equivalent to $\ell$ being in
    $\langle L\rangle$ (up to proper homotopy). The example depicted in
    Figure \ref{fig:diagclosure} shows that this is not the case: condition
    (2) holds, but $\ell$ cannot be properly homotoped into $\langle
    L\rangle$.
  \end{remark}
  
  The reason why we are interested in Lemma
  \ref{lem:lines-outside-surfaces} is the following relationship between
  $\diagclosure{\cdot}$ and limits of concatenations of curves.

  \begin{lemma}\label{lem:diagclosure&limits}
    Suppose $f$ is well tempered. Let $\alpha_1,\dots,\alpha_k$ be curves.
    If $\beta$ is a concatenation of $\alpha_1,\dots,\alpha_k$, then
    \[\L^\pm(\beta) \subset \Big\langle \Big\langle
    {\bigcup_{i=1}^k\L^\pm(\alpha_i)} \Big\rangle\Big\rangle.\]
  \end{lemma}
  
  \begin{proof}
    
    Let $\ell\in\L^+(\beta)$ and $\gamma$ a curve intersecting $\ell$. By
    Lemma \ref{lem:tempered-intersection} there are infinitely many $n\geq 0$,
    such that $i(\gamma,f^n(\beta))\neq 0$. Since $\beta$ is a finite
    concatenation of the $\alpha_i$'s, there is some $j\in\{1,\dots,n\}$
    and infinitely many $n\geq 0$ such that $i(f^n(\alpha_j),\gamma)\neq
    0$. Applying Lemma \ref{lem:tempered-intersection} again, we obtain some
    $\ell'\in\L^+(\alpha_j)$ with $i(\ell', \gamma)\neq 0$. Since $\gamma$
    was arbitrary, Lemma \ref{lem:lines-outside-surfaces} implies that
    $\ell\subset \diagclosure{\bigcup_{i=1}^k\L^\pm(\alpha_i)}$ up to
    proper homotopy. The same argument applies to $\L^-(\beta)$. \qedhere
  \end{proof}

  We end this section with a variant of Lemma
  \ref{lem:lines-outside-surfaces}.

  \begin{lemma}\label{obstruction} 
    Let $A \subset B$ be two finite collections of curves and lines such
    that any two elements intersect only finitely many times. Then there
    exist finitely many curves $\delta_1,\cdots,\delta_s$ disjoint from
    lines in $A$, with the following property: if a line $\ell\subset
    \diagclosure{B}$ cannot be homotoped into $\diagclosure{A}$, then
    $i(\ell,\delta_i)\neq 0$ for some $i$.
  \end{lemma}
  
  \begin{proof} 

    As in the proof of Lemma \ref{lem:lines-outside-surfaces}, we can write
    $\diagclosure{B}$ as the union of a finite-type subsurface $K$ with
    compact boundary and pairwise disjoint half-strips $S_1, \dots, S_k$.
    Since $A\subset B$, we can homotope $\diagclosure{A}$ so that it is the
    union of a finite type subsurface $K' \subset K$ with compact boundary,
    together with some of the half-strips $S_i$; up to renumbering, say
    $S_1,\dots, S_l$.
  
    For each $l< j\leq k$, choose a curve $\delta_{j-l}$ associated to the
    half-strip $S_j$ as in the proof of Lemma
    \ref{lem:lines-outside-surfaces}.  We append to this collection the set
    of boundary curves $\delta_{k-l+1},\dots, \delta_s$ of $K'$.
  
    Now suppose a line $\ell$ in $\diagclosure{B}$ cannot be homotoped into
    $\diagclosure{A}$. If $\ell$ goes through a half-strip not contained in
    $\diagclosure{A}$, then $\ell$ intersects one of the curves $\delta_i$
    with $i\leq k-l$. Otherwise, $\ell\ssm K$ is contained in
    $\diagclosure{A}$, so $\ell\cap K$ cannot be contained in $K'$ and
    hence it intersects some $\delta_i$ with $i>k-l$. \qedhere

  \end{proof}

\section{Well-tempered maps: canonical decomposition}

  \label{sec:decomposition}

  By a \emph{decomposition} of $S$ we mean a collection $\{X_i\}$ of
  pairwise disjoint subsurfaces of $S$ (each of which may be disconnected),
  such that we can find representatives $\tilde{X_i}$ with $S = \bigcup
  \tilde{X_i}$ and such that the only overlap between $\tilde{X_i}$ and
  $\tilde{X_j}$, $i \ne j$, are common boundary components. When $f$ is a
  homeomorphism of $S$, a decomposition $\{X_i\}$ of $S$ is called
  \emph{$f$-invariant} if $f(X_i)$ is isotopic to $X_i$ for all $i$. For a
  subsurface $X \subset S$, we say that $f$ \emph{returns to} $X$ if there
  is a power $k>0$ such that $f^k(X)$ is isotopic to $X$. Any such $k>0$ is
  called a \emph{return time} of $f$ to $X$. Note that all return times of
  $f$ to $X$ are multiples of the smallest return time.   
  
  Our main proposition in this section is the following. 
  
  \begin{prop} \label{prop:decomposition}

    Suppose $f$ is well tempered. Then there is an $f$-invariant
    decomposition of $S$ into three subsurfaces, $\Sper$, $\Sinf$, and
    $S_0$, with the following properties.
    \begin{itemize}
      \item A curve in $S$ is periodic if and only if it can be
        homotoped into $\Sper$. 
      \item A curve in $S$ is wandering if and only if it can be
        homotoped into $\Sinf$. 
      \item $S_0$ contains no essential, nonperipheral curves.  
    \end{itemize}
    We will call $\Sper$, $\Sinf$, and $S_0$ the \emph{canonical
    decomposition} of $f$ on $S$.

  \end{prop}
  
\begin{figure}[htp!]
  \begin{center}
  \includegraphics[scale=.4]{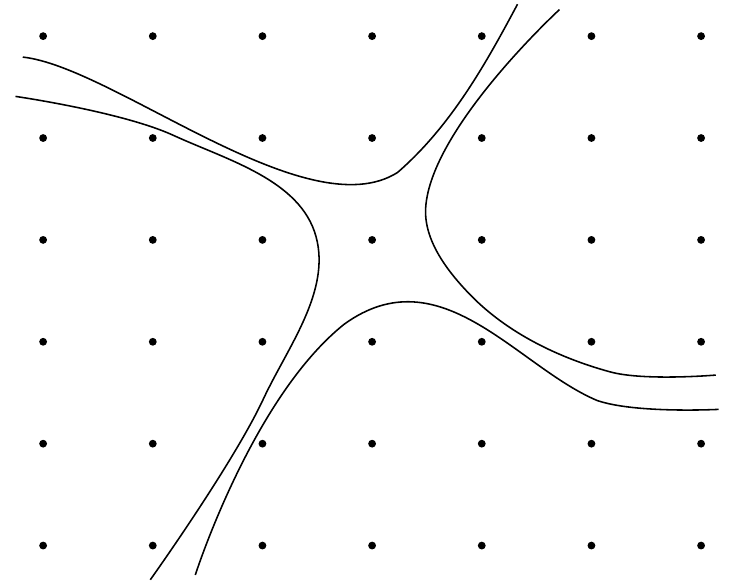}
  \caption{The four invariant lines giving the decomposition of Proposition
    \ref{prop:decomposition} for the map $(x,y)\mapsto (2x+y,x+y)$ on
    $S=\R^2\ssm\Z^2$.}\label{fig:2111}
  \end{center}
  \end{figure}

  An explicit example of a well-tempered map is the map induced by the matrix
  $\left(
  \begin{array}{cc}
  2 & 1\\
  1 & 1
  \end{array}
  \right)$ on $S=\R^2\ssm\Z^2$. In this case, the decomposition is given by
  four invariant lines, which cobound a subsurface $X=S_0$ homeomorphic to
  a once-punctured disk with four points removed from the boundary. The
  rest of the surface is $\Sinf$, consisting of four invariant components
  (see Figure \ref{fig:2111}).

  For any surface $S$ and any well-tempered map $f$, the decomposition is
  defined as follows. Let $\Cper$ be the collection of $f$-periodic
  curves, and $\Cinf$ the collection of $f$-wandering curves. Let $\Sper$
  be the subsurface spanned by $\Cper$ (whose existence is proved in Lemma
  \ref{lem:SperSubs}), $\Sinf$ the subsurface spanned by $\Cinf$ (whose
  existence is proved in Lemma \ref{lem:SinftySubs}). 
 
  Fix a hyperbolic metric on $S$. We choose an almost geodesic
  representative for the components of $\Sper$ and $\Sinf$ as follows. For
  each non-annular component $X$, let $X^*$ be the almost geodesic
  representative of $X$. Recall this means that we first take the geodesic
  representative $X^\circ$ of the interior of $X$ and then remove a small
  regular neighborhood of a boundary component of $X^\circ$ that is
  homotopic to another boundary component of $X^\circ$. We also assume that
  if two components share any boundary component, their almost geodesic
  representatives do as well. By abuse of notation, we continue to use
  $\Sper$ and $\Sinf$ for the union of the representatives of their
  components. Define $$S_0:=\overline{S\ssm(\Sper\cup\Sinf)},$$ whose
  isotopy class is independent of the choice of the hyperbolic metric on
  $S$.

  As just seen, a fundamental step in proving Proposition
  \ref{prop:decomposition} is to show that, for a well-tempered map, $\Cper$
  and $\Cinf$ span subsurfaces $\Sper$ and $\Sinf$. We then need to show
  that $S_0$ has no interesting topology. 

  All three proofs proceed by contradiction. Namely, we will show that if
  the required statement is false, then Lemma \ref{lem:not-good} can be
  used to produce an essential curve $\alpha$ with infinite limit set
  $\L(\alpha)$. The arguments for $\Sper$ and $\Sinf$ are independent. The
  proof for $S_0$ is more technical: it relies on the fact that $S_0$ is
  already a subsurface, being the complement of two subsurfaces, and on the
  observation that no curve in $S_0$ can be periodic or wandering, hence
  any curve in $S_0$  must limit onto a non-empty collection of lines.  

  Let us now prove that $\Sinf$ and $\Sper$ are subsurfaces.

  \begin{lemma}\label{lem:SinftySubs}
    If \(f\) is well tempered, then $\Cinf$ spans a subsurface $\Sinf$. Every
    curve in $\Sinf$ is wandering, and if $\gamma$ is any curve with
    $\L(\gamma) \ne \emptyset$, then $\L(\gamma)$ does not intersect
    $S_\infty$ essentially.
  \end{lemma}

  \begin{proof}

    Suppose there is no subsurface spanned by $\Cinf$. By Lemma
    \ref{lem:periodic&wandering}, \(\Cinf\) is closed under concatenation,
    so we can apply Lemma \ref{lem:not-good} and find curves $\alpha$,
    $\beta_j$ and $\gamma_j$ and pairwise disjoint arcs $\tau_j \subset
    \alpha$ as in the lemma. Recall that $\beta_j \in \Cinf$, but $\alpha,
    \gamma_j \notin \Cinf$. In particular, \(\L(\gamma_j)\) and
    $\L(\alpha)$ are nonempty. Since $\gamma_j$ is homotopic to $\beta_j
    \cup \tau_j \cup \beta_{j+1}$ (and none of these is periodic, by
    assumption and Corollary \ref{cor:not-wandering}), Lemma
    \ref{lem:diagclosure&limits} gives
    \[\L(\gamma_j) \subset \diagclosure{\L(\alpha)\cup\L(\beta_{j}) \cup
    \L(\beta_{j+1})}=\diagclosure{\L(\alpha)}, \] where the last equality
    uses the fact that the $\beta_i$ are wandering, hence
    $\L(\beta_i)=\emptyset$. Since $\L(\alpha)$ is a finite collection of
    lines with finite pairwise intersections, Lemma
    \ref{lem:lines-outside-surfaces} provides a finite union of curves
    $\mu$ such that every line in $\diagclosure{\L(\alpha)}$ intersects
    $\mu$. Hence there is a curve $\delta\subset \mu$ that intersects
    infinitely many $\L(\gamma_j)$. By Lemma \ref{lem:duality},
    \(\L(\delta)\) intersects infinitely many \(\gamma_j\). Since
    $\L(\delta)$ is finite, some line \(\ell\in\L(\delta)\) intersects
    infinitely many \(\gamma_j\). 

    If \(i(\ell,\beta_{j})\neq 0\), then $i(\delta,\L(\beta_{j}))\neq 0$,
    which is a contradiction since $\beta_{j}$ is wandering; the same
    applies to $\beta_{j+1}$. Therefore \(\ell\) must intersect infinitely
    many of the arcs \(\tau_j\), and hence intersects \(\alpha\) infinitely
    many times. This is impossible since $\ell$ is a proper line. This
    contradiction shows that $\C_\infty$ spans a subsurface $\Sinf$.

    For the second statement, if $\gamma$ is a curve such that $\L(\gamma)$
    intersects $\Sinf$ essentially, then there exists $\beta\subset\Sinf$
    such that $i(\beta,\L(\gamma)) \ne 0$. But every curve in $S_\infty$ is
    a concatenation of wandering curves, hence $\beta$ is wandering by
    Lemma \ref{lem:periodic&wandering}. This contradicts Lemma
    \ref{lem:not-wandering}.\qedhere

  \end{proof}

  We also need the following generalization.

  \begin{lemma} \label{SL}
    Suppose \(f\) is well tempered. Let $L^+$ and $L^-$ be two finite
    collections of lines with finite pairwise intersections, and such that
    each of $L^+$ and $L^-$ has finite intersection with every $\mathcal
    L^\pm(\alpha)$. Define $$\mathcal C_L:=\{\beta\mid \mathcal
    L^\pm(\beta)\subseteq \diagclosure{L^\pm}\}.$$ Then $\mathcal C_L$
    spans a subsurface $S_L$ of $S$.
  \end{lemma}

  \begin{proof}
    Note that $\mathcal C_L$ contains all wandering curves, and is closed
    under concatenation by Lemma \ref{lem:diagclosure&limits}. 
    
    Assume that \(\mathcal{C}_L\) does not span a subsurface. Then we have
    curves $\alpha,\beta_j,\gamma_j$ and arcs $\tau_j \subset \alpha$ as in
    Lemma \ref{lem:not-good}, where $\beta_j \in \mathcal{C}_L$; $\alpha,
    \gamma_j \notin \mathcal{C}_L$; and $\gamma_j$ is homotopic to $\beta_j
    \cup \tau_j \cup \beta_{j+1}$. In this situation, 
    \[ \L^\pm(\gamma_j) \subset
    \diagclosure{\L^\pm(\alpha)\cup\L^\pm(\beta_{j}) \cup
    \L^\pm(\beta_{j+1})}\subseteq\diagclosure{\L^\pm(\alpha)\cup L^\pm}.\]
    Let $\{\delta_l^+\}_l$ and $\{\delta_m^-\}_m$ be the curves from Lemma
    \ref{obstruction} for the collections $L^\pm\subset L^\pm\cup \mathcal
    L^\pm(\alpha)$. After passing to a subsequence and possibly replacing
    $f$ by $f^{-1}$ and switching $+$ and $-$, we may assume that for every
    \(j\) there is at least one line in $\mathcal L^+(\gamma_j)$ that is
    not homotopic into $\diagclosure{L^+}$, hence intersects one of the
    $\{\delta_l^+\}_l$. Passing to a further subsequence, we may assume
    that the same curve $\delta \in \{ \delta^+_l\}$ intersects infinitely
    many $\mathcal L^+(\gamma_j)$. By construction, $\delta$ does not
    intersect any of $\mathcal L^+(\beta_j)$. It follows by Lemma
    \ref{lem:duality} that $\mathcal L^-(\delta)$ must intersect $\alpha$
    infinitely many times, which is impossible since the map is well tempered. \qedhere
  \end{proof}
  
  \begin{lemma}\label{lem:SperSubs}
    If \(f\) is well tempered, then $\Cper$ spans a subsurface $\Sper$, and
    every curve in $\Sper$ is periodic.
  \end{lemma}

  \begin{proof}
    By Lemma \ref{lem:periodic&wandering}, $\Cper$ is closed under
    concatenation. Suppose, for contradiction, there is no subsurface
    spanned by $\Cper$. Then we can find curves $\alpha$, $\beta_j$ and
    $\gamma_j$ and arcs $\tau_j \subset \alpha$ as in Lemma
    \ref{lem:not-good}. Recall that $\beta_i \in \Cper$, but $\alpha,
    \gamma_j \notin \Cper$. Moreover, since $\alpha$ intersects periodic
    curves, it is not wandering either. Thus \(\L(\alpha)\) is non-empty
    and consists of lines. 

    Each line in \(\L(\alpha)\) is the limit of a uniformly bounded number
    of strands of $f^{n_i}(\alpha)^*$. Let $N$ be the cardinality of
    $\L(\alpha)$, counted with multiplicity, and set $J=N+1$. Consider the
    finite collection of curves $\gamma_1,\ldots,\gamma_J$. By replacing
    $f$ by a power, we may assume $f(\beta_{j}) = \beta_{j}$ for all $j \le
    J+1$. 

    Recall that $\gamma_j$ is homotopic to \(\beta_{j}\cup
    \tau_j\cup\beta_{j+1}\). For each $n$, let $\tau_j^n$ be the geodesic
    subarc of $f^n(\alpha)^*$ which together with
    $f^n(\beta_{j})^*=\beta_{j}$ and $f^n(\beta_{j+1})^*=\beta_{j+1}$ is
    homotopic to $f^n(\gamma_j)^*$ (see Remark \ref{rmk:tau^n}). Since
    $\gamma_j$ is not periodic, the arcs $\tau_j^n$ cannot be periodic. As
    there are only finitely many isotopy classes of arcs of bounded length
    between \(\beta_{j}\) and \(\beta_{j+1}\), the length of $\tau_j^n$
    cannot be bounded as $n \to \infty$. 

    Moreover, the angles of intersections of \(\tau_j^n\) with $\beta_{j}$
    and $\beta_{j+1}$ are bounded away from zero; otherwise some sequence
    of iterates of \(\alpha\) would wrap around $\beta_{j}$ or
    $\beta_{j+1}$ more and more times, contradicting the tempered
    condition. Thus any convergent subsequence of \(\tau_j^n\) contains a
    ray starting at \(\beta_{j}\) and a ray starting at \(\beta_{j+1}\). By
    construction, each such ray is contained in a line of \(\L(\alpha)\). 

    Choose a subsequence $f^{n_{i}}$ so that each (oriented)
    \(\tau_j^{n_i}\) converges. Then the union of all these limit sets
    contains at least $2J=2N+2$ rays. If a collection of $k$ such limiting
    rays is nested within the same line in \(\L(\alpha)\), then the
    multiplicity of that line is at least $k$. Taking into account both
    orientations of each line, we see that there can be at most $2N$ rays,
    a contradiction. This shows the existence of $\Sper$. \qedhere

  \end{proof}

  \begin{lemma} \label{lem:S0}
  
    If $f$ is well tempered, then $\Sper$ and $\Sinf$ are disjoint
    subsurfaces. 

  \end{lemma}

  \begin{proof}

    Since they are subsurfaces spanned by curves, if they are not disjoint,
    by Corollary \ref{cor:spanned-disjoint} there are curves $\alpha$ in
    $\Sper$ and $\beta$ in \(S_\infty\) with $i(\alpha,\beta) > 0$,
    contradicting  Corollary \ref{cor:not-wandering}.\qedhere

  \end{proof}
 
  \begin{lemma}\label{prop:no-ess-nonper}

    If \(f\) is well tempered, $S_0 :=\overline{S\ssm(\Sper\cup\Sinf)}$ has
    no essential, nonperipheral curves. 
 
  \end{lemma}
  \begin{proof}
    Suppose, by contradiction, that there is an essential, nonperipheral
    curve $\delta$ in $S_0$. Then $\delta$ is neither wandering nor
    periodic, so $L:=\L(\delta)$ is a nonempty collection of lines. Let
    $L^\pm:=\L^\pm(\delta)$ and let $S_L$ be the subsurface from Lemma
    \ref{SL} spanned by the collection of curves $\{\gamma\}$ with
    $\L^\pm(\gamma)\subset \diagclosure{L^\pm}$. Put $S_L$ and $S_0$ in
    minimal position and let $\tilde{S}_L=S_0\cap S_L$. Note that
    $S_\infty\subseteq S_L$ and $\Sper\cap S_L$ is a compact surface, so
    the boundaries of $S_L$ and $S_0$ intersect in only finitely many
    points.
   
    Up to replacing $f$ by its inverse, we may assume \(L^+\) is nonempty.
    By Lemma \ref{lem:lines-outside-surfaces}, there exists a finite
    collection of curves $\{\delta_1,\ldots,\delta_k\}$ intersecting all
    lines in \(\diagclosure{L^+}\). Then Lemma \ref{lem:duality} implies
    that for every curve \(\gamma\) in \(S_L\), either
    \[ 
    \text{(a)}\,\, \L^+(\gamma)=\emptyset, \quad \text{or} \quad \text{(b)}
    \,\, i(\gamma, \L^-(\delta_i))\neq 0 \,\, \text{for some $i$}.
    \]

    \emph{Claim:} There exists an essential, nonperipheral curve in
    $\tilde{S}_L$ satisfying (a).
    
    Indeed, otherwise the finite collection of curves and lines in
    $\bigcup_i \L^-(\delta_i)$ fills every curve in $\tilde{S}_L$. Since
    there are only finitely many intersection points between these curves
    and lines, it follows that each component of $\tilde{S}_L$ has finite
    type. Furthermore, since \(\tilde{S}_L\) is \(f\)-invariant, if one of
    its components is not periodic, then its orbit is a subsurface of
    \(\tilde{S}_L\). In particular the orbit of the component cannot
    accumulate anywhere; that is, it is wandering. Now consider the
    component that contains $\delta$. If it is wandering, \(\delta\) is
    wandering, which is impossible. Hence this component must be periodic.
    But then the first return map is a tempered mapping class of a
    finite-type surface; in particular, it is periodic and so \(\delta\) is
    periodic as well, which is again impossible. This proves the claim.
      
    By the claim, there exists an essential, nonperipheral curve $\delta'$
    in $S_0$ that is forward wandering. Replace $f$ by $f^{-1}$ and run the
    whole construction on $\delta'$, obtaining the subsurfaces $S_{L'}$ and
    $\tilde{S}_{L'}$, where $L' = L^{\pm}(\delta')$. Now, every curve
    \(\gamma\) in $S_{L'}$ is \(f\)-forward wandering (since
    \(\L^+(\gamma)\subset\diagclosure{\L^+(\delta')}=\emptyset\)). The
    above claim produces another essential, nonperipheral curve in
    $\tilde{S}_{L'}$ which is \(f\)-backward wandering, and thus wandering,
    a contradiction. \qedhere
      
  \end{proof}

  A consequence of the lemma is the following.

  \begin{cor} \label{cor:pants}
    Any negative Euler characteristic component of $S_0$ is a pair of
    pants.
  \end{cor}

  \begin{proof}
    Let $Y$ be a component of $S_0$ with $\chi(Y) <0$.  If $\chi(Y) < -1$,
    then $Y$ contains an essential nonperipheral curve, a contradiction.
    Hence $\chi(Y)=-1$. If $Y$ were not a pair of pants, then it would be a pair
    of pants with some boundary points removed, and such a surface still
    contains an essential nonperipheral curve. This contradiction proves
    that $Y$ is a pair of pants. \qedhere
  \end{proof}

  \subsection{Examples} 

  \label{sec:examples}
 
  We end this section with some examples of tempered but not well tempered
  maps for which the conclusions of Proposition \ref{prop:decomposition}
  are false.

\begin{figure}[htp!]
  \begin{center}
  \begin{overpic}[scale=1.2]{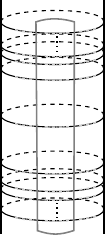}
  \put(46,88){$1$}
  \put(46,8){$-1$}
  \end{overpic}
  \caption{A sketch of the surface in Example \ref{ex:Spernotsurf}, with
    a periodic curve (in gray) and the Cantor sets in light
    gray.}\label{fig:cylinder}
  \end{center}
  \end{figure}

  \begin{example} \label{ex:Spernotsurf}
    This example shows the existence of a tempered map $f$ such that
    $\L(\alpha)$ is locally finite for every curve $\alpha$, but the
    collection of periodic curves does not span a subsurface. 

    Consider the cylinder $\Sigma:=S^1\times \R$ and fix a Cantor set
    $K\subset S^1$. Let $S$ be the surface obtained as
    \[S:=\Sigma\ssm\bigcup_{n=-\infty}^\infty K\times\{t_n\},\]
    where $t_{\pm\infty}=\pm 1$, $t_0=0$, $t_n=-1+\frac{1}{2^{-n}}$ for
    $n<0$ and $t_n=1-\frac{1}{2^{n}}$ for $n>0$ (see Figure
    \ref{fig:cylinder}).

    Fix a monotone non-decreasing homeomorphism $s:\R\to\R$ such that
    \begin{align*}
    s(t_{\pm\infty})&=t_{\pm\infty}\\
    s(t_k)&=t_{k+1} \hspace{.5cm} \forall k\in\Z.
    \end{align*}
    Then the map $\text{id}\times s:\Sigma\to\Sigma$ induces a tempered
    homeomorphism $f$ of $S$, which has no wandering curves. The periodic
    curves are also fixed curves, which are those that can be homotoped to
    be vertical between height $-1$ and height $1$. In particular $\Cper$
    does not span a subsurface.

  \end{example}

\begin{figure}[htp!]
  \begin{center}
  \includegraphics[width=.5\textwidth]{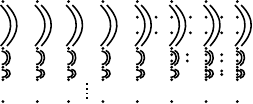}
  \caption{The modified surface in Example
    \ref{ex:Sinfnotsurf}.}\label{fig:connectedSinf}
  \end{center}
  \end{figure}

  \begin{example}\label{ex:Sinfnotsurf}

    This example shows that there exists a tempered map $f$ such that
    $\L(\alpha)$ is locally finite for every curve $\alpha$, but $\Cinf$
    does not span a subsurface. It also shows that such an example can be
    constructed so that a single component of $\Sinf$ is not a subsurface.
 
    In the plane, for any $k\in\N$, let $\ell_k$ be the boundary of
    \[ [k+1,\infty) \times
    \left[\frac{1}{k+1}+\varepsilon_k,\frac{1}{k}-\varepsilon_k\right], \]
    where $\varepsilon_k=\frac{1}{10k^2}$ (see Figure \ref{fig:lightblue}
    for a schematic picture). Let $S$ be the plane punctured at the points
    $\left(m,\frac{1}{n}\right)$, for $n\in\N$ and $m\in\Z$, at the points
    $(0,m)$, for $m\in\Z$, and for every line $\ell_k$, $k\in\N$, at a
    sequence of points on the line not accumulating anywhere.

   \begin{figure}[htp!]
   \begin{center}
   \begin{overpic}[width=.8\textwidth]{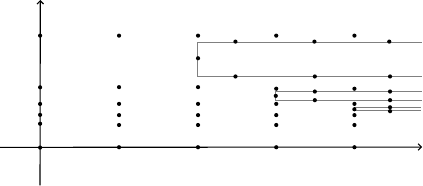}
   \put(1,5){$(0,0)$}
   \put(25,5){$(1,0)$}
   \put(43,5){$(2,0)$}
   \put(62,5){$(3,0)$}
   \put(81,5){$(4,0)$}
   \put(1,34){$(0,1)$}
   \put(-2,23){$(0,1/2)$}
   \put(-2,18){$(0,1/3)$}
   \put(101,33){$\ell_1$}
   \put(101,21){$\ell_2$}
   \put(101,16){$\ell_3$}
   \end{overpic}
   \caption{The lines $\ell_k$ in the construction of Example
     \ref{ex:Sinfnotsurf}.}\label{fig:lightblue}
   \end{center}
   \end{figure}
   The homeomorphism $f$ is defined as follows: 
   \begin{itemize}
   \item on a small regular neighborhood (of width less than
     \(\varepsilon_k\)) of each horizontal line $y=\frac{1}{n}$ ($n\in \N$)
       and $y=0$, the map is induced by $(x,y)\mapsto (x+1,y)$;
   \item for each $k$, it performs a puncture shift supported on a small
     regular neighborhood of $\ell_k$, tapered to the identity on the rest
       of the surface.
   \end{itemize}

    One can show that $f$ is tempered. Moreover, $F(\Cinf)$ is the disjoint
    union of one strip per line $y=\frac{1}{n}$, $n\in\N$, and of the
    strips $(k,\infty) \times \left(
    \frac{1}{k+1}+\frac{\varepsilon_k}{2},\frac{1}{k}-\frac{\varepsilon_k}{2}\right)$.
    The closure of each component is a subsurface, but the horizontal
    strips accumulate onto $y=0$, so $\Cinf$ does not span a subsurface.
    Note moreover that there are no periodic curves.

    Finally, we can modify the surface by adding handles connecting
    consecutive horizontal strips in a translation-invariant way (see
    Figure \ref{fig:connectedSinf} for a schematic picture). The induced
    map on this new surface is still tempered and has no periodic curves. In
    this case, $F(\Cinf)$ consists of one component per line $\ell_k$ and a
    single nonplanar component. This nonplanar component is not homotopic
    to a subsurface.

  \end{example}

\section{Characterization of translations}

  \label{sec:translation}

  Recall a homeomorphism $f: X \to X$ is called a \emph{translation} if it
  generates an infinite cyclic group acting properly on $X$. In this
  section, we give a characterization of the mapping class of a translation
  in terms of the dynamics of its action on curves. The main application
  will be to a well-tempered map $f$ and its action on the invariant
  subsurface $\Sinf$ in its canonical decomposition. Therefore, in this
  section, we need to work with surfaces with boundary, and we will prove a
  more general version of Theorem \ref{thm:translationintro}. Everything in
  this section is independent of the previous sections.

  Extending the notation for curves, when $\delta$ is an essential arc in a
  hyperbolic surface with totally geodesic boundary, we denote by
  $\delta^*$ the geodesic arc homotopic to $\delta$ which is
  perpendicular to the boundary at both endpoints.

  \begin{lemma}\label{lem:e_pm}

    Let $X$ be a connected surface, possibly with boundary, spanned by its
    curves. Let $f$ be a homeomorphism of $X$ such that every curve in $X$
    is $f$-wandering. Then there are ends $e_\pm$ of $X$ (possibly $e_+ =
    e_-$) such that for every curve $\alpha$ in $X$, $f^{n}(\alpha)^*$
    converges to $e_\pm$ as $n\to \pm\infty$. 

  \end{lemma}

  We call $e_+$ and $e_-$ the $f$-\emph{attracting} and
  $f$-\emph{repelling} ends of $X$, respectively.

  \begin{proof}

    Let $\alpha$ be a nonperipheral curve in $X$. We first show that
    $f^n(\alpha)^*$ converges to an end. Let
    \[\alpha^*=\alpha_0,\alpha_1,\dots\alpha_k=f(\alpha)^*\] be a chain of
    geodesic curves in \(X\) connecting \(\alpha^*\) and \(f(\alpha)^*\),
    i.e.\ \(i(\alpha_i,\alpha_{i+1})\neq 0\) for every \(i\). For any
    compact set $K \subset X$, there exists $n_0$ such that for all $n \ge
    n_0$ and all $i \in \{1,\ldots,k\}$, the curves $f^n(\alpha_i)^*$ are
    disjoint from $K$. Since consecutive ones intersect, they must lie in
    the same connected component of $X \setminus K$. In particular,
    $f^n(\alpha)^*$ and $f^{n+1}(\alpha)^*$ lie in the same complementary
    component of $K$, for all $n\ge n_0$. As this holds for all compact
    sets of $X$, it follows that $f^n(\alpha)^*$ converges to a single end
    $e_+$ as $n\to\infty$. The same argument applied to $f^{-1}$ shows that
    $f^n(\alpha)^*$ converges to an end $e_-$ as $n\to -\infty$.

    Now let $\beta$ be another nonperipheral curve in $X$. Choose a chain
    of curves connecting $\alpha^*$ to $\beta^*$ and repeat the same
    argument as above to deduce that $f^n(\beta)^\ast\to e_\pm$ as
    $n\to\pm\infty$. 

    Finally, if $\gamma$ is a peripheral curve, then it can be expressed
    as a concatenation of nonperipheral curves (since $X$ is spanned by
    its curves), all of which converge under forward (respectively
    backward) iteration to $e_+$ (respectively $e_-$), and hence so does
    $\gamma$.\qedhere

  \end{proof}
 
  \begin{lemma}\label{lem:delta}

    Let $X$ be a connected surface, possibly with boundary, and let
    $e_+,e_-$ be two distinct ends of $X$. Then there is either a curve or
    a finite collection of disjoint arcs in $X$ separating $e_+$ from
    $e_-$.

  \end{lemma}

  \begin{proof}
    Let $DX$ be the double of $X$ along its boundary (and set $DX=X$ if
    $\partial X=\emptyset$). Let $i : X \hookrightarrow DX$ be the
    inclusion map and $i_*: \Ends(X) \to \Ends(DX)$ the induced map. Note
    that $i_*$ is injective: if $e_1 \ne e_2$ are ends of $X$, there is a
    compact subset $K$ separating them, and therefore $i_*(e_1)$ and
    $i_*(e_2)$ are separated by the double of $K$ in $DX$.

    Since $DX$ has no boundary, there is a curve $\hat{\delta}\subset DX$
    separating $i_*(e_+)$ and $i_*(e_-)$. Assume $\hat{\delta}$ is in
    minimal position with respect to $i(\partial X)$ and let
    $\delta:=i^{-1}(\hat{\delta})$. Then $\delta$ is either a curve in $X$
    or a finite collection of pairwise disjoint arcs. Moreover, $\delta$
    separates $e_+$ and $e_-$, because if not, there would exist a line
    $\ell$ in $X$ from $e_+$ to $e_-$ disjoint from $\delta$, and hence a
    line $i(\ell)$ from $i_*(e_+)$ to $i_*(e_-)$ disjoint from
    $\hat{\delta}$, contradicting that $\hat{\delta}$ separates the
    ends.\qedhere
  \end{proof}

  \begin{prop}\label{prop:wandering-arcs}

    Let $X$ be a connected bordered surface spanned by its curves. Let $f$
    be a homeomorphism of $X$ such that every curve is $f$-wandering and
    every boundary line is $f$-periodic. Then every essential arc $\gamma$
    in $X$ is $f$-wandering. Moreover, $\gamma$ converges to the
    $f$-attracting end under forward iteration and to the $f$-repelling
    end under backward iteration. 

  \end{prop}

  \begin{proof}
    Fix a complete hyperbolic structure on $X$ with totally geodesic
    boundary and positive injectivity radius, i.e.\ a positive lower bound
    on the lengths of essential closed curves in $X$. For any arc $\gamma$,
    denote by $\gamma^*$ its geodesic representative which is orthogonal to
    the boundary of $X$. We must show that $f^n(\gamma)^*$ leaves every
    compact set of $X$.   

    It is enough to show the statement for arcs whose endpoints lie on
    noncompact boundary components of $X$. Indeed, suppose first that
    $\gamma$ has endpoints on peripheral curves $\alpha$ and $\beta$,
    possibly $\alpha=\beta$. Then the regular neighborhood of $\alpha \cup
    \beta \cup \gamma$ contains an essential curve. Since every curve under
    forward iteration converges to the same attracting end of $X$, $\gamma$
    is also forward wandering. Similarly for backward wandering. If
    $\gamma$ has one endpoint on a peripheral curve $\alpha$ and the other
    endpoint on a boundary line $\ell$, then a regular neighborhood of
    $\alpha \cup \gamma$ contains an arc $\gamma'$ with both endpoints on
    $\ell$. If $\gamma'$ also converges to the same end of $X$ as curves
    under iteration, then the same is true for $\gamma$.

    Henceforth, we may assume that $\gamma$ is an arc with both endpoints
    lying on noncompact boundary components of $X$. As for curves, denote
    by $\L^\pm(\gamma)$ the accumulation sets of $\gamma_n^*$ as $n \to \pm
    \infty$. The argument of Lemma \ref{lem:duality} applies equally to
    arcs and gives the same duality: \[ i(\L(\gamma), \beta) \ne 0
    \Longleftrightarrow i(\L(\beta), \gamma) \ne 0,\] where $\beta$ is any
    curve or arc in $S$. Since $X$ is spanned by its curves, any geodesic
    segment of $\L(\gamma)$ that meets the interior of $X$ essentially
    intersects some curve. But this is impossible by the duality above,
    since curves are wandering. Hence, if $\L(\gamma) \ne \emptyset$, then
    it is a union of boundary components. It follows that the endpoints of
    $\gamma_n^*$ must escape to infinity along $\partial X$; otherwise
    $\L(\gamma)$ would contain an arc or a ray perpendicular to $\partial
    X$. We now consider the following cases.  

    \emph{Case 1.} Both endpoints of $\gamma$ lie on the same oriented
    boundary line $\ell$ and $f$ preserves the oriented line $\ell$. 

    In this case, let $\alpha$ be the curve obtained by closing up $\gamma$
    by an arc in $\ell$. Also let $\beta$ be the possibly immersed curve
    obtained by joining $\gamma$ to $f(\gamma)$ by two arcs in $\ell$. The
    argument of Lemma \ref{lem:periodic&wandering} applies equally to
    immersed curves, so $f^n(\alpha)^*$ and $f^n(\beta)^*$ both converge to
    the same end of $X$ as $n \to \pm \infty$. 
    
    Suppose for contradiction that $\L^+(\gamma)\neq\emptyset$, i.e.\
    $\gamma$ is not forward wandering. Then there is a compact subset $K
    \subset X$ and a sequence $n_i\to\infty$ such that $f^{n_i}(\gamma)^*
    \cap K \ne \emptyset$ for all $i$. Up to homotopy, we can write
    $f^n(\gamma)^*$ as the concatenation $\tau_n * f^n(\alpha)^* *
    \overline\tau_n$, where $\tau_n$ is a geodesic arc perpendicular to
    $\ell$ at one endpoint  and to $f^n(\alpha)^*$ at the other.

    \medskip \noindent \textbf{Claim 1.} For all sufficiently large $i$,
    the arc \(\tau_{n_i}\) intersects a fixed neighborhood $L$ of \(K\).

    \begin{proof} 
      We first show that the distance between the endpoints of
      $f^{n_i}(\gamma)^*$ tends to 0 as $i\to\infty$. Otherwise, after
      passing to a subsequence, the curve $f^{n_i}(\alpha)^*$ and the
      arc $f^{n_i}(\gamma)^*$ would lie
      uniformly close to each other. Since
      $f^{n_i}(\gamma)^*$ meets $K$ and $\alpha$ is forward wandering,
      this is a contradiction.

      The length of \(\tau_{n_i}\) must therefore tend to infinity.
      Otherwise, after passing to a subsequence, the concatenation
      $\tau_{n_i} * f^{n_i}(\alpha)^* * \overline\tau_{n_i}$ would lie
      uniformly close to the curve $f^{n_i}(\alpha)^*$, and hence would
      also escape to infinity, contrary to $f^{n_i}(\gamma)^*$ meeting $K$. 

      Now, the length of \(f^{n_i}(\alpha)^*\) is bounded below by the
      injectivity radius assumption. Since the endpoints of
      $f^{n_i}(\gamma)^*$ become arbitrarily close, basic hyperbolic
      geometry implies that $f^{n_i}(\gamma)^*$ lies within a uniform
      neighborhood of $\tau_{n_i} * f^{n_i}(\alpha)^* *
      \overline\tau_{n_i}$. Since $f^{n_i}(\gamma)^*$ meets \(K\) but
      $f^{n_i}(\alpha)^*$ escapes to infinity, the arc \(\tau_{n_i}\) must
      intersect a fixed neighborhood $L$ of \(K\) for all sufficiently
      large $i$. \qedhere
    \end{proof}

    By enlarging $L$ if necessary, capping off its inessential boundary
    components of $L$ with punctured disks and tightening its essential
    boundary components, we may assume that $L$ is a totally geodesic
    subsurface of finite type. Each boundary component is then either a
    closed geodesic or a piecewise geodesic loop with pieces alternating
    between segments of $\partial X$ and geodesic arcs perpendicular to
    $\partial X$. Note that for such $L$, if a geodesic arc with endpoints
    outside $L$ intersects $L$, then it cannot be homotoped off of $L$ by a
    homotopy that keeps both endpoints outside $L$.

    \medskip \noindent \textbf{Claim 2.} For all sufficiently large $n$,
    the arc $\tau_n$ is homotopic to $\tau_{n+1}$ keeping the endpoints
    disjoint from $L$.

    \begin{figure}[htp!]
    \begin{center}
    \begin{overpic}[scale=.7]{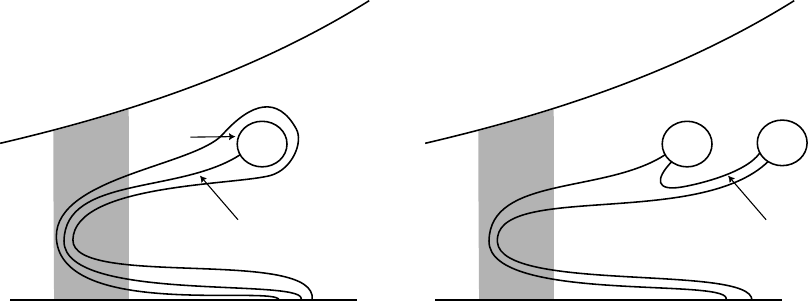}
    \put(30,7){\(\tau_{n}\)}
    \put(30,25){\(f^{n}(\gamma)^*\)}
    \put(11,20){\(f^{n}(\alpha)^*\)}
    \put(3,12){\(L\)}
    \put(20,-5){\(\ell\)}

    \put(70,17){\( \tau_{n}\)}
    \put(80,24){\(f^{n}(\alpha)^*\)}
    \put(78,5){\(\tau_{n+1}\)}
    \put(102,18){\(f^{n+1}(\alpha)^*\)}
    \put(95,7){\(\eta_n\)}
    \put(56,12){\(L\)}
    \put(73,-5){\(\ell\)}
    \end{overpic}
    \vspace{.5cm}
    \caption{Schematic pictures in Case 1}\label{fig:balloons}
    \end{center}
    \end{figure}

    \begin{proof}

      Recall that $\beta$ is obtained by joining $\gamma$ with $f(\gamma)$
      along some arcs in $\ell$. Let $\eta_n$ be the geodesic arc joining
      $f^n(\alpha)^*$ and $f^{n+1}(\alpha)^*$, perpendicular to both, and
      homotopic to $\overline{\tau}_n * l_n * \tau_{n+1}$, where $l_n
      \subset \ell$ is the arc joining $\tau_n$  and $\tau_{n+1}$, via a
      homotopy keeping the endpoints on the two loops. Then $f^n(\beta)^*$
      is homotopic to $f^n(\alpha)^* * \eta_n * f^{n+1}(\alpha)^*
      *\overline\eta_n$. Note that the arcs $\eta_n$ must leave every
      compact set, for otherwise on a subsequence $n_j$ where $\eta_{n_j}$
      intersect a fixed compact set, the loop $f^{n_j}(\beta)^*$ and the
      concatenation $f^{n_j}(\alpha)^* * \eta_{n_j} * f^{n_j+1}(\alpha)^*
      *\overline\eta_{n_{j+1}}$ would be uniformly close to each other and
      $f^{n_j}(\beta)^*$ would also intersect a (slightly larger) fixed
      compact set.

      By Claim 1 and the fact that the endpoint of $\tau_n$ on $\ell$
      escapes to $\infty$, we can choose $i$ sufficiently large so that
      $\tau_{n_i}$ intersects $L$, and for $n\geq n_i$, $f^n(\alpha)^*$,
      $f^n(\beta)^*$, $\eta_n$, and $l_n$ are all disjoint from a large
      neighborhood of $L$.

      We now argue by induction on $n$ with base case $n=n_i$. Assume
      that for some $n\geq n_i$, $\tau_n$ intersects $L$. Since $\eta_n$ is
      homotopic to $\overline{\tau}_n * l_n * \tau_{n+1}$ we have a
      geodesic quadrilateral with sides $\overline{\tau}_n$, $l_n$,
      $\tau_{n+1}$, $\overline\eta_n$ and this gives a homotopy
      between $\tau_n$ and
      $\tau_{n+1}$ that keeps the endpoints outside $L$. In
      particular, $\tau_{n+1}$ also intersects $L$. This finishes the
      induction. \qedhere 

    \end{proof}

    From Claim 2, for all sufficiently large $n$, the arcs $\tau_n$ are
    all homotopic keeping the endpoints outside $L$. Thus for large
    $n$ and much larger $m\gg n$ we have a quadrilateral with sides
    $\overline\tau_n,l_{nm},\tau_m,\sigma_{nm}$ where
    $l_{nm}\subset\ell$ is a segment and $\sigma_{nm}$ is a path
    disjoint from $L$. We can replace $\sigma_{nm}$ by a geodesic and
    it will still be disjoint from $L$, and we have a geodesic
    quadrilateral with right angles at $l_{nm}$. When $m\gg n$ the segment
    $l_{nm}$ is long, far from $L$, and $\sigma_{nm}$ comes very close
    to its midpoint. This gives a homotopy of $\tau_n$ to a subarc of
    $\sigma_{nm}$ that keeps one endpoint fixed and moves the other
    along a half of $l_{nm}$ and then by
    a very small amount, so certainly the endpoints remain outside $L$
    throughout the homotopy.

    This contradiction shows that $\gamma$ is forward wandering. A similar
    proof also shows that $\gamma$ is backward wandering. This finishes
    this case.

    \emph{Case 2.} Both endpoints of $\gamma$ are on a periodic line
    $\ell$. Then apply Case 1 to the power of $f$ that preserves the
    oriented line $\ell$.

  \emph{General case.} We have seen that the limit set of an arc $\gamma$,
  if non-empty, can only consist of boundary lines. If $\ell$ is a boundary
  line in $\L(\gamma)$, let $\beta$ be an essential arc with both endpoints
  on $\ell$. By the previous cases, $\beta$ is wandering, so by Lemma
  \ref{lem:duality}
  \[ i(\L(\beta), \gamma)=0  \Longleftrightarrow i(\beta, \L(\gamma)) = 0.\]
  However, $\beta$ intersects $\ell \in \L(\gamma)$, a contradiction. 

  Finally, since every pair of arcs can be connected by a sequence of
  curves, it follows that every arc converges to the $f$-attracting end
  under forward iteration and to the $f$-repelling end under backward
  iteration. This finishes the proof.
  \qedhere 
  
  \end{proof}
  
  \begin{cor}\label{cor:double-is-wandering}
    Let $X$ be a connected bordered surface spanned by its curves. Let $f$
    be a homeomorphism of $X$ such that every curve is $f$-wandering and
    every boundary line is $f$-periodic. Let \(DX\) be the double of $X$,
    and let \(Df\) be the double of $f$ to $DX$. Then:
    \begin{enumerate}
      \item every curve in $DX$ is $Df$-wandering;
      \item if $e_+$ and $e_-$ are the $f$-attracting and $f$-repelling
        ends of $X$, then each determines a unique end $\hat{e}_+$ and
        $\hat{e}_-$ of $DX$, and these are precisely the
        $Df$-attracting and $Df$-repelling ends of $DX$,
        respectively.
  \end{enumerate}     
  \end{cor} 

  \begin{proof}
    
    For (1), every curve of $DX$ that is not contained in $X$ meets each
    copy of $X$ in a collection of arcs. By the previous proposition, every
    arc in this collection is wandering towards the same end of $X$.
    Therefore every curve in $DX$ is wandering. 

    For (2), by Lemma \ref{lem:e_pm}, there are unique ends \(\hat{e}_\pm\)
    towards which iterates of curves converge. As iterates of curves in
      \(X\) converge to \(e_\pm\), the claim follows. \qedhere

  \end{proof}

  \begin{lemma} \label{lem:covering}

    Let $X$ be a connected surface, possibly with boundary, and let $\phi$
    be a mapping class on $X$. Suppose that for some $n >
    0$ there exists a representative $g \in \phi^n$ that is a translation
    on $X$. Then $\phi$ is also represented by a translation. In
    particular, there is a hyperbolic metric on $X$ with totally geodesic
    boundary for which some representative $f \in \phi$ is an isometric
    translation. Moreover, we can choose the hyperbolic metric so that the
    first return of $f$ to any boundary component of $X$ has translation
    length $1$.  

  \end{lemma}
   
   \begin{proof}

     Fix a representative $f \in \phi$. Note that since $f^n \sim g$, the
     maps $f$ and $g$ commute up to isotopy. Consider the covering map \[
       p: X \to Q=X/ \langle g^2 \rangle.\] Note that $Q$ may be either of
     finite or of infinite type. Our goal is to show $\phi$ induces a
     well-defined mapping class $\bar{\phi}$ on $Q$, using the curve graph
     of $Q$. By replacing $n$ by a positive multiple if necessary, we may
     assume that $Q$ has high enough complexity so that the automorphism
     group of its curve graph coincides with the extended mapping class
     group of $Q$ \cite{Iva97, Kor99, Luo00, BDR20}. In what follows, we
     write $h \sim h'$ to mean that $h$ and $h'$ are isotopic.
    
     We first treat the case in which $Q$ has no boundary. For convenience,
     fix a hyperbolic metric on the further quotient $Q/\langle g \rangle$
     and equip $X$ with the lifted metric, so that $g$ acts on $X$ by an
     isometry. Let $\alpha \subset Q$ be a simple closed geodesic. If
     $\alpha$ lifts, then $p^{-1}(\alpha)=\bigcup \alpha_m$ is a disjoint
     union of simple closed geodesics, where $\alpha_m = g^{2m}(\alpha_0)$
     for some chosen lift $\alpha_0$. Let $\beta_m$ be the geodesic
     representative of $f(\alpha_m)$. Since $fg \sim gf$, we have $\beta_m
     = g^{2m}(\beta_0)$. It follows that there is a unique simple closed
     geodesic $\beta \subset Q$ such that $\beta = p(\beta_1)$. If $\alpha$
     does not lift, then $p^{-1}(\alpha)$ is a finite collection of lines,
     say $\{l_1,\ldots,l_k\}$, which are cyclically permuted by the action
     of $g^2$. Let $t_i$ be the geodesic representative of $f(l_i)$. As for
     curves, the relation $fg \sim gf$ implies that $g^2$ cyclically
     permutes $\{t_1,\ldots,t_k\}$. Hence there is again a unique simple
     closed geodesic $\beta \subset Q$ such that $\beta = p(\bigcup t_i)$.
     This allows us to define $f_*(\alpha) =\beta$. It is straightforward
     to check that $f_*$ preserves disjointness, and hence determines an
     automorphism of the curve graph of $Q$. Hence $f_*$ is induced by an
     extended mapping class\footnote{I.e.\ the homotopy class of a
     homeomorphism, which could be orientation-reversing.} $\bar{\phi}$ of
     $Q$. By construction, $\bar{\phi}$ depends only on the mapping class
     of $f$. 

     Since $f^n \sim g$, the mapping class $\bar{\phi}$ is periodic, and in
     fact $\bar{\phi}^n=\bar g$, where $\bar g$ is the involution on $Q$
     induced by $g$. By Proposition \ref{prop:Nielsenrealization}, there is
     a hyperbolic metric on $Q$ for which $\bar{\phi}$ admits an isometric
     representative. Lift both the metric and the isometry to $X$, choosing
     the lift $\tilde f: X \to X$ such that $\tilde f^n \sim g$. In
     particular, $\tilde f$ is an isometric translation of $X$. 

     We now claim $\tilde f \sim f$. Let $F_0$ be a fundamental domain for
     $g$, and $F_i=g^i(F_0)$. Since $f$ and $\tilde f$ induce the same
     mapping class on $Q$, our choice of the lift implies that $f^{-1}
     \circ \tilde f$ is isotopic to the identity on each $F_{2i} \cup
     F_{2i+1}$. Moreover, since both $\tilde f$ and $f$ commute with $g$ up
     to isotopy, $f^{-1} \circ \tilde f$ is also isotopic to the identity
     on $F_{2i+1} \cup F_{2i+2}$. To complete the proof, we use the
     Alexander method; see \cite{hmv_Alexander} for a version for
     infinite-type surfaces. For each $j$, choose a pants decomposition of
     $F_j$ which includes all the peripheral curves. For every
     nonperipheral pants curve $\alpha$, choose a transversal curve lying
     in $F_j$. If $\alpha$ is peripheral and lies in $F_j \cap F_{j+1}$,
     then pick a transversal curve lying in $F_j \cup F_{j+1}$. Doing this
     for every $j$ gives a system of curves fixed by $f^{-1} \circ \tilde
     f$ up to isotopy. Therefore, by the Alexander method, $f^{-1} \circ
     \tilde f$ is isotopic to the identity. 

     Now suppose that \(Q\) has boundary. Double \(Q\) and \(X\) along
     their boundaries to obtain \(DQ\) and \(DX\), each equipped with an
     involution that exchanges the two copies. Also, double the maps 
     \(Df,Dg\colon DX\to DX\), and we still have \((Df)^n\sim
     Dg\), with the homotopy commuting with the involution of $DX$. Since
     \(DQ=DX/\langle Dg^2 \rangle\) has no boundary, the previous argument
     applies to \(Df\colon DX\to DX\). Thus \(Df\) induces a periodic
     mapping class \(\hat\phi\) on \(DQ\). Because \(Df\) is obtained by
     doubling \(f\), the mapping class \(\hat \phi\) commutes with the
     involution  \(\tau\colon DQ\to DQ\). By
     Proposition~\ref{prop:Nielsenrealization}, after
     choosing a suitable hyperbolic metric on \(DQ\), we may represent
     \(\hat\phi\) by an isometry that commutes with \(\tau\). It follows
     that this isometry preserves the fixed-point set of \(\tau\), namely
     the doubled boundary \(\partial Q\), and hence restricts to an
     isometry $\bar{f}$ of $Q$ representing $\bar{\phi}$. In particular,
     the boundaries of $Q$ are totally geodesic. Moreover, since any
     hyperbolic metric on $DQ$ preserved by $\tau$ is the doubling of a
     hyperbolic metric on $Q$, we may prescribe, for each $\bar{f}$-orbit
     of compact boundary components of $Q$, the common boundary length to
     be any positive number $L$. Lift both the metric and the isometry to
     \(X\), choosing the lift \(\tilde f\colon X\to X\) such that $\tilde
     f^{\,n}\sim g$. By a similar argument as above, \(\tilde f\sim f\). 
     
     Finally, let $l$ be a boundary component of $X$, and suppose that
     $\tilde f^k$ is the first return to $l$. Then $l$ projects to a
     compact boundary component of $Q$ of length $L$. By choosing $L=2n/d$,
     where $d=\text{gcd}(2n,k)$, we can arrange that the first return of
     $\tilde f$ to $l$ is a translation by length 1. \qedhere

  \end{proof}

  \begin{thm}\label{thm:translation}
    Let $X$ be a connected possibly bordered surface spanned by its curves.
    Let $f$ be a homeomorphism of $X$ such that every curve is
    $f$-wandering and every boundary line is $f$-periodic. Then there is
    a hyperbolic metric on $X$ with totally geodesic boundary and an
    isometric translation on $X$ isotopic to $f$.
  \end{thm}

  \begin{remark}
    If the surface has empty boundary, Theorem \ref{thm:translation}
    reduces to Theorem \ref{thm:translationintro}.
  \end{remark}

  \begin{proof}
    By Lemma \ref{lem:e_pm}, there are ends $e_+$ and $e_-$ of $X$ to which
    all curves converge under forward/backward iteration. In particular,
    $f(e_\pm)=e_\pm$. We now distinguish two cases.  

    \emph{Case 1: $e_+\neq e_-$.} Let $\delta$ be a curve or a minimal
    multi-arc separating $e_+$ and $e_-$ as in Lemma \ref{lem:delta}. Thus
    we have $X=X_-\cup_\delta X_+$ where $X_{\pm}$ is a neighborhood of
    $e_{\pm}$. Choose \(n\neq 0\) so that $f^n(X_+)\subset X_+$. Each
    $\delta_i=f^{ni}(\delta)^*$ is a multiarc or a curve separating $e_+$
    and $e_-$ with complementary components nested. This shows that the
    dual graph is a line on which $f^n$ acts as a translation. We may
    therefore isotope $f^n$ to a translation $g$ of $X$, with the region
    between $\delta_0$ and $\delta_1$ serving as a fundamental domain for
    the action of $g$. This case is finished by Lemma \ref{lem:covering}.
   
    \emph{Case 2: $e_+=e_-$.} We first treat the case in which $X$ has no
    boundary, and then reduce the bordered case to it by doubling. 

    The goal is to construct an exhaustion of $X$ by $f$-invariant
    subsurfaces, each of which satisfies the hypothesis of Lemma
    \ref{lem:covering}. Choose a curve $\alpha_1$ and then choose a finite
    sequence of curves
    \[\alpha_1,\alpha_2,\ldots, \alpha_k = f(\alpha_1)\] such that
    $i(\alpha_i,\alpha_{i+1}) > 0$. Since every curve in $X$ is wandering,
    the union of the $\langle f \rangle$-orbits of the curves $\alpha_i$
    forms a locally finite graph $\Gamma$. Therefore, they span a
    subsurface $T_1$ which is $f$-invariant up to isotopy. We claim $f$ is
    isotopic to a translation of $T_1$. Since the curves are all wandering,
    there exists $n_0\ge 1$ such that $\alpha_i$ is disjoint from
    $f^n(\alpha_j)$, for all $n \ge n_0$ and $i,j \in \{1,\ldots,k\}$. It
    follows that the induced action by $\langle f \rangle \cong \mathbb{Z}$
    on $\Gamma$ is proper and cocompact, so $\Gamma$ is two-ended. The
    regular neighborhood $N(\Gamma)$ is also two-ended, and since curves in
    $N(\Gamma)$ are $f$-wandering, these ends are the attracting and
    repelling endpoints. Applying Case 1 to $N(\Gamma)$ we obtain a
    translation $f_1$ of $N(\Gamma)$ isotopic to $f|_{N(\Gamma)}$. The
    subsurface $T_1$ is obtained from $N(\Gamma)$ by capping off some
    compact boundary components with disks or once-punctured disks, and
    possibly attaching half-planes along noncompact boundary components.
    Any disk or punctured-disk attachments occur in $\mathbb{Z}$-orbits,
    hence the translation extends over them. Furthermore, since $T_1$ is
    spanned by finitely many curves, it has only finitely many noncompact
    boundary components. So some power of $f_1$ fixes each of them and acts
    on each as a translation. It follows that some power of $f_1$ also
    extends across any attached half-plane as a translation. This shows
    some power of $f$ is isotopic to a translation of $T_1$, and hence the
    same is true of $f$ by Lemma \ref{lem:covering}, proving the claim.
    Repeating this procedure for a larger collection of curves yields a
    sequence of subsurfaces $T_1 \subset T_2 \subset \cdots$ exhausting
    $X$, each of which is $f$-invariant up to isotopy, and $f$ is isotopic
    to a translation on each $T_i$. See Figure \ref{fig:Tis}. 

    \begin{figure}[htp!]
    \begin{center}
    \begin{overpic}[width=.5\textwidth]{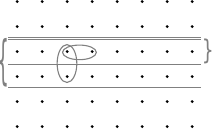}
    \put(102,33){\color{MyGray}$T_1$}
    \put(-7,29){\color{MyGray}$T_2$}
    \end{overpic}
    \caption{An example of the exhaustion $T_1\subset T_2\subset\dots$,
      where the map is the horizontal shift to the right, $T_1$ is the
      surface spanned by the orbit of the ``horizontal'' curve and $T_2$
      the surface spanned by the orbit of the ``horizontal'' and the
      ``vertical'' curve.}\label{fig:Tis}
    \end{center}
    \end{figure}
    
    We now construct, for each $i$, a hyperbolic metric on $T_i$ with
    totally geodesic boundary extending the hyperbolic metric on $T_{i-1}$
    together with a map $f_i \sim f$ preserving $T_i$ such that
    \begin{itemize}
      \item $f_i|_{T_i}$ is an isometric translation;
      \item the first return of $f$ to any noncompact boundary of $T_i$  is
        a translation of length $1$; and
      \item $f_i|_{T_{i-1}} = f_{i-1}|_{T_{i-1}}$.
    \end{itemize}
    The base case $i=1$ comes from Lemma \ref{lem:covering}. Now assume $i
    \ge 2$, and choose $g \sim f$ such that $g(T_i) = T_i$ and
    $g|_{T_{i-1}} = f_{i-1}|_{T_{i-1}}$. Let $Y$ be the closure of a
    component of $T_i \setminus T_{i-1}$. Since $T_i$ is spanned by
    finitely many $f$-orbits of curves, there is a smallest $k>0$ such that
    $g^k(Y)=Y$. Let $\mathcal{Y}=\bigsqcup_{i=1}^k g^i(Y)$, regarded as a
    possibly disconnected surface. The restriction $g|_{\mathcal{Y}}$ is
    isotopic to a translation on $\mathcal Y$, in the sense that the first
    return to each component is a translation. Therefore, there is a
    hyperbolic metric on $\mathcal Y$ together with an isometry
    $g_{\mathcal Y}  \sim g|_{\mathcal Y}$, such that the first return of
    $g_{\mathcal Y}$ to any noncompact boundary of any $g^i(Y)$ is a
    translation by length $1$. Gluing the hyperbolic metric on $T_{i-1}$ to
    the one on $\mathcal{Y}$ via the identifications coming from $X$ yields
    a hyperbolic metric on $T_{i-1} \cup \mathcal{Y}$. Moreover, along any
    common boundary between $T_{i-1}$ and $\mathcal Y$, the maps $f_{i-1}$
    and $g_{\mathcal Y}$ have first returns that are translations of length
    $1$ in the same direction. Hence these maps glue to an isometric
    translation on $T_{i-1} \cup \mathcal{Y}$ isotopic to $g|_{T_{i-1} \cup
    \mathcal{Y}}$. Performing this construction for each $g$-orbit of
    components of $T_i \setminus T_{i-1}$ yields a hyperbolic metric on
    $T_i$ extending the one on $T_{i-1}$, together with an isometry on
    $T_i$ extending the translation on $T_{i-1}$. We can now isotope $f$ to
    a map $f_i$ that agrees with this translation on $T_i$. 

    Passing to the limit, we obtain a hyperbolic metric on $X$ with respect
    to which $f$ is isotopic to an isometric translation.
   
    Finally, in the case that $X$ is bordered, let $DX$ be the double of
    $X$, and let $Df$ be the doubled homeomorphism. By Corollary
    \ref{cor:double-is-wandering}, every curve in $DX$ is $Df$-wandering,
    and the common attracting and repelling end $e$ of $X$ determines an
    end $\hat e$ of $DX$ such that every curve in $DX$ converges to $\hat
    e$ under both forward and backward iteration of $Df$. Thus $DX$ with
    the map $Df$ satisfies the hypotheses of the boundaryless case above.
    This gives a hyperbolic metric on $DX$ and an isometric translation $g$
    isotopic to $Df$. All of these choices can be made equivariantly with
    respect to the involution exchanging the two copies of $X$. Hence the
    resulting metric is invariant under this involution, and its
    restriction to $X$ is a hyperbolic metric with totally geodesic
    boundary. With respect to this restricted metric, $g|_X$ acts by a
    translation on $X$ and is isotopic to $f$. \qedhere

  \end{proof}
    
  \begin{remark}
    In the previous proof, in the case $e_+=e_-$, passing to the double is
    necessary if $X$ has nonempty boundary. This is because the subsurfaces
    $T_i$ may only exhaust the interior of $X$, in which case we do not
    obtain a hyperbolic structure on all of $X$ (some boundary components
    would lie at infinity with respect to the metric on the interior of
    $X$).
  \end{remark}
    
  We end this section with an example showing that the periodic boundary
  line hypothesis is necessary in Proposition \ref{prop:wandering-arcs},
  Corollary \ref{cor:double-is-wandering} and Theorem
  \ref{thm:translation}. 
    
  \begin{example}\label{ex:not-translation-double}

    Let \(R\) be the Roller compactification (see \cite{roller_poc}) of
    \(\R^2\) with its standard cube complex structure. Let \(X\) be the
    surface obtained by removing from \(R\) the closure of all points in
    \(\R^2\) with integer coordinates. The homeomorphism \(f\) of \(X\) is
    the continuous extension of the map \((x,y)\mapsto (x+1,y+1)\) of
    \(\R^2\). Then every curve of \(X\) and every boundary line is
    wandering. On the other hand, there are arcs from boundary to boundary
    which are not wandering, such as the closure of the line
    \(y=\frac{1}{2}\subset \R^2\ssm\Z^2\). In particular, the map induced
    by $f$ on $DX$ is not a translation.

  \end{example}

  \begin{figure}[htp!]
  \includegraphics[scale=2]{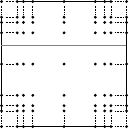}
  \caption{The surface $X$ from Example
    \ref{ex:not-translation-double}, with a non-wandering horizontal arc}
  \end{figure}

\section{Structure theorem of well-tempered maps}

  \label{sec:structure}
  
  The goal of this section is to prove the following classification result for
  well-tempered maps:
  
  \begin{thm}\label{thm:decomposition} Let $f$ be a well-tempered homeomorphism of a
    surface $S$ of infinite type. Then, up to modifying $f$ by an isotopy,
    there is a hyperbolic metric on \(S\) and a canonical decomposition of
    $S$ into three $f$-invariant subsurfaces $\Sper$, $\Sinf$ and $S_0$
    with totally geodesic boundary such that, for every component $X$ of
    $\Sper$, $\Sinf$ and $S_0$, the map $f$ returns to $X$ and the first
    return is isotopic to: \begin{itemize} \item a periodic isometry if
    $X\subseteq \Sper$ or $X\subseteq S_0$; \item an isometric translation
    if $X\subseteq \Sinf$. \end{itemize} Furthermore, each component $X$ of
    $S_0$ is topologically trivial, i.e.\ $X$ contains no essential
    nonperipheral curves. \end{thm}
    
  We begin by developing several properties of the canonical decomposition
  associated to a well-tempered map $f$. In particular, we describe how
  components of $\Sper$, $\Sinf$, and $S_0$ can neighbor each other. We
  also show that, for each component $X$ of the decomposition, the map $f$
  returns to $X$. We then combine these results with those of Section
  \ref{sec:translation} and the work of Afton-Calegari-Chen-Lyman
  \cite{accl_Nielsen} to complete the proof of Theorem
  \ref{thm:decomposition}. Finally, we show how to deduce Theorem
  \ref{thm:mainintro} from Theorem \ref{thm:decomposition}.

  \subsection{Properties of the canonical decomposition}

  For a well-tempered map, recall the almost geodesic representatives of
  $\Sper$, $\Sinf$, and $S_0$ defined in Section \ref{sec:decomposition}.
  An essential component of $\Sper$ or $\Sinf$ either has infinite type or
  has negative Euler characteristic; in particular, such a component always
  contains an essential curve. 
  
  Two components of the decomposition are \emph{adjacent} if they have
  boundary components that are properly isotopic (equivalently, they have
  representatives sharing at least one boundary component).

  A component $X$ is \emph{self-adjacent} along $\alpha$ if $\alpha$ is
  properly homotopic to two distinct boundary components of $X$.
  Equivalently, $X$ is self-adjacent if there is an inessential component
  $Y$ which shares two boundary components with $X$.

  In this section we will show various properties of the canonical
  decomposition. The first result excludes strip components and
  self-adjacency (and will follow from Lemma \ref{lem:no-self-adjacent} and
  Corollary \ref{cor:no-strips}).

  \begin{prop}\label{prop:properties}
    No component of the canonical decomposition is a strip, nor is any
    component self-adjacent.
  \end{prop}

  The next proposition describes the components of \(\Sinf\), and it is
  proved in Lemmas \ref{lem:no-compact} and Lemma
  \ref{lem:Sinf-not-wandering}:

  \begin{prop}\label{prop:propertiesSinf}
  Let \(X\) be a connected component of \(\Sinf\). Then:
  \begin{enumerate}
    \item no boundary component of \(X\) is compact;
    \item every boundary component of \(X\) is periodic;
    \item \(f\) returns to \(X\) and there is a hyperbolic structure on
      \(X\) such that the first return map is isotopic to an isometric
      translation.
  \end{enumerate}
  \end{prop}

  The following properties of \(\Sper\) are shown in Lemmas
  \ref{lem:per-on-per} and \ref{lem:Sper-neighbor}:

  \begin{prop}\label{prop:propertiesSper}
    No two components of \(\Sper\) are adjacent. If \(X\) is a component of
    \(\Sper\), \(f\) returns to \(X\) and there is a hyperbolic structure
    on \(X\) such that the first return map is isotopic to a finite-order
    isometry.
  \end{prop}

  Finally, we describe the properties of \(S_0\):

  \begin{prop}\label{prop:propertiesS0}
   The following are the topological possibilities for
      a component $X$ of $S_0$:

      \begin{enumerate}
        \item A closed disk with at least three points removed from the
          boundary. In this case, at most one neighbor of $X$ is in $\Sper$.
        \item A once-punctured closed disk with at least one point removed
          from the boundary. In this case, all neighbors of $X$ are in $\Sinf$.
        \item A cylinder with one compact boundary component and at least one point
          removed from the other boundary (a crown). In this case, $X$ has
          exactly one neighbor in $\Sper$ along its compact boundary component. 
          \item A pair of pants. In this case, all neighbors of \(X\) are
            annular components of \(\Sper\).
      \end{enumerate}
  \end{prop}

\begin{figure}[htp!]
      \begin{center}
        \begin{overpic}{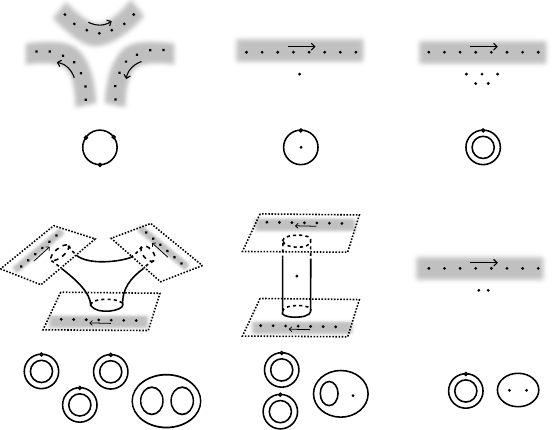}
          \put(-2,67){$f_1:$}
          \put(5,49){$S_0\simeq$}
          \put(35,67){$f_2:$}
          \put(40,49){$S_0\simeq$}
          \put(69,67){$f_3:$}
          \put(74,49){$S_0\simeq$}
          \put(-7,28){$f_4:$} 
          \put(-6,5){$S_0\simeq$}
          \put(37,28){$f_5:$}
          \put(38,5){$S_0\simeq$}
          \put(67,28){$f_6:$}
          \put(70,5){$S_0\simeq$}
        \end{overpic}
        \caption{Examples of well-tempered mapping
        classes}\label{fig:componentsofS0}
      \end{center}
    \end{figure}

  Note that all topological possibilities from Proposition
  \ref{prop:propertiesS0} occur for components of $S_0$, as shown in Figure
  \ref{fig:componentsofS0}.   In each case, the map is a puncture shift in
  each shaded strip, in the direction indicated by the arrow, and it is the
  identity outside the strips.
    
    The rest of the section is dedicated to showing the previous
    propositions, whose proofs are intertwined. We start by showing that
    certain boundary components have to be periodic.

    \begin{lemma} \label{lem:jointly-wandering1}
      
      If two components of $S_\infty$ share a common boundary component
      $\alpha$, then $\alpha$ is periodic. Similarly, if a component of
      $\Sinf$ is self-adjacent along $\alpha$, then $\alpha$ is periodic.
    
    \end{lemma}

    \begin{proof}
      Let $X$ and $Y$ be (not necessarily distinct) components of
      $S_\infty$ with a common boundary component $\alpha$. Since $X$ and
      $Y$ are essential, we can find a curve $\beta\subset X \cup Y$ that
      intersects $\alpha$ essentially. Because $\beta$ is not contained in
      $\Sinf$, it is not wandering, so $\L(\beta)\neq \emptyset$. Choose
      \(\ell\in \L(\beta)\) and let \(\gamma\) be a curve intersecting
      \(\ell\). Since \(i(\gamma,\L(\beta))\neq 0\), Lemma
      \ref{lem:duality} implies \(i(\L(\gamma),\beta)\neq 0\), hence
      \(\gamma\) is not wandering. By Lemma \ref{lem:SinftySubs},
      \(\L(\gamma)\) cannot intersect any component of $S_\infty$ (in
      particular, \(X\) and \(Y\)) essentially. On the other hand,
      $\L(\gamma)$ intersects $\beta$. This is only possible if $\alpha \in
      \L(\gamma)$. By $f$-invariance, the entire orbit of $\alpha$ lies in
      $\L(\gamma)$. As \(f\) is well tempered, the orbit of \(\alpha\) is
      finite, i.e.\ \(\alpha\) is periodic. \qedhere  

    \end{proof}

    \begin{lemma} \label{lem:jointly-wandering2}
       
      Let $X$ be a component of $S_0$. If $X$ either
      \begin{itemize}
        \item shares at least one boundary component with $S_\infty$ but
          contains a non-contractible curve not homotopic to
          that boundary component;
          or \item shares at least two boundary components with $\Sinf$,
      \end{itemize}
      then $X$ is not wandering.

    \end{lemma}

    \begin{proof}
       
      In the first case, let $\gamma$ be a non-contractible curve in $X$
      and let $Y$ be a component of $\Sinf$ sharing a boundary component
      $\alpha$ with $X$, where $\gamma$ is not homotopic to $\alpha$. Since
      $Y$ is essential, we can find a curve $\beta \subset X \cup Y$
      crossing $\alpha$, for instance by joining a curve in $Y$ to $\gamma$
      along an arc and taking the boundary of a regular neighborhood of
      their union. In the second case, let $Y$ and $Y'$ be components of
      $\Sinf$ (possibly $Y=Y'$) sharing boundary components $\alpha$ and
      $\alpha'$, respectively, with $X$. We can then find a curve $\beta
      \subset X \cup Y \cup Y'$ intersecting $\alpha$ and $\alpha'$, by
      joining a curve in $Y$ to a curve in $Y'$ by an arc that crosses
      $\alpha$ and $\alpha'$. 

      We claim that if $X$ is wandering, then $\beta$ is wandering, which is
      a contradiction. The argument is similar to above. Indeed, suppose
      $\beta$ is not wandering. Then $\L(\beta)\neq \emptyset$, but by Lemma
      \ref{lem:SinftySubs}, $\L(\beta)$ cannot intersect $f^n(Y)$ or
      $f^n(Y')$ essentially. On the other hand, $f^n(\beta)$ is contained in
      $f^n(X \cup Y)$ or $f^n(X \cup Y \cup Y')$, so some line $\ell \in
      \L(\beta)$ must be homotopic into $f^n(X)$ for some $n$. Since
      $\L(\beta)$ is $f$-invariant and $X$ is wandering, infinitely many
      iterates of $\ell$ would lie in $\L(\beta)$, contradicting the
      finiteness of $\L(\beta)$. \qedhere

    \end{proof}

    We now show that all boundary components of $\Sinf$ are lines.
    \begin{lemma} \label{lem:no-compact}
      
      No component $X$ of $\Sinf$ has a compact boundary component. 
      
    \end{lemma}

    \begin{proof}

      We argue by contradiction. Suppose that $\alpha$ is a compact boundary
      component of $X$. Since $X$ is spanned by curves, $\alpha$ is a
      wandering curve. Therefore the component $Y$ adjacent to $X$ along
      $\alpha$ cannot belong to $\Sper$ (since boundary curves of $\Sper$ are
      periodic) or to $S_\infty$ (by Lemma \ref{lem:jointly-wandering1}).
      Hence $Y$ is a component of $S_0$. 

      If there is no $n>0$ such that $f^n(Y)$ is isotopic to $Y$, then $Y$ is
      wandering (as $S_0$ is an $f$-invariant subsurface). By Lemma
      \ref{lem:jointly-wandering2}, $Y$ cannot contain a non-contractible
      curve different from $\alpha$. If $\alpha$ were the only boundary
      component of $Y$, then $Y$ is a closed disc with at most one puncture,
      which is impossible. Thus $Y$ has at least two boundary components,
      both wandering. Hence these two boundary components belong to $\Sinf$,
      contradicting Lemma \ref{lem:jointly-wandering1}.
      
      Therefore there exists $n>0$ such that $f^n(Y)$ is isotopic to $Y$.
      Since $\alpha$ is wandering, $Y$ contains infinitely many iterates of
      $\alpha$, and hence contains an essential nonperipheral curve. This
      contradicts Proposition \ref{prop:no-ess-nonper}, completing the proof.
      \qedhere
      
    \end{proof}
    
    Using the previous results, we can deduce that no component of $S_0$ is
    wandering; moreover, for each such component the first return map has
    finite order.

    \begin{lemma}\label{lem:per-on-S0} 
      
      For every component $X$ of $S_0$, $f$ returns to $X$ and there is a
      hyperbolic structure in which every compact boundary component has
      length one, such that the first return $f^k$ is isotopic to a periodic
      isometry of $X$. 
    
    \end{lemma}

\begin{figure}[htp!]
  \begin{center}
  \begin{overpic}{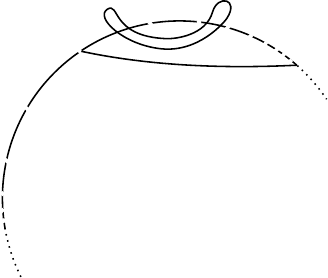}
    \put(70,67){\(Z_0\)}
    \put(26,73){\(\ell\)}
    \put(50,58){\(\ell_0\)}
    \put(-6,14){\(e_-\)}
    \put(95,64){\(e_+\)}
    \put(66,86){\(\alpha\)}
    \put(30,30){\(X\)}
    \end{overpic}
    \caption{The situation in the proof of Lemma
      \ref{lem:per-on-S0}, for \(e_-\neq \e_+\).}
    \label{fig:wanderingboundary}
  \end{center}
  \end{figure}

    \begin{proof}
    
      Recall that $X$ has no essential, nonperipheral curves, and by our
      construction of $\Sper$ and $\Sinf$, $X$ cannot be a closed disk or a
      closed disk with one puncture. If $f$ does not return to $X$, then
      since $S_0$ is a subsurface and $f$-invariant, $X$ is wandering. In
      particular, $X$ cannot border $\Sper$, so it cannot have a compact
      boundary component by Lemma \ref{lem:no-compact}. In this case, $X$
      is a disk with at most one puncture and points removed from the
      boundary, which yields a contradiction by Lemma
      \ref{lem:jointly-wandering2}. Hence $f$ returns to $X$. 

      Up to replacing $f$ with a power, we may assume $f(X)=X$. If $X$ has
      finitely many boundary components, then $f$ has finite order.
      Otherwise, $X$ is a disk with infinitely many points removed from the
      boundary and possibly a puncture or an open disk removed from its
      interior. Consider the action of $f$ on the noncompact boundary
      components of $X$. If the $f$-orbit of every component is infinite,
      then $f$ is semi-conjugate to an irrational rotation of the disk,
      which violates the fact that \(f\) is well tempered. Thus, there
      exists a finite orbit; in that case, all periodic components of
      $\partial X$ have the same period, and the same holds for the
      periodic ends. Up to passing to a power, we can assume all periodic
      boundary components and all periodic ends are fixed. 

      Suppose by contradiction that there is a wandering boundary component
      \(\ell\). By the structure of orientation-preserving homeomorphisms
      of the circle, there are ends \(e_\pm\) (possibly \(e_+=e_-\)) of
      \(X\) such that \(f^n(\ell)\to e_\pm\) as \(n\to\pm\infty\). Let
      \(Y\) be the component of the canonical decomposition adjacent to
      \(X\) along \(\ell\). Since one of its boundary components is
      wandering, \(Y\subset \Sinf\). Let \(\alpha\) be a curve in
      \(\interior(X)\cup\interior(Y)\cup f(\interior(Y))\cup \ell\cup
      f(\ell)\) which intersects \(\ell\) essentially (see Figure
      \ref{fig:wanderingboundary}). Note that \(\alpha\) is not wandering,
      so \(\L(\alpha)\neq \emptyset\). Without loss of generality, suppose
      \(\L^+(\alpha)\neq \emptyset\). Let $\ell_0$ be the line in $X$
      spanned by $e_+$ and the end of $\ell$ such that $f^i(\ell)$, $i \ge
      0$, lies on one side of $Y \ssm \ell_0$, which also contains
      \(\alpha\). Let $Z_0$ be the subsurface corresponding to that side
      (see Figure \ref{fig:wanderingboundary}), and let $Z_i = f^i(Z_0)$.
      Then $Z_i$ converges to $e_+$ as $i\to\infty$. But
      \(\L^+(\alpha)\subset Z_i\cup\Sinf\) for every \(i\geq 0\), and
      \(\L^+(\alpha)\) cannot intersect \(\Sinf\) essentially, so we get a
      contradiction. Let $k$ be the first return time to $X$. Then we
      obtain a hyperbolic metric on $X$ with respect to which $f^k$ is
      isotopic to a finite order isometry by looking at the quotient
      $X/\langle f^{k} \rangle$, choosing a hyperbolic structure with
      corners on the quotient and lifting it to $X$.  \qedhere

  \end{proof}

  We can now establish the properties of the first return map for
  components of $\Sper$.

  \begin{lemma}\label{lem:per-on-per}
          
    Let $X$ be a component of $\Sper$. Then $f$ returns to $X$ and there is
    a hyperbolic metric on $X$ with compact boundary components of length
    one such that the first return $f^k$ is isotopic to a periodic isometry
    of $X$. 
          
  \end{lemma}

  \begin{proof}
          
    Since $X$ is connected and every curve in $X$ is $f$-periodic, $f$
    must return to $X$. Up to replacing $f$ by a power we may assume $f(X)$
    is isotopic to $X$. If $X$ has bounded topology (i.e.\ $X$ is finite
    type or the double of $X$ is of finite type), then $f$ is isotopic to a
    periodic map of $X$. Otherwise, choose a finite collection $\C$ of
    orbits of curves in $X$ and let $F=F(\C)$ be the subsurface spanned by
    $\C$. After possibly enlarging $\C$, we can assume $F$ is connected and
    has negative Euler characteristic. Since $\C$ is $f$-invariant so is
    $F$; moreover, since all curves in $\C$ are $f$-periodic, $f|_F$ is
    isotopic to a periodic map of $F$ by Lemma
    \ref{lem:periodic-finite-type}. Let $n > 0$ be the period of $f|_F$.
    For any curve $\alpha$ not in $F$, we can find a larger connected
    $f$-invariant finite-type subsurface $F'$ containing both $F$ and
    $\alpha$. By the same reasoning, $f$ is also isotopic to a periodic map
    on $F'$, but since $F' \supset F$, $f|_{F'}$ also has period $n$ (by
    classical Nielsen-Thurston theory). In particular, $f^n(\alpha)$ is
    isotopic to $\alpha$ for all $\alpha$.  By the Alexander method
    (\cite{hmv_Alexander}), this implies that $f|_X$ is isotopic to a
    periodic map of $X$. The fact that we can realize the first return
    $f^k$ to $X$ by a periodic isometry is a direct consequence of
    Proposition \ref{prop:Nielsenrealization} (see \cite{accl_Nielsen}).
    \qedhere
          
  \end{proof}

  As a consequence, we can describe more precisely the properties of
  pair-of-pants components of \(S_0\).

  \begin{cor}\label{cor:changing-the-decomposition}
    Each pair of pants component of $S_0$ is adjacent only to annular
    components of \(\Sper\). Furthermore, no component of $S_0$ is
    self-adjacent.
  \end{cor}

  \begin{proof}
  
    Since a pair of pants component $X$ of $S_0$ has only compact boundary
    components, it can only border components of \(\Sper\). Let \(Y\) be
    one such component; if $Y$ is not an annulus, we can construct a curve
    $\alpha\subset X\cup Y$ intersecting the boundary of $X$ essentially.
    As some power of $f$ is periodic on both $X$ and $Y$, and as by the
    tempered condition there is no twisting allowed along the boundary,
    $\alpha$ is periodic, a contradiction.
          
    If a component $X$ of $S_0$ were self-adjacent along a boundary
    component $\alpha$, then $\alpha$ must be a curve--otherwise there is a
    component $Y$ which is a strip with both boundary components in $X$.
    Such a component $Y$ must belong to $S_0$, since it has no curves, but
    then $X\cup Y$ would be part of the same component of $S_0$. If
    $\alpha$ is a curve, then by Lemma \ref{lem:per-on-S0} $\alpha$ is
    periodic, so the annulus $Z$ about $\alpha$ is a component of $\Sper$.
    Since $X$ has at least two compact boundary components, it is a pair of
    pants, but then as before we could construct a curve $\beta$ contained
    in the genus-one surface $X\cup Z$ which intersects $\alpha$
    essentially and is periodic, a contradiction. \qedhere

  \end{proof}

  Our next step is to exclude (self-)adjacency of components of $\Sper$.

  \begin{lemma} \label{lem:Sper-neighbor}
          
    No two components of $\Sper$ can be adjacent, nor can a component
    of $\Sper$ be self-adjacent. 
          
  \end{lemma}

  \begin{proof}
          
    First suppose $X$ and $Y$ are essential components of $\Sper$ and
    $\alpha$ is a common boundary component. Choose a regular neighborhood
    $N(\alpha)$ about $\alpha$ contained in $X \cup Y$. There exists $n$
    such that $f^n$ is isotopic to the identity on $X$ and $Y$. Isotope
    $f^n$ so that it takes $N(\alpha)$ to $N(\alpha)$ fixing $\partial
    N(\alpha)$ pointwise, and further isotope $f^n$ so that it is the
    identity on truncated subsurfaces $X \setminus N(\alpha)$ and $Y
    \setminus N(\alpha)$. If $N(\alpha)$ is a strip, then we can further
    isotope $f^n$ (rel boundary) inside of $N(\alpha)$ to the identity map.
    If $\alpha$ is a curve, then since $f$ is (well) tempered, it cannot twist
    about $\alpha$, so $f^n$ is also isotopic (rel boundary) to the
    identity map on $N(\alpha)$. But this implies $f^n$ is isotopic to the
    identity on $X \cup Y$, and in particular, we can find an essential
    curve $\beta$ in $X \cup Y$ crossing $\alpha$ which is $f$-periodic. 
    
    The proof that a component $X$ has no self-adjacency along a periodic
    component $N(\alpha)$ is similar. In this case, we can isotope $f^n$ to
    be identity on $X \cup N(\alpha)$, which contradicts that $X$ and
    $N(\alpha)$ are components of the decomposition. \qedhere 
          
  \end{proof}

  With all these results at hand, we can prove the wanted properties of the
  first return map of components of $\Sinf$.
    
  \begin{lemma}\label{lem:Sinf-not-wandering}

    For every component $X$ of $\Sinf$, each boundary component is
    periodic. In particular, $f$ returns to $X$, and there is a hyperbolic
    metric on $X$ such that the first return $f^n$ is isotopic to an
    isometric translation on $X$.

  \end{lemma}
  
  \begin{proof} 
    
    Let $\alpha$ be a boundary component of $X$. The neighbor $Y$ of $X$
    along \(\alpha\) cannot be in \(\Sinf\) (by Lemma
    \ref{lem:jointly-wandering1}). So it is in $\Sper$ or $S_0$. In both
    cases, the first return map to \(Y\) is periodic, and so \(\alpha\)
    must be periodic. In
    particular, $f$ must return to $X$. By Theorem \ref{thm:translation},
    the first return $f^n$ up to isotopy preserves a hyperbolic metric on
    $X$ on which it acts as an isometric translation. \qedhere

  \end{proof}

  Our next goal is to show that there is no self-adjacent component of the
  decomposition. We need a preliminary result.
    
  \begin{lemma} \label{lem:S0-neighbor}
    
    A component $X$ of $S_0$ which is not a pair of pants can share at most
    one boundary component with $\Sper$. In this case, if $\alpha$ is the
    common boundary between $X$ and $\Sper$, then every non-contractible
    curve in $X$ is homotopic to $\alpha$. 

  \end{lemma}

  \begin{proof}

    We now know $X$ is $f$-periodic (by Lemma \ref{lem:per-on-S0}), so the
    proof is similar to Lemma \ref{lem:Sper-neighbor}. The assumptions on
    $X$ are to rule out the existence of a periodic curve $\beta$ crossing
    the common boundary of $X$ and $\Sper$ essentially. \qedhere  

  \end{proof}

  We can now show the no self-adjacency result.

  \begin{lemma}\label{lem:no-self-adjacent}

    No component of the canonical decomposition is self-adjacent.

  \end{lemma}

  \begin{proof}
    Suppose $X$ is a self-adjacent component. By Corollary
    \ref{cor:changing-the-decomposition} and Lemma \ref{lem:Sper-neighbor},
    $X$ is a component of $\Sinf$. Since $\Sinf$ has no compact boundary
    components by Lemma \ref{lem:no-compact}, there is a strip component
    $Y$ of $S_0$, with core line $\ell$, so that $X$ is the only neighbor
    of $Y$. Let $\ell_1$ and $\ell_2$ be the boundary components of $Y$. By
    Lemma \ref{lem:Sinf-not-wandering}, there is some $n\geq 1$ so that
    $f^n$ fixes $\ell_1$ (and therefore $\ell_2$, since they are
    homotopic). Since $Y$ is a strip, an orientation on $\ell$ induces
    orientations on $\ell_1$ and $\ell_2$. Now pick a curve $\alpha$
    contained in $X\cup Y$ and intersecting both $\ell_1$ and $\ell_2$
    once. Then $\alpha\cap X$ is an arc from $\partial X$ to $\partial X$,
    and by Proposition \ref{prop:wandering-arcs} it is wandering. This
    implies that $f^n$ acts as a positive translation on $\ell_1$ (with
    respect to the chosen orientation) if and only if it acts as a positive
    translation on $\ell_2$. But then we can homotope $f^n$ on the strip,
    leaving it invariant on its boundary, so that it is a translation, and
    therefore $\alpha$ is wandering, a contradiction. \qedhere 
  
  \end{proof}

  An easy consequence of the previous lemma is the fact that there are no
  components of the decomposition which are strips:

  \begin{cor}\label{cor:no-strips}
    No component of $S_0$, $\Sper$ or $\Sinf$ is a strip.
  \end{cor}

  \begin{proof}
    Since by construction every component of $\Sper$ and $\Sinf$ contains
    an essential (possibly peripheral) curve, if there is a strip, it must
    be a component of $S_0$. But by construction, a strip can arise only if
    a component of $\Sper$ or $\Sinf$ is self-adjacent, which is impossible
    by Lemma \ref{lem:Sper-neighbor} and \ref{lem:no-self-adjacent}.\qedhere
  \end{proof}

  We can now describe the possible topological types for the components of
  $S_0$:
  
  \begin{proof}[Proof of Proposition \ref{prop:propertiesS0}]
    
    If \(X\) has negative Euler characteristic, we are in case (4) by
    Corollaries \ref{cor:pants}  and \ref{cor:changing-the-decomposition}. 
  
    Suppose then that the Euler characteristic is non-negative. Note that
    $X$ is not a sphere, being a subsurface of $S$. It also cannot be a
    disk or a once-punctured disk, by the construction of $\Sinf$ and
    $\Sper$. 
    
    If $\chi(X)=1$, then $X$ is a disk with a closed totally disconnected
    set removed from its boundary. Let $n \in \{0,1,\ldots,\infty\}$ be the
    cardinality of this set. If $n=0$ or 1, then $X$ is a disk or a
    half-plane, which is impossible. By Corollary \ref{cor:no-strips}, $n
    \ne 2$, so $n\geq 3$. By Lemma \ref{lem:S0-neighbor}, $X$ has at most
    one neighbor in $\Sper$.

    If $\chi(X)=0$, then the interior of $X$ is an open annulus. We already
    know that $X$ cannot be the once-punctured disk, and it cannot be a
    closed annulus. Thus, $X$ is either a once-punctured disk or a closed
    annulus, with at least one point removed from its boundary. If $X$ is a
    once-punctured disk with points removed from its boundary, then $X$ has
    a non-contractible curve which is nonperipheral, so it cannot border
    $\Sper$. If it is an annulus with points removed from both boundary
    components, then $X$ contains an essential nonperipheral curve,
    contradicting Proposition \ref{prop:no-ess-nonper}. Finally, if $X$ is
    an annulus with points removed from only one boundary component, then
    the core curve of $X$ is not homotopic to the non-compact boundaries of
    $X$. Hence $X$ has at most one neighbor in $\Sper$, and therefore by
    Lemma \ref{lem:no-compact} it has exactly one neighbor in \(\Sper\),
    joined along the compact boundary component. \qedhere 

  \end{proof}

  Since a well-tempered map returns to each component of $S_0$, $\Sper$ and
  $\Sinf$, one might wonder if there is a uniform return time for all
  components. The answer is negative, as shown by the following lemma.
  Moreover, given a component of $\Sinf$, there is not always a uniform
  power of $f$ fixing all of its boundary components at once.
  
  \begin{lemma}\label{lem:no-uniform-return}
    There is a well-tempered map such that there is no uniform return time
    for all components of the canonical decomposition. Furthermore, the map
    can be chosen so that there is a component of $\Sinf$ with no uniform
    return time for its boundary components.
  \end{lemma}
  
  \begin{proof}

    Consider the bordered surface $$P= \left(\R_{\geq 0}\times
    \R\right)\ssm\bigcup_{n\geq 1}B_n,$$ where $B_n$ is the open disk of
    center $(n,0)$ and radius $\frac{1}{4}$. Let $a_n$ be the oriented
    segment of the $x$-axis from $\partial B_n$ to $\partial B_{n+1}$ and
    $a_0$ the oriented segment of the $x$-axis from the origin to $B_1$.
    Define the homomorphism $\psi:\pi_1(P)\to \Z$ by
    $$\psi(\alpha)=\sum_{n\geq 0}2^n\hat{\iota}(\alpha,a_n),$$ where
    $\hat{\iota}(\cdot,\cdot)$ denotes the algebraic intersection form.
  	
    \begin{figure}[htp!]
      \begin{center}
        \begin{overpic}[scale=1.4]{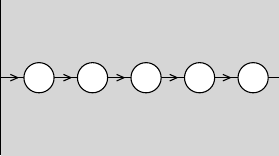}
          \put(2,32){$a_0$}
          \put(21,32){$a_1$}
          \put(40,32){$a_2$}
          \put(59,32){$a_3$}
          \put(78,32){$a_4$}
          \put(11,26){$B_1$}
          \put(30,26){$B_2$}
          \put(49,26){$B_3$}
          \put(68,26){$B_4$}
          \put(87,26){$B_1$}
        \end{overpic}
        \caption{The surface $P$ with the arcs $a_n$ in the proof of Lemma
        \ref{lem:no-uniform-return}.}
      \end{center}
    \end{figure}
    
    Let $p:\Sigma\to P$ be the associated cover of the kernel of $\psi$,
    and let $f$ be a generator of the deck group. Then $\Sigma$ is a
    bordered surface, whose boundary is a union of lines, given by lifts of
    the $y$-axis and of the curves $\partial B_n$. Note that, for any $n$,
    the number of lifts of $\partial B_n$ is $2^n-2^{n-1}$ and they are
    cyclically permuted by $f$. So there is no uniform power of $f$ fixing
    all boundary components of $\Sigma$.
    
    To turn this into an example of a well-tempered map on a borderless
    surface, we simply glue to each boundary component of $\Sigma$ a copy
    of the thrice-punctured half-plane $H$. We get a surface $S$ and a
    homeomorphism $\hat{f}$ of $S$ so that $\hat{f}|_\Sigma=f$ and for
    every copy of $H$, the first return map is the identity.
    
    The map $\hat{f}$ is well tempered and:
    \begin{itemize}
      \item $\Sinf=\Sigma$ and $\hat{f}|_{\Sinf}=f$, so there is no
        bound on the first return time of the boundary components of
        $\Sinf$;

      \item $\Sper$ is the union (over all attached copies of $H$) of
        components homeomorphic to a disk with three punctures, spanned by
        the curves in each copy of $H$. The first return time of a
        component equals the first return time of the boundary component of
        the corresponding copy of $H$, so these return times are also
        unbounded;

      \item $S_0$ is the union (over all attached copies of $H$) of
        components homeomorphic to an annulus with a point removed from the
        boundary. Again there is no uniform bound on the first
        return times of components of $S_0$. \qedhere
    \end{itemize} 
  	
  \end{proof}

  \begin{figure}[htp!]
  \begin{center}
  \includegraphics{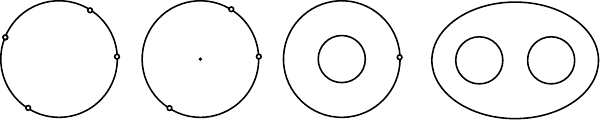}
  \caption{The topological possibilities for the components of $S_0$.}
  \end{center}
  \end{figure}  

  \subsection{Proof of the classification theorem}
  
  Theorem \ref{thm:decomposition} is now an easy consequence of the results
  proved so far.
  
  \begin{proof}[Proof of Theorem \ref{thm:decomposition}]
  
    Let $S_0$, $\Sper$ and $\Sinf$ be given by Proposition
    \ref{prop:decomposition}. The hyperbolic structure on $S$ is given by
    gluing the hyperbolic structures provided by Propositions
    \ref{prop:propertiesSinf}, \ref{prop:propertiesSper} and
    \ref{prop:propertiesS0}. Since there are no strips (Corollary
    \ref{cor:no-strips}), the resulting hyperbolic structure contains no
    funnels or half-planes. Moreover the same results guarantee that for
    every component of the decomposition the first return map of $f$ exists
    and is isotopic to a periodic isometry or an isometric translation, as
    required. \qedhere

  \end{proof}
  
  We can also deduce Theorem \ref{thm:mainintro}:
  
  \begin{proof}[Proof of Theorem \ref{thm:mainintro}]

    Let \(\lambda=\partial \Sinf\). Since \(\Sinf\) is an \(f\)-invariant
    subsurface and has no compact boundary components (Proposition
    \ref{prop:propertiesSinf}), \(\lambda\) is a one-manifold without
    compact components. Let \(X\) be a complementary component of
    \(\lambda\). If \(X\) is a component of \(\Sinf\), \(\Sper\) or
    \(S_0\), we know by Theorem \ref{thm:decomposition} that the first
    return map exists and is either a translation or a finite-order map.
  
    Suppose then that \(X\) is the union of multiple components of the
    canonical decomposition, which are in \(\Sper\) or \(S_0\) by the
    definition of \(\lambda\). Note that there needs to be at least one
    component of \(\Sper\), since by construction no two components of
    \(S_0\) can be adjacent. In particular, a first return map exists.
  
    If \(X\) contains a unique component \(Y\) of \(\Sper\), then \(X\) is
    of the form \(Y\cup\bigcup_{i\in I}Z_i\), where the \(Z_i\) are
    components of \(S_0\), each of which is adjacent to \(Y\). Consider
    \(n>0\) so that \(f^n\) is isotopic to the identity on \(Y\). Then
    \(f^n\) also fixes at least one boundary component of each \(Z_i\). By
    the classification of the topological possibilities of the components
    of \(S_0\) (Proposition \ref{prop:propertiesS0}), together with the
    fact that \(f\) is (well) tempered and thus there is no twisting, we
    get that \(f^n\) is also isotopic to the identity on every \(Z_i\). In
    particular the first return map to \(X\) is periodic.
  
    If \(X\) contains multiple components of \(\Sper\), they all must be
    annuli, otherwise we could construct a periodic curve not contained in
    \(\Sper\). By Proposition \ref{prop:propertiesS0}, \(X\) contains then
    a pair of pants component \(Y\) of \(S_0\), the annular components of
    \(\Sper\) adjacent to it, together with crown components of \(S_0\),
    one per annular component of \(\Sper\) in \(X\). So topologically \(X\)
    is a pair of pants with at least one point removed from each boundary
    component. The three essential curves of \(X\) are periodic, so we can
    choose a power \(n\) of \(f\) under which they are all sent to
    themselves. As there is no twisting along these curves, \(f^n\) is the
    identity on \(X\), and again the first return map to \(X\) is periodic.
    \end{proof}

  \subsection{Constructing well-tempered maps}

  Theorem \ref{thm:mainintro} shows that every well-tempered map is
  obtained by gluing together translations and periodic maps supported on
  disjoint subsurfaces, not accumulating anywhere. If we want to construct
  well-tempered maps by following this procedure, we need to be careful
  when two supporting subsurfaces share a compact boundary component: we
  have to make sure that the gluing does not create a Dehn twist.

  On the other hand, there is no need for special care when gluing along
  lines. More precisely: 

  \begin{lemma}
    Let $\{\Sigma_i\}_{i\in I}$ be a decomposition of a surface $S$ into
    essential subsurfaces without self-adjacency, not accumulating anywhere
    and without compact boundary components. For every $i$, let $f_i$ be
    either a periodic map or a translation of $\Sigma_i$, and assume that
    for every line $\ell\subset\partial\Sigma_i\cap\partial\Sigma_j$,
    $f_i(\ell)=f_j(\ell)$. Then there is a well-tempered mapping class $f$ on
    $S$ such that the restriction of $f$ on each $\Sigma_i$ is properly
    homotopic to $f_i$.
  \end{lemma}

  \begin{proof}
    Define $S'$ to be the surface, homeomorphic to $S$, given by inserting
    strips at each line in the boundary of the decomposition. More
    precisely, we take the disjoint union of the $\Sigma_i$ and of a strip
    $S_\ell=\ell\times[0,1]$ for every $\ell$ in the boundary, and if
    $\ell\subset\partial\Sigma_i\cap\partial\Sigma_j$, we glue
    $\ell\times\{0\}$ with the copy of $\ell$ in $\Sigma_i$ and
    $\ell\times\{1\}$ with the copy of $\ell$ in $\Sigma_j$. Define then
    $f_\ell$ on $S_\ell$ by linearly interpolating $f_i$ and $f_j$:
    $$f_\ell(p,t)=(1-t)f_i(p)+t f_j(p)$$ where we have fixed an
    identification of $\ell$ with $\R$. Then we obtain a homeomorphism $f$
    of $S'$ by gluing the $f_i$ and the $f_\ell$. We claim that $f$ is
    well tempered.

    Indeed, let $\alpha$ and $\beta$ be two curves in $S'$. By the
    properties of the decomposition, and up to modifying the curves by a
    homotopy, there are only finitely many indices $j\in I$ and lines
    $\ell$ in the boundary of the decomposition such that
    $\alpha\cap\Sigma_j$, $\beta\cap\Sigma_j$, $\alpha\cap S_\ell$ or
    $\beta\cap S_\ell$ are nonempty. Since each $f_j$ is well tempered, the
    quantities $$i(f_j^n(\alpha\cap\Sigma_j),\beta\cap\Sigma_j)$$ are
    uniformly bounded. Moreover, $$i(f_\ell^n(\alpha\cap S_\ell), \beta\cap
    S_\ell)$$ is bounded by the number of intersections of $\alpha$ and
    $\beta$ with the boundary of the decomposition, hence it is also
    uniformly bounded. Therefore $i(f^n(\alpha),\beta)$ is uniformly
    bounded, so $f$ is tempered. A similar argument shows the finiteness of
    the limit set of $\alpha$. Hence $f$ is well tempered. \qedhere

  \end{proof}

  Finally, note that if we glue periodic maps and translations and we
  obtain a well-tempered map, the canonical decomposition might not
  coincide with the collection of supporting subsurfaces of the maps we are
  gluing. For instance, let $S=\R^2\ssm \Z^2$, and \[ \Sigma_1=S\cap
  \left\{ y\geq\frac{3}{4} \right\} \text{ and } \Sigma_2=S\cap \left\{
    y\leq\frac{1}{4} \right\}.\] Let $f_1$ be the map $(x,y)\mapsto
  (x+1,y)$ restricted to $\Sigma_1$, and let $f_2$ be the map $(x,y)
  \mapsto (x+2,y)$ restricted to $\Sigma_2$. We interpolate the maps $f_1$
  and $f_2$ in the intermediate strip $\{ 1/4 \le y \le 3/4\}$. Then the
  induced map on $S$ is a translation; in particular, $S=\Sinf$.


  \bibliographystyle{plain}
  \bibliography{bibliography}

\end{document}